%% file: fan.tex
\RequirePackage{fix-cm} 
\documentclass[a4paper,11pt,final]{amsart}

\input{header-fan.tex}

\title{The heart fan of an abelian category}
\author{Nathan Broomhead, David Pauksztello, David Ploog, Jon Woolf}
\dedicatory{Dedicated to the memory of Helmut Lenzing}

\renewcommand{\bm}[1]{#1}

\makeatletter
\@namedef{subjclassname@2020}{%
  \textup{2020} Mathematics Subject Classification}
\makeatother

\begin{document}

\begin{abstract}
We apply convex geometry (cones, fans) to homological input (abelian categories, hearts of bounded t-structures) to construct a new invariant of an abelian category, its heart fan.
This can be viewed as a `universal phase diagram' for Bridgeland stability conditions with the given heart.
When the abelian category is the module category of a finite-dimensional algebra, the heart fan is complete and contains the $g$-fan as the subfan of full-dimensional cones. 
The heart fan is also closely related to the wall-and-chamber structure for King semistability.
\end{abstract}

\keywords{Abelian category, triangulated category, heart, convex cone, fan, g-fan}

\subjclass[2020]{18G80, 16E35, 14F08, 18E40, 52A20}


\maketitle

\begin{center}
\begin{minipage}{0.7\textwidth}
\small
\tableofcontents
\end{minipage}
\end{center}

\addtocontents{toc}{\protect\setcounter{tocdepth}{0}}
\section*{Introduction}
\addtocontents{toc}{\protect\setcounter{tocdepth}{1}}

\noindent
Much of homological algebra is phrased in terms of abelian and triangulated categories, and these are widely used in algebraic geometry, the representation theory of finite-dimensional algebras, symplectic geometry and constructible topology. We apply convex geometry to study abelian and triangulated categories. In the first three sections we explain how to naturally construct a convex cone $C(\HH)$ --- the \defn{heart cone} --- from an abelian category $\HH$ and a fan $\fan(\HH)$ --- the \defn{heart fan} --- from the heart $\HH$ of a bounded t-structure (henceforth, a `bounded heart' or often just `heart') in a triangulated category $\DD$.

The following result defines the heart cones and the heart fan they make up in the simplest situation in which the abelian and triangulated categories have Grothendieck groups that are free of finite rank; later, we weaken this requirement by fixing a homomorphism from the Grothendieck group to a finite rank free abelian group $\LL$ . 
There is a standard partial order on bounded hearts in which $\HH[1] \geq \HH$ given by inclusion of the co-aisles of the associated t-structures.

\begin{introtheorem}[Theorems~\ref{thm:heart fan} and \ref{thm:finite heart fan}]
\label{introthm:heart fan}
Let $\HH$ be a bounded heart in a triangulated category $\DD$. Assume that the Grothendieck group $\LL \coloneqq K(\HH)$ is free of finite rank and let $\VSdual \coloneqq \Hom{\LL}{\R}$.
  \begin{enumerate}
  \item The set $C(\HH) = \{ v\in\VSdual \mid v([h]) \geq0 ~\forall h\in\HH \}$ is a closed, strictly convex cone in $\VSdual$.
  \item There is a fan $\fan(\HH)$ in $\VSdual$ generated by $\{ C(\KK) \mid \HH[1] \geq \KK \geq \HH\}$.
  \item The fan $\fan(\HH)$ does not depend on the ambient triangulated category $\DD$.
  \item If $\HH$ is a length category then the fan $\fan(\HH)$ is complete and the following are equivalent:
  \begin{enumerate}[label=(\roman*)]
  \item $\fan(\HH)$ is finite.
  \item Every heart $\KK$ with $\HH[1] \geq \KK \geq \HH$ is length.
  \item There are finitely many hearts $\KK$ with $\HH[1] \geq \KK \geq \HH$.
  \end{enumerate}
    \end{enumerate}
\end{introtheorem}

We illustrate this theorem with two examples: finite-dimensional modules over the Kronecker algebra, $\HH_1 = \mod{\kk(\Kronecker)}$ and coherent sheaves on the projective line, $\HH_2 = \coh(\PP^1)$. Their Grothendieck groups are free of rank two, so all cones live in $\VSdual \cong \R^2$. The heart cone $C(\HH_1)$ is spanned by a basis whereas $C(\HH_2)$ is a ray. See Figure~\ref{fig:fans} for further examples. 

\begin{figure}[h]
\begin{center}
\newcommand{\raycolor}{violet!30!red!70} 
\begin{adjustbox}{valign=t}
  \input{intro_Kronecker}
\end{adjustbox}
\hfill
\begin{adjustbox}{valign=t,raise=-0.5ex}
  \input{intro_projective_line}
\end{adjustbox}
\label{fig:heart fans introduction}
\end{center}
\caption{%
  Heart fans of $\HH_1 = \mod{\kk(\Kronecker)}$ and of $\HH_2 = \coh(\PP^1)$.
  The heart $\HH_1$ is length and its heart fan is complete even though it has infinitely many cones.
  The heart fan of the noetherian but non-artinian category $\HH_2$ is supported on a half-plane.}
\label{fig:intro}
\end{figure}

While all heart fans in our examples have convex support, we do not know if this always holds.
Extreme counterexamples could arise from abelian categories of Grothendieck rank $\geq 2$  with only the trivial torsion pairs, \ie having no hearts between $\HH$ and $\HH[1]$. 

\begin{description}
\item[Convex support question] Is the support of every heart fan convex?
\item[Lonely hearts question] Are there abelian categories of Grothendieck rank $\geq2$ having no non-trivial torsion pairs?
\end{description}

While the definition of the heart fan in Theorem~\ref{introthm:heart fan} is correct, we in fact start out integrally, \ie with monoids in a lattice. Moreover, we define a preliminary convex-geometric object, cofans, which can be seen as fans before dualisation.
We construct the heart fan in a three-step process: we first define the \defn{integral heart cofan} $\cofanZ(\HH)$ which induces the \defn{real heart cofan} $\cofanR(\HH)$ and finally by dualisation the (real) heart fan $\fan(\HH)$.

Cofans occur implicitly in toric geometry where, however, there is no need to stress them: cofan and associated fan determine each other if all cones are polyhedral. This breaks already in our most basic examples, leading us to develop cofans. Readers only interested in the heart fan can make a beeline for it by jumping to Theorem-Definition~\ref{thm:heart fan} where we give a self-contained proof that it is indeed a fan, possibly after looking at our \hyperlink{sub:notation}{notation} and Definition~\ref{def:effective and heart cones}.

We propose cofans as the natural convex-geometric object associated to $c$-vectors appearing in cluster theory and representation theory of finite-dimensional algebras.
In representation theory, $c$-vectors correspond to the dimension vectors of two-term simple-minded collections and the heart cofan dualises to (a completion of) the $g$-fan; see Section \ref{sec:g-fans}. However, our construction does not require that the corresponding silting theory exists, giving us a notion of a `virtual $g$-fan' when it does not. \label{todo:virtual1}
In particular, because our construction applies to any bounded heart, it applies more generally than to our main examples, which come from quivers with relations and coherent sheaves on smooth projective varieties.
By analogy with toric geometry, the $c$-vector side is algebraic and the $g$-vector side geometric: the homological algebra is on the $c$-side, the geometric objects such as fans and stability spaces are on the $g$-side.

\subsection*{Motivation and relation to other constructions}

Our first observation was that the dual cones of nearby tilted hearts fit neatly into a fan, akin to how fans of toric varieties arise.
From that initial insight we set out to explore the connection between homological algebra and convex geometry systematically. Our work is also motivated by connections to current research which are explained in Sections \ref{sec:stability conditions}, \ref{sec:wall-and-chambers} and \ref{sec:g-fans}. These can be read independently of each other.

Firstly, the heart fan $\fan(\HH)$ is closely connected with Bridgeland stability, and can be viewed  as a \defn{universal phase diagram} for stability conditions with heart $\HH$: the charge $Z$ of a stability condition $(P,Z)$ with heart $P(0,1]=\HH$ determines a real two-dimensional slice through $\fan(\HH)$ which is a planar fan whose rays record the phases of semistable objects. This relationship is apparent in the rank two examples in Figure~\ref{fig:fans}. More precisely, there is a well-defined map from rays in the support of the heart fan to wide subcategories in $\HH$ whose image on a ray in the slice cut out by $Z$ is the corresponding subcategory $P(\phi)$ of semistable objects. The details are explained in Section~\ref{sec:stability conditions}. This is part of a larger project; see below.

Now let us specialise to a category $\HH=\mod{A}$ of finite-dimensional modules over a finite-dimensional algebra $A$. In the literature, there are three constructions of convex-geometric invariants of $A$, labelled by different algebraic data in each case:
\begin{itemize}[left = 0em]
\item The \defn{wall-and-chamber structure} defined by King-semistability of modules; see Section~\ref{sec:wall-and-chambers} and \cite{BST}, \cite{King}.
\item The \defn{$\bm{g}$-fan} constructed via two-term silting complexes, or equivalently, support $\tau$-tilting pairs in the module category; see Section~\ref{sec:g-fans} and \cite{AHIKM}, \cite{Asai21}, \cite{Demonet-Iyama-Jasso}.
\item The \defn{Hall algebra scattering diagram} defined in \cite[Thm.~6.5]{bridgeland scattering}.
\end{itemize}
In the first two cases the convex-geometric object is a subset of the heart fan of $\mod{A}$, as defined in Theorem~\ref{introthm:heart fan}, but with a different labelling by algebraic data. The following statement outlines these relationships, which are explained in full detail in Sections~\ref{sec:wall-and-chambers} and \ref{sec:g-fans} respectively.

\begin{introtheorem}[Theorems~\ref{thm: wall and chamber} and \ref{thm:g-fan}]
  \label{introthm:applications}
  Let $A$ be a finite-dimensional algebra.
  \begin{enumerate}
  \item Each $v\in \VSdual$ determines a full subcategory of $v$-semistable modules in $\mod{A}$ which is constant on the relative interior of each cone of the heart fan; the stability space $\cD(M) =\{ v\in \VSdual \mid M\ \text{is $v$-semistable}\}$ of $M\in \mod{A}$ is the support of a subfan. The wall-and-chamber structure has walls the codimension one stability spaces and chambers the relative interiors of the full-dimensional cones in the heart fan.
      \item The $g$-fan of $A$ is the subfan of the heart fan of $\mod{A}$ consisting of all full-dimensional cones together with their faces. The cones are labelled by two-term presilting complexes in $\KK^b(\proj{A})$. 
      \end{enumerate}
\end{introtheorem}

\noindent
As an example, the $g$-fan of the Kronecker algebra is the heart fan shown above except for the highlighted ray which is a maximal cone of dimension one. The wall-and-chamber structure consists of the rays in the heart fan, but with a different labelling by algebraic data.
See Example~\ref{ex:3-Kronecker fans} for the heart, stability and $g$-fans of the 3-Kronecker quiver.

Section~\ref{sec:wall-and-chambers} applies to a length abelian category with finite many simple objects, extending the established theory beyond the case of finite-dimensional algebras. Similarly, we define the \defn{virtual $\bm{g}$-fan} of an abelian category $\HH$ with free Grothendieck group of finite rank as $\fan(\HH)\full$, the subfan of the heart fan generated by full cones. Theorem~\ref{introthm:applications}(2) justifies this definition, as it returns the classical $g$-fan in the special case $\HH = \mod{A}$.

The precise relationship of the heart fan of $\HH=\mod{A}$ to the Hall algebra scattering diagram of $A$ is less clear to us, because the latter is constructed by a limiting process which is hard to relate to our direct construction of $\fan(\HH)$. However, Reading \cite[Theorem 3.1]{Reading} shows that a consistent scattering diagram determines a complete fan, which in simple examples is the heart fan. It would be interesting to know in what generality this holds, and indeed whether the full scattering diagram (including the wall-crossing data) can be constructed from $\fan(\HH)$ together with the map from its support to wide subcategories of $\HH$.

\subsection*{Further work}

In a sequel to this article, we construct a `multifan' $\fan(\DD)$ for a triangulated category $\DD$. Intuitively a multifan should be thought of as a `fan over a vector space' rather than a `fan in a vector space'; multifans formalise the gluing of compatible fans. To briefly explain, consider again the examples in Figure~\ref{fig:intro}. As is well-known, $\Db(\HH_1) \cong \Db(\HH_2) \eqqcolon \DD$. In particular, the ray corresponding to $\coh(\PP^1)$ occurs in the heart fan of $\HH_1$, and is the highlighted diagonal ray emanating from the origin to the top left. In fact, the violet sections in each heart fan are identified under a basis change in $K(\DD)$. The multifan is obtained by identifying common heart cones in the disjoint union of the heart fans of all bounded hearts of $\DD$. 

In addition, we define a tangent multifan $T\fan(\DD)$ by gluing together the tangent spaces of each heart cone. The tangent multifan describes the local combinatorics and geometry of $\fan(\DD)$. Its realisation as a space can be interpreted as the `space of lax stability functions' on all hearts in $\DD$; here by `lax' we mean that the stability function may map non-zero objects of the heart to zero. The Bridgeland stability space $\stab{\DD}$ embeds as the open subset consisting of the stability functions with support and Harder--Narasimhan properties.

We want to draw attention to Parth Shimpi's upcoming \cite{Shimpi} on contraction algebras of 3-fold flops. Part of his work is the computation of heart fans in Grothendieck rank 3.

\subsubsection*{Acknowledgments}
It is a great pleasure to thank
Klaus Altmann,
Lukas Bonfert, 
Lutz Hille,
Andreas Hochenegger,
Arunima Ray,
Parth Shimpi
and
Greg Stevenson
for their comments.
This project has been supported by EPRSC grant no.\ EP/V050524/1 of the second author.
We dedicate this work to the memory of Helmut Lenzing who sadly passed away while this manuscript was being prepared.

\subsection*{Notation and terminology.}
\hypertarget{sub:notation}{}
Fix a lattice, \ie a free abelian group of finite rank, $\LL$. Denote by $\LLdual = \Hom{\LL}{\Z}$ the dual lattice, by $\VS = \LL \otimes \R$ and $\VSdual = \LLdual\otimes\R = \Hom{\LL}{\R}$ the associated real vector space and its dual respectively.

Let $\HH$ be a category having a Grothendieck group $K(\HH)$, \eg $\HH$ is abelian or has a triangulated structure.
We call $\HH$ a \defn{lattice category} if $K(\HH)$ is a lattice, \ie a finite rank, free abelian group.
Throughout, we fix a homomorphism $\ll \colon K(\HH) \to \LL$.
If $\HH$ is a lattice category, $\ll$ will be taken as identity unless explicitly stated otherwise.
In many other settings, a natural choice is the numerical Grothendieck group, $\LL = K(\HH)/\ker(\chi)$, obtained by factoring out the kernel of the Euler pairing $\chi$. This makes sense when $\chi$ is well defined and its left and right kernels agree; see \cite{Tabuada} for details.
The datum $\ll\colon K(\HH) \to \LL$ is inspired by the theory of stability conditions.

By a \defn{category over $\LL$}, we mean $\HH$ together with a homomorphism $\ll\colon K(\HH) \to \LL$. This is occasionally denoted $\HH/\LL$, trusting that no confusion with quotient categories is possible.

Elements of $\VSdual$ are denoted $v, w$ or similar.
We write $\ll(h)$ for $\ll([h])$ where $h\in \HH$ and $[h]$ is its class in $K(\HH)$, and often suppress $\LL$ and especially $\ll$ from notation. In particular, we write $v(h) \coloneqq v(\ll([h]))$ for $v\in \VSdual$ and $h\in\HH$.

In the examples, varieties and algebras are over a field $\kk$ which we assume for convenience to be algebraically closed although this assumption can be weakened in places. By a `module' we always mean a finite-dimensional left module.

\section{Convex geometry: cones, faces, fans, cofaces and cofans}
\label{sec:convex}

\noindent
Here we introduce the convex-geometric notions required later on; this section contains no homological algebra. We work integrally, \ie with cones in a finite rank lattice $\LL$. In classical language, this means that our cones are rational. However, extra care is needed because the cones need not be finitely generated.
We first introduce terminology for faces and fans; this is in parallel with the classical notions, \eg from toric variaties or $g$-fans. Starting in Subsection~\ref{sub:cofaces}, we introduce the notions of \hyperlink{sub:cofaces}{coface} and \hyperlink{sub:cofans}{cofan} which are new.

\subsection{Cones}
An \defn{integral cone} in $\LL$ is a non-empty subset $\sigma \subseteq \LL$ which is closed under sums, \ie $\sigma+\sigma \subseteq \sigma$.
In other language, a cone is a submonoid of $\LL$ but we prefer the more geometric term `integral cone'; see Remark~\ref{rem:real cones}. When the ambient lattice is clear from the the context, we just say `cone'.
We write $\cones{\LL}$ for the set of cones in $\LL$.

Minkowski sum and difference $\sigma\pm\sigma'$ of cones $\sigma, \sigma' \subseteq \LL$ are again cones.
The $\Z$-linear hull of $\sigma$ is the subgroup $\sigma-\sigma$ and the \defn{dimension} $\dim(\sigma) \coloneqq \rk(\sigma-\sigma)$ is its rank. A cone is \defn{full} if $\dim(\sigma) = \rk(\LL)$. 

A cone $\sigma$ is \defn{strictly convex} if it contains no non-zero subgroup, \ie $\sigma \cap (-\sigma) = \{0\}$.
The (integral) cone \defn{generated} by a subset $S\subseteq \LL$ is the set	
\[
E(S) \coloneqq \bigg\{ \sum_{i=0}^m n_is_i \mid n_i\in \N, s_i \in S, m\in \N \bigg\}
\]
of non-negative integral combinations of elements of $S$. A cone is \defn{polyhedral} if it is generated by a finite subset and \defn{smooth} if it is generated by a subset of a basis of $\LL$.

Given a subset $S \subseteq \LL$, its \defn{dual cone} is $\cdual{S} \coloneqq \{v \in \LLdual \mid v|_S \geq 0\}$; this is indeed a cone in the dual lattice $\LLdual = \Homm{\Z}{\LL}{\Z}$. Our unusual notation is motivated by the adjointness properties of Subsection~\ref{sub:left adjoint}. The \defn{orthogonal} of $S$ is the subgroup $S\orth=\{ v\in \LLdual \mid v|_S=0\} \subseteq \LLdual$. We set $C(S)_\Z \coloneqq \cdual{E(S)} \subseteq \LLdual$ to be the dual of the integral cone generated by $S$. 

\begin{remark}
\label{rem:real cones}
In convex geometry, a cone is a subset $\sigma$ of an $\R$-vector space that is closed under sums and under non-negative scalars, \ie $\sigma+\sigma\subseteq\sigma$ and $\Rpos\cdot\sigma \subseteq\sigma$. For distinction with integral cones, we call these \defn{real cones} in this article.
If a lattice spanning the vector space is specified then a real cone is called \defn{rational} if it is generated by a subset of the lattice.

We denote by $\Cones{\VS}$ the set of real cones in $\VS = \LL\otimes\R$. The convex hull of an integral cone $\sigma \subseteq \LL$ is the real cone $\sigma_\R \coloneqq \conv(\sigma) \subseteq \VS$. Real cones of this form are precisely the rational cones in $\VS$. All statements in this section have obvious analogues for real cones; see Subsection~\ref{sub:real cones} for more details.
The real cone generated by a subset $S$ is denoted $E(S)_\R$.
\end{remark}

\subsection{Faces}
\label{sub:faces}

A \defn{face} $\tau$ of a cone $\sigma$ is a subcone closed under summands, \ie a non-empty subset $\tau\subseteq\sigma$ such that if $a,b\in\sigma$ then $a+b\in \tau \iff a,b \in \tau$. This implies $0\in \tau$. The face relation is reflexive, anti-symmetric and transitive. Let $\fCones{\LL}$ denote the poset whose elements are the cones in $\LL$ with $\tau \faceof\sigma$ whenever $\tau$ is a face of $\sigma$. A \defn{facet} is a face of codimension 1.

For a cone $\sigma$ let $\Faces{\sigma} = \fCones{\LL}_{\leq \sigma}$ be the subposet consisting of $\sigma$ and all its faces. Any intersection of faces of a cone is again a face. The infimum of $\tau,\tau'\in \Faces{\sigma}$ is $\tau\cap\tau'$. Given a subset $\alpha \subseteq \sigma$ of a cone, the \defn{minimal face containing $\alpha$}, denoted $\mf{\alpha}$, is the intersection of all faces of $\sigma$ containing $\alpha$. The supremum of two faces $\tau,\tau'\in \Faces{\sigma}$ is the minimal face containing $\tau+\tau'$. In particular, $\Faces{\sigma}$ is an order-theoretic lattice. The following explicit description of the minimal face will be useful later.

\begin{lemma}
\label{lem:min face}
The minimal face in $\Faces{\sigma}$ containing a subcone $\alpha\subseteq\sigma$ is $\mf{\alpha} = \sigma \cap (\alpha-\sigma)$.
\end{lemma}

\begin{proof}
As $\alpha \subseteq \sigma$ are cones, so are $\alpha-\sigma$ and $\sigma\cap(\alpha-\sigma)$. To check that $\sigma\cap(\alpha-\sigma)$ is a face of $\sigma$, let $b,b'\in\sigma$ with $b+b' = a-c$ where $a\in\alpha$ and $c\in\sigma$. Then $b = a - (b'+c) \in \alpha-\sigma$ and similarly for $b'$. Hence $\sigma\cap(\alpha-\sigma) \faceof \sigma$ is a face.

Suppose $\tau$ is a face of $\sigma$ containing $\alpha$. If $b = a-c \in \sigma\cap(\alpha-\sigma)$ then $b+c = a \in \alpha \subseteq \tau$. Hence $b\in\tau$ and so $\sigma\cap(\alpha-\sigma) \subseteq \tau$. Thus it is the minimal face containing $\alpha$.
\end{proof}

\begin{remark}
\label{rmk:min face linear}
The \defn{minimal face} of a cone $\sigma$ is $\mface{\sigma} \coloneqq \mf{0} = \sigma\cap(-\sigma)$. This is the maximal subgroup contained in $\sigma$. The cone $\sigma$ is strictly convex if and only if $\mface{\sigma} = \{0\}$.
\end{remark}

\subsection{Fans}
\label{sub:fans}

\begin{definition}
\label{def:fan}
A \defn{fan} $\fan$ in $\LL$ is a set of cones in $\LL$ such that (a) $\fan$ is closed under taking faces and (b) any two cones in $\fan$ intersect in a common face. 
\end{definition}

A fan $\fan$ in $\LL$ is called
\begin{itemize}
\item \defn{polyhedral} or \defn{smooth} if all cones in $\fan$ have this property;
\item \defn{finite} if $\fan$ is a finite set of cones;
\item \defn{complete} if the union of all cones in $\fan$ is $\LL$.
\end{itemize}

\begin{example}
If $\sigma$ is a cone then $\Faces{\sigma}$ is a fan. It is polyhedral (and also finite) or smooth precisely when $\sigma$ is so. Its support is the cone $\sigma$; it is complete only when $\sigma=\LL$.
\end{example}

\begin{remark}
\label{rem:checking fan on maximal cones}
Any fan $\fan$ is generated by its maximal cones, \ie $\fan = \bigcup_{\sigma \in M} \Faces{\sigma}$ where $M\subseteq\fan$ is the subset of maximal cones with respect to the face relation. Conversely, given a set $N$ of cones in $\LL$, the union $\bigcup_{\sigma \in N} \Faces{\sigma}$ is a fan if the intersection of any two cones in $N$ is a common face. It suffices to check this for any two generating cones in $N$. 
\end{remark}

\subsection{Cofaces}
\label{sub:cofaces} \hypertarget{sub:cofaces}{}

We introduce the cofaces of a cone as localisations at a subcone: a commutative monoid with the cancellation property, \eg a cone, naturally embeds in its (Grothendieck) group of differences.

The terminology `coface' is motivated by Lemma~\ref{lem:faces and cofaces} below which shows that every coface is the localisation at a face, and yields a bijection between the faces and cofaces of a cone. 

\begin{definition}
\label{def:coface}
A \defn{coface} of $\sigma$ is the localisation $\sigma-\kappa$ of $\sigma$ at a subcone $\kappa$. We write $\sigma-\kappa\cofaceof \sigma$.
\end{definition}

While faces are subsets of a cone, cofaces are supersets: the coface relation $\sigma \cofaceof \sigma'$ implies an inclusion in the reverse order $\sigma \supseteq \sigma'$.
In Lemma~\ref{lem:coface poset}, we will show that the coface relation $\cofaceof$ partially orders the cones in $\LL$. Thus we have two posets,
\[
\fCones{\LL}\coloneqq (\cones{\LL},\faceof) \quad \text{and} \quad \cCones{\LL} \coloneqq (\cones{\LL},\cofaceof) 
\]
with the same underlying set, ordered respectively by the face and coface relations.
In analogy with the subposet of faces, we also define the subposet of cofaces of the cone $\sigma$:
\[
\Faces{\sigma} = \fCones{\LL}_{\leq\sigma} \quad \text{and} \quad \Cofaces{\sigma} \coloneqq \cCones{\LL}_{\cofaceof\sigma}
\] 

\begin{lemma}
\label{lem:coface faces}
Let $\kappa \subseteq \sigma$ be cones in $\LL$. Localisation at $\kappa$ induces a poset isomorphism 
\[
 \Faces{\sigma}_{\supseteq \kappa} \isom \Faces{\sigma-\kappa},\quad \tau \mapsto \tau-\kappa
\]
 from the poset of faces of $\sigma$ containing $\kappa$ to that of faces of $\sigma-\kappa$. The inverse is $\nu \mapsto \nu \cap \sigma$.
\end{lemma} 

\begin{proof}
Given faces $\alpha$, $\beta$ with
  $\kappa \subseteq \alpha \faceof \beta \faceof \sigma$,
we first show $\alpha-\kappa \faceof \beta-\kappa$. This implies that the map of the lemma is well defined (use $\beta=\sigma$) and monotone.
The inclusion $\alpha-\kappa \subseteq \beta-\kappa$ is obvious.
If $(b-k)+(b'-k')=a-k''$ for some $a\in \alpha$, $b,b'\in \beta$ and $k,k',k''\in \kappa$ then $b+b'+k'' = a+k+k'\in \alpha$. Hence $b,b'\in \alpha$ and $b-k, b'-k'\in \alpha-\kappa$, and thus $\alpha-\kappa \faceof \beta-\kappa$.

Now let $\alpha \faceof \sigma-\kappa$; we will show $\alpha =(\alpha \cap \sigma)-\kappa$. Clearly $\alpha\cap \sigma \faceof \sigma$ is a face. If $k \in \kappa$ and $a\in \alpha\cap \sigma$ then $k +(a-k) = a \in \alpha$ and so $k, a-k \in \alpha$. Hence $\kappa \subseteq \alpha\cap \sigma$, and also $(\alpha\cap \sigma)-\kappa \subseteq \alpha$.
For the reverse inclusion, suppose $a \in \alpha$. Then $a=b-k$ for some $b \in \sigma$ and $k \in \kappa$. Then $a+k \in \alpha \cap \sigma$, so $a = (a+k)-k \in (\alpha \cap \sigma)-\kappa$. 

Finally, we show $\alpha=(\alpha-\kappa)\cap \sigma$ for any face $\alpha$ of $\sigma$ containing $\kappa$. Clearly $\alpha \subseteq (\alpha-\kappa)\cap \sigma$. In the other direction, if $a\in (\alpha-\kappa)\cap \sigma$ then $a+k \in \alpha$ for some $k\in \kappa$, and thus $a\in \alpha$. Hence $(\alpha-\kappa)\cap \sigma \subseteq \alpha$ and we have equality.
\end{proof}

\begin{lemma}
\label{lem:coface poset}
The set $\cCones{\LL}$ of cones in $\LL$ with the coface relation $\cofaceof$ is a partial order. 
\end{lemma}

\begin{proof}
The relation $\cofaceof$ is reflexive because $\sigma = \sigma-\zerocone \cofaceof \sigma$. It is anti-symmetric because if $\sigma \cofaceof \tau$ then $\sigma\supseteq \tau$. Finally, it is transitive because if $\sigma'=(\sigma-\kappa)-\nu$ for cones $\kappa \subseteq \sigma$ and $\nu \subseteq \sigma-\kappa$ then $\sigma' = \sigma - \sigma\cap(-\sigma')$. To see this note first that $\sigma' \supseteq \sigma - (\sigma\cap(-\sigma'))$. Then for the other containment, choose $a'\in \sigma'$ and write it in the form
$
a' = (a-k)-(b-l) = (a+l) -(b+k)
$
where $a,b \in \sigma$ and $k,l\in \kappa$ with $b-l\in \nu$. Then observe that $b+k = (b-l)+(k+l)\in -\sigma'$, because $\nu,\kappa \subseteq -\sigma'$. Hence $b+k\in \sigma\cap(-\sigma')$ and so $a' = (a+l) -(b+k) \in \sigma - \sigma\cap(-\sigma')$.
\end{proof}

\begin{lemma}
\label{lem:faces and cofaces}
Let $\kappa \subseteq \sigma$ be cones in $\LL$. Then $\sigma-\kappa = \sigma-\mf{\kappa}$ is the localisation at the minimal face $\mf{\kappa}$ of $\sigma$ containing $\kappa$, and the following map is a poset isomorphism:
\[
\Faces{\sigma}\opp \isom \Cofaces{\sigma},\quad \tau \mapsto \sigma-\tau .
\]
The inverse map is $\kappa\mapsto \mface{\kappa}\cap\sigma$ where $\mface{\kappa} = \kappa\cap(-\kappa)$ is the minimal face of $\kappa$.
\end{lemma}

\begin{proof}
Clearly $\sigma-\kappa\subseteq \sigma-\mf{\kappa}$. Given $a\in \mf{\kappa}$, there is $b\in \sigma$ with $a+b\in\kappa$ by Lemma~\ref{lem:min face}. Therefore $-a = b-(a+b) \in \sigma -\kappa$. Hence $\sigma-\mf{\kappa}\subseteq \sigma-\kappa$ and we have equality. This shows that the map is surjective.

Suppose $\tau,\tau'\in \Faces{\sigma}$ with $\sigma-\tau=\sigma-\tau'$. Then for any $t\in \tau$ there are $t'\in \tau'$ and $a\in \sigma$ with $-t=a-t'$. Hence $t+a=t'\in \tau'$ and so $t\in \tau'$ because $\tau'$ is a face. Thus $\tau\subseteq \tau'$ and by symmetry $\tau=\tau'$. Thus the map is bijective.

The map is monotone because $\tau\leq \tau'$ implies $\sigma-\tau' = (\sigma-\tau)-\tau' \cofaceof \sigma-\tau$. To obtain the inverse note that if $\kappa=\sigma-\tau$ then $\mface{\kappa} = \mface{(\sigma-\tau)} = \tau-\tau$ by Lemma~\ref{lem:coface faces}. Since $\tau$ is a face of $\sigma$ we conclude $\tau = \mface{\kappa} \cap \sigma$. This map is also monotone.
\end{proof}

\begin{lemma}
\label{lem:poset height}
The posets $\fCones{\LL}$ and $\cCones{\LL}$ have height $\rk(\LL)$, \ie the longest chains have $\rk(\LL)+1$ elements.
\end{lemma}

\begin{proof}
For $\fCones{\LL}$ this follows because a proper face has strictly lower dimension. To see this suppose $\tau \leq \sigma$ and $\dim \tau = \dim \sigma$. Then $\tau-\tau$ is a finite index subgroup of $\sigma-\sigma$. Thus, given $a\in \sigma$ there are $b,b'\in \tau$ such that $na=b-b'$ for some $n\in \N$. Rearranging, $na+b'=b \in \tau$ and since $\tau$ is a face of $\sigma$ we conclude that $a\in \tau$. Therefore $\tau=\sigma$.

For $\cCones{\LL}$ it follows from the above and Lemma~\ref{lem:faces and cofaces}. (Note that the dimension of a coface is the same as that of the ambient cone.)
\end{proof}

\subsection{Cofans}
\label{sub:cofans} \hypertarget{sub:cofans}{}

These are the coface-analogues of fans: a cofan is a collection of cones that is (a) closed under cofaces and such that (b) the sum of two cones is a coface of each --- it is the equivalent of a fan with `face' replaced by `coface' and `intersection' replaced by `sum'. We rephrase this definition below using a formulation more convenient for proofs.

\begin{definition}
\label{def:cofan}
A \defn{cofan} in $\LL$ is a set $\cofan$ of cones in $\LL$ such that
\begin{enumerate}[label = (\alph*)]
\item if $\tau$ is a face of $\sigma \in \cofan$ then $\sigma-\tau \in \cofan$;
\item if $\sigma, \sigma'\in \cofan$ then $\sigma+\sigma' = \sigma-\tau$ for a face $\tau\faceof \sigma$.
\end{enumerate}
\end{definition}

\begin{example}
\label{ex:cofaces is cofan}
If $\sigma$ is a cone in $\LL$ then $\Cofaces{\sigma}$ is a cofan: it is closed under localisation at faces because the coface relation is transitive. If $\tau,\tau'\leq\sigma$ are faces then the sum
\[
(\sigma-\tau)+(\sigma-\tau') = \sigma -(\tau+\tau')=\sigma-\mf{\tau+\tau'} = (\sigma-\tau) - (\mf{\tau+\tau'}-\tau)
\]
is a coface of $\sigma-\tau$, as required for (b). Here $\mf{\tau+\tau'}-\tau$ is a face of $\sigma-\tau$ by Lemma~\ref{lem:coface faces} since $\mf{\tau+\tau'} \in \Faces{\sigma}_{\supseteq \tau}$.
\end{example}

\begin{remark}
\label{rmk:generated cofans}
Any cofan $\cofan$ is generated by its maximal cones, \ie $\cofan = \bigcup_{\sigma \in M} \Cofaces{\sigma}$ where $M\subseteq\cofan$ is the subset of cones which are maximal for the coface relation. Under inclusion of subsets, these are the minimal cones in $\cofan$. Conversely, a set $N$ of cones in $\LL$ generates a cofan $\bigcup_{\sigma \in N} \Cofaces{\sigma}$ if, and only if, the sum of any two cones in $N$ is a common coface.
\end{remark}

\begin{lemdef}
\label{lem:cofans induce fans}
If $\cofan$ is a cofan in $\LL$ then $\fan \coloneqq \bigcup_{\sigma\in \cofan} \Faces{\cdual{\sigma}}$ is a fan in $\LLdual$. We call $\fan$ the \defn{associated fan}.
\end{lemdef}

\begin{proof}
Since the face relation is transitive, $\fan$ contains all faces of its cones. To show that it is a fan it therefore suffices to show that any two cones intersect in a common face. In fact it suffices to check this for maximal cones, \ie cones which are not proper faces of other cones; see Remark~\ref{rem:checking fan on maximal cones}. Maximal cones have the form $\cdual{\sigma}$ for some $\sigma\in\cofan$. Given $\sigma,\sigma'\in\cofan$, their intersection is
$
\cdual{\sigma} \cap \cdual{\sigma'} = \cdual{\sigma+\sigma'} = \cdual{\sigma-\tau} = \cdual{\sigma} \cap \tau\orth
$
for some face $\tau \faceof \sigma$ by the cofan property. Therefore the intersection is a common face as required.
\end{proof}

We define the following attributes of cofans by analogy with fans. A cofan $\cofan$ in $\LL$ is
\begin{itemize}
\item \defn{polyhedral} if each cone in $\cofan$ is polyhedral;
\item \defn{smooth} if each cone in $\cofan$ is a coface of a smooth cone (not necessarily in $\cofan$); 
\item \defn{complete} if for each $v\in \LLdual$ there is some $\sigma\in \cofan$ with $v|_\sigma \geq 0$; 
\item \defn{finite} if $\cofan$ is a finite set of cones.
\end{itemize}

If $\cofan$ is polyhedral then the dual cone $\cdual\sigma$ of each $\sigma \in \cofan$ is polyhedral (by Gordan's Lemma), and therefore the associated fan $\fan$ is polyhedral. Likewise the associated fan of a smooth cofan is smooth. If $\cofan$ is complete then so is $\fan$.
The associated fan of a finite cofan may be infinite.
%

Fans can be pulled back, and cofans pushed forward, along homomorphisms. 

\begin{lemma}
\label{lem:cofans under linear maps}
Let $f\colon \LL \to \LL'$ be a homomorphism of lattices. Then
\begin{enumerate}
\item If $\cofan$ is a cofan in $\LL$ then $f(\cofan) = \{ f(\sigma) \mid \sigma\in\cofan \}$ is a cofan in $\LL'$.
\item If $\fan'$ is a fan in $\LL'$ then $f^{-1}(\fan') = \{ f^{-1}(\sigma') \mid \sigma'\in\fan' \}$ is a fan in $\LL$.
\item If $\fan$ is the fan associated to $\cofan$ then $(f^*)^{-1}(\fan)$ is the fan associated to $f(\cofan)$.
\end{enumerate}
\end{lemma}

\begin{proof}
The image and preimage of a cone under a homomorphism are also cones. Moreover, it is easy to check that the image of a coface is a coface of the image, and that the preimage of a face is a face of the preimage. In fact, more generally,
\[
\Cofaces{f(\sigma)} = f(\Cofaces{\sigma}) \quad \text{and} \quad \Faces{ f^{-1}(\sigma)} = f^{-1} (\Faces{\sigma}).
\]
For the first, note that if $\tau \faceof f(\sigma)$ then $f(\sigma)-\tau = f( \sigma-f^{-1}(\tau)\cap \sigma)$ is the image of a coface of $\sigma$. For the second, suppose $\tau \faceof f^{-1}(\sigma)$ is a face. Then we claim that $\tau = f^{-1}( \mf{f(\tau)} )$ is the preimage of the minimal face of $\sigma$ containing $f(\tau)$. First note that $\tau = f^{-1}(f(\tau))$ because $\tau$ contains the maximal subgroup in $f^{-1}(\sigma)$, and this contains $\ker f$. Then, using the formula for the minimal face from Lemma~\ref{lem:min face}, we have
\[
f^{-1}( \mf{f(\tau)} ) = f^{-1}( \sigma \cap(f(\tau)-\sigma) ) = f^{-1}(\sigma) \cap f^{-1}( f(\tau)-\sigma) = f^{-1}(\sigma) \cap ( \tau- f^{-1}(\sigma)).
\]
This is the minimal face of $f^{-1}(\sigma)$ containing $\tau$, \ie $\tau$ itself as claimed. 
Therefore $f(\cofan)$ is closed under taking cofaces, and $(f^{-1})(\fan')$ closed under taking faces.

To show that sums in $f(\cofan)$ are common cofaces, let $\sigma,\sigma' \in \cofan$ and $\tau\faceof\sigma$ with $\sigma+\sigma' = \sigma-\tau$. Then we have
$
f(\sigma)+f(\sigma') = f(\sigma+\sigma') = f(\sigma-\tau) = f(\sigma) -f(\tau) .
$
Thus $f(\cofan)$ is a cofan.

To show that intersections in $(f^{-1})(\Delta')$ are common faces, take $\sigma,\sigma'\in \fan'$ and compute
$
f^{-1}(\sigma)\cap f^{-1}(\sigma') = f^{-1}(\sigma\cap \sigma') = f^{-1}(\tau)
$
for some face $\tau\leq \sigma$. Thus $f^{-1}(\fan')$ is a fan.

Finally, since $\cdual{f(\sigma)} = (f^*)^{-1}(\cdual{\sigma})$ for any cone $\sigma$, the fan associated to $f(\cofan)$ is indeed
\newcommand{\bFaces}[1]{\catt{Faces}\big(#1\big)}
\[
  \bigcup_{\sigma\in \cofan} \bFaces{ \cdual{f(\sigma)} } 
= \bigcup_{\sigma\in \cofan} \bFaces{ (f^*)^{-1}(\cdual{\sigma}) } 
= (f^*)^{-1} \bigcup_{\sigma\in \cofan} \Faces{ \cdual{\sigma} } .
\qedhere
\]
\end{proof}

\subsection{Left adjoint to dualisation}
\label{sub:left adjoint}

One reason we consider two different partial orders (face and coface relations) on the set of cones in $\LL$ is that the dualisation map
\[ \rdual\colon \cCones{\LL} \to \fCones{\LLdual}, \quad \sigma \mapsto \rdual(\sigma) = \{ v \in \LLdual \mid v|_\sigma \geq 0 \} 
\]
 is monotone.
%
%
%
It also takes sums to intersections, \ie preserves limits, but (nevertheless) does not have a left adjoint. However, its restriction $\rdual \colon \cofan \to \fan$ to a map from a cofan to its associated fan is well defined (denoted by the same symbol) and does have a left adjoint
\[
 \ldual{\cofan} \colon \fan \to \cofan, \quad \tau \mapsto \min\{ \sigma \in \cofan \mid \tau \faceof \rdual(\sigma) \}
\]
which we call the \defn{lowest coface map}.
The existence of the minimum defining $\ldual{\cofan}(\tau)$ is the key part of the following result; we include the cofan in the notation because the adjoint depends upon it.
Neither of the dualisation maps is injective or surjective in general.

\begin{proposition}
\label{prop:left adjoint to dual}
Let $\cofan$ be a cofan with associated fan $\fan$. Then:
\begin{enumerate}
\item The map $\ldual{\cofan}$ is left adjoint to $\rdual$.
\item If $\cofan' \subseteq \cofan$ is a subcofan then $\ldual{\cofan'} = \ldual{\cofan}|_{\fan'}$ where $\fan'$ is the fan associated to $\cofan'$.
\item If $\tau \leq \rdual(\sigma)$ for some $\sigma \in \cofan$ then $\ldual{\cofan}(\tau)=\sigma-\sigma\cap\tau\orth$.
\end{enumerate}
\end{proposition}

\begin{proof}
Fix $\tau \in \fan$. To define $\ldual{\cofan}$, we first find $\sigma \in \fan$ such that $\rdual(\sigma)$ is minimal with $\tau \faceof \rdual(\sigma)$. We then show that among the $\sigma \in \cofan$ with the same dual cone $\rdual(\sigma)$ there is a minimal such.

By definition $\tau \faceof \rdual(\sigma)$ for some $\sigma \in \cofan$. If $\tau\faceof \rdual(\sigma')$ too then $\tau \faceof \rdual(\sigma)\cap \rdual(\sigma')= \rdual(\sigma+\sigma')$ and $\sigma+\sigma'\in \cofan$ as cofans are closed under sums. Thus we can construct a descending chain of cones in $\fan$ of which $\tau$ is a face. By Lemma~\ref{lem:poset height}, this chain terminates. Therefore we can find $\sigma\in \cofan$ such that $\rdual(\sigma)$ is minimal with $\tau\faceof\rdual(\sigma)$.

If $\sigma'\in \cofan$ has $\rdual(\sigma')=\rdual(\sigma)$ then $\sigma+\sigma' \cofaceof \sigma, \sigma'$ is a common coface, and $\rdual(\sigma+\sigma')=\rdual(\sigma)$. Thus we can construct a descending chain of cofaces of $\sigma$ with the same dual. Again by Lemma~\ref{lem:poset height}, this chain must terminate. This shows that $\ldual{\cofan}$ exists. It is immediate that $\ldual{\cofan}(\tau) \cofaceof \sigma \iff \tau \faceof \rdual(\sigma)$, so that the maps are adjoint.

The description $\ldual{\cofan}(\tau) =\min\{ \sigma \in \cofan \mid \tau \faceof \rdual(\sigma)\}$, together with the fact that cofans are closed under taking cofaces, implies $\ldual{\cofan'} = \ldual{\cofan}|_{\fan'}$ for any subcofan $\cofan'$ of $\cofan$. 

For the last part, suppose $\tau \leq \cdual{\sigma}$ and consider the subcofan $\cofan' \coloneqq \Cofaces{\sigma}$. If $\kappa$ is a face of $\sigma$ with $\tau \leq \rdual(\sigma-\kappa)= \rdual(\sigma)\cap\kappa\orth$ then $\kappa\subseteq \sigma\cap \tau\orth$. Since $\sigma\cap\tau\orth \in \Faces{\sigma}$ we get the desired formula:
$
\ldual{\cofan}(\tau)=\ldual{\Cofaces{\sigma}}(\tau) = \min\{ \sigma-\kappa \mid \tau \faceof \rdual(\sigma-\kappa)\} =\sigma-\sigma\cap\tau\orth.
$
\end{proof}

\subsection{Exposed faces, dual faces and dual face fans}
An adjunction between posets is known as a \defn{Galois connection}. Every Galois connection restricts to an order isomorphism between the image subposets, called a Galois correspondence. To describe these for the adjunction between a cofan and its associated fan, we need the following definitions. Readers mostly interested in the (co)fan of an abelian category may proceed to Section~\ref{sec:effective cones}.

We discuss dual face fans for two reasons: First, we don't have categorical descriptions for the cones in a heart fan but we do for the associated dual face fan (which is a subset of the heart fan); see Remark~\ref{rem:dual face fan}. This is used in Remark~\ref{rmk:labelled fan} (thick label map) and Section~\ref{sec:wall-and-chambers}. Second, we want to study how the standard case of polyhedral fans unfolds in greater generality.

\begin{definition}
\label{def:exposed and dual faces}
An \defn{exposed face} of a cone $\sigma$ is a subset $\sigma\cap v\orth$ for some $v\in \rdual(\sigma)$.
Let $\ExFaces{\sigma}$ be the subposet of $\Faces{\sigma}$ consisting of exposed faces.

A \defn{dual face} of $\rdual(\sigma)$ is a subset of the form $\rdual(\sigma-\kappa)$ for some face $\kappa\faceof\sigma$.
Let $\DFaces{\cdual{\sigma}}$ be the subposet of $\Faces{\cdual{\sigma}}$ of dual faces; it depends on $\sigma$, and not just on its dual $\cdual{\sigma}$.
\end{definition}

The sets of exposed and dual faces are closed under intersection. For exposed faces this is because $(\sigma \cap v\orth)\cap(\sigma\cap w\orth)=\sigma\cap(v+w)\orth$ for $v,w\in \cdual{\sigma}$. For dual faces it follows because
\[
\cdual{\sigma-\kappa} \cap \cdual{\sigma-\kappa'} = \cdual{(\sigma-\kappa) + (\sigma-\kappa')}=\cdual{\sigma - (\kappa+\kappa')}
\]
for faces $\kappa$ and $\kappa'$ of $\sigma$.

An exposed face is a face, but not necessarily conversely; see Example~\ref{ex:non-exposed face}. By the first part of the next lemma, every dual face is exposed; in Example~\ref{ex:cofaces adjunction} we will see that every dual face is dual to the localisation at an exposed face. However, $\rdual(\sigma)$ can have exposed faces which are not dual faces; see Example~\ref{ex:fan:countable_semisimple}.

\begin{lemma}
\label{lem:dual faces}
Let $\sigma\subseteq\LL$ be a cone, $\kappa, \kappa' \subseteq \sigma$ subcones and $\tau \faceof \sigma$ a face. Then:
\begin{enumerate}
\item The dual face $\rdual(\sigma-\kappa) = \rdual(\sigma)\cap k\orth$ for some $k\in \kappa$.
\item Dual faces $\rdual(\sigma-\kappa) = \rdual(\sigma-\kappa')$ agree if, and only if, the minimal exposed faces containing $\kappa$ and $\kappa'$ agree. 
\item There is a finite chain $\tau=\sigma_n \pfaceof \cdots \pfaceof \sigma_1 \pfaceof \sigma_0 =\sigma$
      with each $\sigma_i \pfaceof \sigma_{i-1}$ an exposed face.
\end{enumerate}
\end{lemma}

\begin{proof}
Parts (1) and (2) both follow from the fact that exposed faces are closed under pairwise intersection and form a finite height poset, hence are closed under arbitrary intersection since descending chains must stabilise. 
For (1), let $\kappa \subseteq \sigma$ be a subcone; then
\[
\cdual{\sigma-\kappa} 
= \cdual{\sigma} \cap \kappa\orth 
= \bigcap_{k\in \kappa} \cdual{\sigma}\cap k\orth 
= \bigcap_{i=0}^m \cdual{\sigma}\cap k_i\orth 
= \cdual{\sigma}\cap (k_0+\cdots+k_m)\orth
\]
for some $k_0,\ldots,k_m\in \kappa$. In particular, every dual face is exposed.

For (2) write $\rdual(\sigma-\kappa) = \rdual(\sigma)\cap \kappa\orth = \{ v\in \rdual(\sigma) \mid \kappa \subseteq \sigma \cap v\orth \}$. Thus $\rdual(\sigma-\kappa) = \rdual(\sigma-\kappa')$ precisely when the sets of exposed faces containing $\kappa$ and $\kappa'$ agree. Since these sets are closed under arbitrary intersection, this is equivalent to asking that the minimal exposed faces containing $\kappa$ and $\kappa'$ are the same.

For (3), if $\tau=\sigma$ we are done, so we assume $\tau$ is a proper face. It suffices to construct a proper \emph{exposed} face $\sigma' \pfaceof \sigma$ containing $\tau$, and then proceed inductively starting with $\tau \faceof \sigma'$. The process terminates because $\Faces{\sigma}$ has finite height. When $\tau \pfaceof \sigma$ is a proper face then
$
\cdual\sigma\cap \tau\orth = \cdual{\sigma-\tau} \neq \cdual{\sigma-\sigma} = \sigma\orth .
$
Choose $v \in \rdual(\sigma) \cap \tau\orth$ with $v\notin\sigma\orth$ so that 
$\sigma'\coloneqq \sigma\cap v\orth$ is a proper exposed face containing $\tau$. 
\end{proof}

It follows from the above lemma that $\DFaces{\cdual{\sigma}}$ is a subposet of $\ExFaces{\rdual(\sigma)}$. Let $\CEFaces{\sigma} \coloneqq \{ \sigma-\kappa \mid \kappa \in \ExFaces{\sigma}\}$ be the subposet of $\Cofaces{\sigma}$ consisting of the localisations at exposed faces of the cone $\sigma$. We refer to these as \defn{co-exposed faces}.

\begin{corollary}
\label{cor:Galois correspondence}
For any cofan $\cofan$, the adjunction $\ldual{\cofan} \dashv \rdual$ restricts to an isomorphism between the images of $\rdual$ and $\ldual{\cofan}$, and there is a commutative diagram of monotone maps
\[
\begin{tikzcd} 
  &[-1.7em] \cofan \ar[shift right=0.5ex,rr,swap,"\rdual"]         \ar[two heads,ddrr] &&
            \fan   \ar[shift right=0.5ex,ll,swap,"\ldual{\cofan}"] \ar[two heads,crossing over,ddll] \\
\\
    \bigcup\limits_{\mathclap{\sigma \in M}} \CEFaces{\sigma} \ar[equals,r]
  & \ldual{\cofan}(\fan)                                  \ar[<->,rr] \ar[hook,uu] &&
    \rdual(\cofan)                                        \ar[hook,uu] \ar[equals,r]
  &[-1.7em] \bigcup\limits_{\mathclap{\sigma\in M}} \DFaces{\cdual{\sigma}}.
\end{tikzcd}
\]
where $M \subseteq \cofan$ is the subset of maximal cones for the coface relation.
\end{corollary}

\begin{proof}
Any Galois connection restricts to a Galois correspondence between the images; in our case this yields the above diagram.

The equality $\rdual(\cofan) = \bigcup_{\sigma\in M} \DFaces{\rdual(\sigma)}$ follows because any cone in $\cofan$ is a coface of some maximal cone. For the other equality, note similarly that any $\tau\in \fan$ is a face of $\rdual(\sigma)$ for some $\sigma \in M$. Thus $\ldual{\cofan}(\tau)=\sigma-\sigma\cap \tau\orth \in \CEFaces{\sigma}$ by the last part of Proposition~\ref{prop:left adjoint to dual} and Lemma~\ref{lem:dual faces}(1). 
Now consider a co-exposed face $\sigma - \kappa$. By the same reasoning as above $\ldual{\cofan}\rdual(\sigma-\kappa)=\sigma-\kappa'$ for some exposed face $\kappa'$ of $\sigma$. But then $\kappa'=\kappa$ by the last part of Lemma~\ref{lem:dual faces} because
$
\rdual(\sigma-\kappa') = \rdual\ldual{\cofan}\rdual(\sigma-\kappa) = \rdual(\sigma - \kappa)
$
by the adjunction triangular identity.
Therefore every co-exposed face of $\sigma$ is in the image of $\ldual{\cofan}$.
\end{proof}

\begin{remark}
\label{rmk:dual faces and fans}
We refer to $\rdual(\cofan)$ as the \defn{dual face fan} associated to the cofan $\cofan$. It is a set of dual cones in $\LLdual$ which is
\begin{enumerate}[label = (\alph*)]
\item closed under taking dual faces and such that
\item the intersection of any two cones is a common dual face of each. 
\end{enumerate}
The first follows because if $\sigma\in \cofan$ and $\kappa \in \Faces{\sigma}$ then $\sigma-\kappa\in \cofan$; the second because if $\sigma,\sigma'\in \cofan$ then $\rdual(\sigma)\cap\rdual(\sigma')= \rdual(\sigma+\sigma') = \rdual(\sigma-\kappa) = \rdual(\sigma'-\kappa')$ for some faces $\kappa\faceof\sigma$ and $\kappa'\faceof\sigma'$. In general $\rdual(\cofan)$ has fewer cones than $\fan$ but the same support:
  $\bigcup_{\sigma\in \cofan} \rdual(\sigma) = \bigcup_{\sigma\in\fan} \sigma$. 
\end{remark}

\begin{example}[Coface cofans]
\label{ex:cofaces adjunction}
The simplest example is the cofan $\Cofaces{\sigma}$ from a single cone $\sigma$. The associated fan is $\Faces{\cdual{\sigma}}$ and $\ldual{\Cofaces{\sigma}}(\tau) = \sigma-\sigma\cap\tau\orth$. 

The corresponding dual face fan is $\DFaces{\rdual(\sigma)}$ and the Galois connection restricts to an isomorphism $\rdual\colon \CEFaces{\sigma} \isom \DFaces{\cdual{\sigma}}$. Composing with localisation is injective by Lemma~\ref{lem:faces and cofaces} and surjective by definition of co-exposed face, hence yields an isomorphism
$\ExFaces{\sigma}\opp \isom \DFaces{\cdual{\sigma}}$, $\kappa \mapsto \rdual(\sigma-\kappa)$
with inverse $\tau \mapsto \sigma\cap \tau\orth$. 

The Galois connection is a poset isomorphism, and so agrees with the Galois correspondence, precisely when all faces of $\sigma$ and of $\cdual{\sigma}$ are exposed.
\end{example}

\begin{example}[Polyhedral cofans: toric varieties and $g$-fans]
\label{ex:polyhedral and toric}
The situation is simpler for a polyhedral cofan $\cofan$, \ie a cofan all of whose cones are polyhedral. Every face of a polyhedral cone is exposed and the dual of a polyhedral cone is again polyhedral. Hence the fan $\fan$ associated to $\cofan$ is polyhedral as well. Therefore $\ldual{\cofan}(\tau) = \cdual{\tau}$ for $\tau\in\fan$, and $\rdual$ is involutive. The Galois connection and correspondence coincide, and the commuting diagram of Corollary~\ref{cor:Galois correspondence} boils down to the single poset isomorphism
 $\rdual \colon \ldual{\cofan}(\fan) = \cofan \isom \fan = \rdual(\cofan)$.

The associated fan construction is the reverse of the process used to construct a toric variety: there one starts with a finite polyhedral fan, constructs the associated cofan of monoids in the dual lattice, and glues together the spectra of the algebras generated by these monoids to obtain the toric variety. The polyhedral condition ensures that the algebras are finitely generated, and the finiteness that the variety is quasi-compact.
%
%
Thus one can construct a toric variety directly from a finite polyhedral cofan; the associated fan of the cofan is the fan of the toric variety. More generally, one can construct a toric scheme from an arbitrary cofan.

The $g$-fan of a finite-dimensional algebra is polyhedral, and even smooth. The cofan underlying the $g$-fan is explained in algebraic terms in Section~\ref{sec:g-fans} where we also show that the heart fan of the module category is a completion of the $g$-fan, which need no longer be polyhedral. 
\end{example}

\subsection{Real cofans and fans}
\label{sub:real cones}

Classical fans of rational cones are the real versions of all of the above. Let $\Cones{\VS}$ be the set of real cones in $\VS = \LL\otimes\R$, \ie subsets closed under addition and scalar multiplication by elements of $\R_{\geq 0}$; see Remark~\ref{rem:real cones}. The notions defined in this section for integral cones extend to real cones.
For example, the dual of a real cone $\sigma \subseteq \VS$ is
  $\cdual{\sigma} \coloneqq \{v\in \LL_\R\vdual \mid v|_\sigma \geq 0\}$.
This is a closed cone and, by the Hyperplane Separation Theorem, $\cdual{\cdual{\sigma}} = \overline{\sigma}$.
Given a subset $S \subseteq \VS$, the dual of the real cone $E(S)_\R$ is denoted $C(S) \coloneqq \cdual{E(S)_\R}$.
A cone in $\VS$ generated by a subset of an $\R$-basis is called \defn{simplicial}. 

Integral and real cones are compatible in a strong sense; this simple behaviour arises because $\R_{\geq 0}$ is a cone in $\R$ with only itself and $\zerocone$ as faces. The compatibility arises from a Galois connection between cones in $\LL$ and in $\VS$ ordered by inclusion:
\[
\begin{tikzcd} 
\Cones{\LL}  \ar[shift left=0.5ex,rr,"(-)_\R"]  &&
\Cones{\VS}  \ar[shift left=0.5ex,ll,"-\cap\LL"]
\end{tikzcd}
\]
where $\sigma_\R \coloneqq \conv(\sigma)$ is the convex hull of $\sigma$ in $\VS$. The unit of this adjunction is the inclusion $\sigma \subseteq \sigma_\R \cap \LL$ of an integral cone in its \defn{saturation}
  $\sat{\sigma} \coloneqq \{ a \in \LL \mid na\in\sigma \text{ for some }  n\in\N_{>0} \}$,
and the counit is the inclusion $(\tau\cap \LL)_\R\subseteq \tau$ of the \defn{rationalisation} of a real cone. The Galois connection therefore restricts to a correspondence between saturated cones (\ie those equal to their saturation) in $\LL$ and rational cones in $\VS$.
The following lemma is immediate.
\begin{lemma}
\label{lem:real-integral-faces-cofaces}
The left adjoint $\sigma \mapsto \sigma_\R$ preserves cofaces and is monotone for the coface relation; the right adjoint $\tau \mapsto \tau\cap \LL$ preserves faces and is monotone for the face relation.
\end{lemma}

\begin{lemma}
\label{lem:integral-to-real-cones}
If $\sigma$ is a cone in $\LL$ then there is a commutative diagram of poset isomorphisms
\[
\begin{tikzcd}
    \Faces{\sigma}       \ar[d,swap,"\tau\,\mapsto\,\sigma-\tau"] \ar[rrr,"\tau\,\mapsto\, \tau_{\,\R}"] 
&&& \Faces{\sigma_\R}     \ar[d,"\nu\,\mapsto\,      \sigma_\R-\nu"] \\
    \Cofaces{\sigma}\opp \ar[rrr,swap,"\sigma-\tau\,\mapsto\, (\sigma-\tau)_\R"] 
&&& \Cofaces{\sigma_\R}\opp 
\end{tikzcd}
\]
\end{lemma}

\begin{proof}
The vertical maps are poset isomorphisms by Lemma~\ref{lem:faces and cofaces}, and its analogue for real cones which has verbatim the same proof. The diagram commutes because $(\sigma-\tau)_\R = \sigma_\R-\tau_\IR$. Hence it suffices to show that $\tau \mapsto \tau_\IR$ is a well-defined poset isomorphism $\Faces{\sigma} \to \Faces{\sigma_\R}$. 

First we verify the map is well defined, \ie that $\tau_\IR$ is a face of $\sigma_\R$. If $\tau=\sigma\cap v\orth$ for some $v\in \cdual\sigma$ is an {\em exposed} face then $\tau_\IR = (\sigma \cap v\orth)_\R = \sigma_\R \cap v\orth$ is an exposed face of $\sigma_\R$. The general case follows inductively because by Lemma~\ref{lem:dual faces}(3) for any face $\tau$ there is a finite chain $\tau = \sigma_n \pfaceof \cdots \pfaceof \sigma_1 \pfaceof \sigma_0 = \sigma$ in $\Faces{\sigma}$ in which $\sigma_i$ is an exposed face of $\sigma_{i-1}$. Since the map preserves inclusions it is monotone for the face relation because a face of a cone is also a face of any subcone containing it.

We claim the inverse is $\Faces{\sigma_\R} \to \Faces{\sigma}$, $\kappa \mapsto \kappa \cap \sigma$.
%
%
It is easily seen to be monotone. Moreover, each face of a rational cone must be rational: if some positive real combination of lattice elements in a cone lies in a face, then each of those lattice elements lies in that face too. It follows that $\kappa = (\kappa \cap \sigma)_\R$. Now consider $\tau \in \Faces{\sigma}$. The composite $\tau_\IR \cap \LL$ is the saturation of $\tau$. Hence $\tau = (\tau_\IR \cap \LL)\cap \sigma = \tau_\IR \cap \sigma$. This completes the proof.
\end{proof}

The next result shows that any integral cofan has a real version whose cones are simply the convex hulls of the cones of the cofan; we call this the \defn{realification} of the integral cofan. Likewise, any real fan has an integral version, but in this case we do not give an explicit description of the cones as this will depend on whether the fan has irrational cones or not.

\begin{corollary}
\label{cor: real cofan}
  If $\cofan$ is a cofan in $\LL$ then $\cofan_\R \coloneqq \{ \sigma_\R \mid \sigma \in \cofan \}$ is a cofan in $\VS$.
  
  If $\fan$ is a fan in $\VS$ then $\fan_\Z \coloneqq \bigcup_{\tau\in \fan} \Faces{\tau \cap \LL}$ is a fan in $\LL$.
\end{corollary}

\begin{proof}
By Lemma~\ref{lem:integral-to-real-cones} the set of cones $ \cofan_\R$ is closed under cofaces. Moreover, sums of its cones are common cofaces because $\sigma_\R+\tau_\IR=(\sigma+\tau)_\R$. Therefore $\cofan_\R$ is a cofan in $\VS$.

By construction $\fan_\Z$ is a set of cones in $\LL$ closed under faces. To show it is a fan, it suffices to check that the intersection of any two maximal cones is a common face. This follows from Lemma \ref{lem:real-integral-faces-cofaces} as $(\tau\cap\LL)\cap(\tau'\cap \LL) = (\tau\cap\tau')\cap\LL$ and $\tau\cap\tau'$ is a common face of $\tau$ and $\tau'$.
\end{proof}

\begin{example}
Consider the cofan $\cofan = \{\sigma, \Z^2\}$ where $\sigma=\{(m,n) \in \Z^2 \mid n>\sqrt{2}m\}$. The associated fan $\fan = \Faces{\cdual\sigma} = \{0\}$ but the fan $\fan_\R$ associated to the realification $\cofan_\R$ is generated by the ray through $(-\sqrt{2},1)$.
This shows that one cannot, in general, construct the associated real fan $\fan_\R$ from the associated integral fan $\fan$.
\end{example}

\begin{remark}
\label{rem: real cofans under linear maps}
The real analogue of Lemma~\ref{lem:cofans under linear maps} about the functoriality of fans and cofans under a base change $f\colon \LL \to \LL'$ holds, with the same proof, for real fans and cofans and a linear map $V \to V'$. 
Moreover, this is compatible with passing to or from the real version when $f_\R \coloneqq f\otimes \R \colon \VS \to \VS'$ in the sense that 
\begin{enumerate}
\item $f(\cofan)_\R = f_\R(\cofan_\R)$ for an integral cofan $\cofan$ in $\LL$,
\item $f^{-1}(\fan'_\Z) = f_\R^{-1}(\fan')_\Z$ for a real fan $\fan'$ in $\VS'$.
\end{enumerate} 
\end{remark}

\begin{remark}[integral cofans and real fans]
\label{rem:philosophy}
Starting from an integral cofan $\cofan$, we have defined the real cofan $\cofan_\R$ and their associated fans $\fan$ and $\fan_\R$. Out of these four objects, only two appear in the subsequent sections: the integral cofan of a heart $\HH$ that for emphasis we denote $\cofanZ(\HH)$ and the associated fan of the real cofan that we denote $\fan(\HH)$. Our philosophy is that the integral cofan carries all the relevant information and is the most basic object whereas the real fan is the familiar convex-geometric object.
\end{remark}

\section{Effective cone and heart cone of an abelian category}
\label{sec:effective cones}

\noindent
Let $\HH$ be an abelian category and $\ll\colon K(\HH) \to \LL$ a homomorphism to a lattice, \ie a free abelian group of finite rank. We call this setup an `abelian category over a lattice' and write $\HH/\LL$.

\begin{remark}
\label{rem:Grothendieck kernel}
Classes $[h] \in K(\HH)$ play a prominent role in this article. The \defn{null subcategory} is $\Groker{\HH} \coloneqq \{ h\in\HH \mid 0 = [h] \in K(\HH) \} \subseteq \HH$. It can happen that $\Groker{\HH} = \HH$, \eg whenever $\HH$ allows countable coproducts, using a variant of the Eilenberg swindle. The categories we have in mind feature $\Groker{\HH}=0$ instead, but see Example~\ref{ex:cone:tilted projective line} for an exception.
Our main examples are:
\begin{itemize}[left=0em]
\item $\HH = \mod{A}$, the category of finitely generated left $A$-modules over a finite-dimensional algebra $A$. In this case, $\LL \coloneqq K(\HH)$ is a free abelian group generated by the pairwise non-isomorphic simple $A$-modules. A module $M\in\HH$ with $0 = [M] \in K(\HH)$ has vanishing dimension vector, hence $M = 0$.
The same reasoning applies to any length category $\HH$.
\item $\HH = \coh(X)$, the category of coherent sheaves on a smooth projective variety $X$.
Then $K(X)=K(\HH)$ may have infinite rank or torsion; let $\LL \coloneqq N(X)$ be the numerical Grothendieck group, \ie the quotient of $K(X)$ by the radical of the Euler pairing. This is a lattice; see \cite{Tabuada-fg} for a generalisation. A coherent sheaf $F\in\HH$ with $0=[F] \in N(X)$ has vanishing Hilbert polynomial by the Hirzebruch--Riemann--Roch theorem, forcing $F=0$. In particular, $N_\HH=0$.
%
\end{itemize}
\end{remark}

\begin{definition}
\label{def:effective and heart cones}
Let $\HH$ be an abelian category over a lattice $\LL$. We call
\[ \begin{array}{r@{\,}ll}
 E(\HH) &\coloneqq \{ a_1\ll(h_1)+\cdots+a_n\ll(h_n) \mid n\in\N, a_i\in\N, h_i\in\HH \} \subseteq \LL & \text{the \defn{effective cone} of $\HH$ and} \\
 C(\HH) &\coloneqq \{ v\in\VSdual \mid v(h) \geq 0 \text{ for all } h\in\HH \} \subseteq \VSdual & \text{the \defn{heart cone} of $\HH$.}
\end{array} \]
\end{definition}

\begin{remark}
\label{rem:effective cone}
In words, $E(\HH)$ is the \emph{integral} cone generated by the image of $\ll \colon K(\HH) \to \LL$.
The real version of the effective cone is its convex hull $E(\HH)_\R = \conv(E(\HH))$ inside $\VS$. This is a linearisation, and often drastic simplification, of $E(\HH)$ itself.
The heart cone $C(\HH)$ is the \emph{real} cone that is dual to the (real or integral) effective cone. In the notation of Section~\ref{sec:convex}, this is $C(\HH) = \cdual{E(\HH)_\R}$.
Our choice to denote the integral cones $E(\HH)$ and the real cones $C(\HH)$ follows the philosophy outlined in Remark~\ref{rem:philosophy}.

We write $v|_\HH \geq0$ to mean $v\in C(\HH)$, expanding on the notational shortcut $v(h) = v(\ll([h]))$. Hence we can rewrite the heart cone as
  $C(\HH) = \{ v\in\VSdual \mid v|_\HH \geq 0 \}$.

The terminology `effective cone' is inspired by effective divisors or cycles in algebraic geometry.
We say `heart cone' as the abelian categories will occur as hearts in triangulated categories.

If $\HH$ is a lattice category, \ie $K(\HH)$ is a lattice, its effective cone is just the Grothendieck monoid. These are subtle invariants. For example, Shunya Saito shows in \cite{Saito} that if $X$ is a noetherian scheme then the Grothendieck monoid of $\coh(X)$ determines the topological space underlying $X$.
\end{remark}

Effective cones need be neither polyhedral nor strictly convex; see Examples~\ref{ex:cone:projective line} and \ref{ex:cone:tilted projective line}.
However, if $\HH$ is a lattice category of finite length with the standard choice $\LL = K(\HH)$ then $E(\HH)$ is smooth, so in particular polyhedral and strictly convex; see Proposition~\ref{prop:full cones}.

Below, $\Groker{\HH,\ll} \coloneqq \{ h\in\HH \mid 0 = \ll(h) \in \LL \}$ is the \defn{null subcategory} of $\HH$ as a category over $\LL$.

\begin{lemma}
\label{lem:cones}
Let $\HH$ be an abelian category over $\LL$.
\begin{enumerate}
\item If $E(\HH)$ is strictly convex then the null subcategory $\Groker{\HH,\ll}$ is a Serre subcategory of $\HH$.
\item If $\ll\colon K(\HH) \to \LL$ is surjective then $E(\HH)$ is full and $C(\HH)$ is strictly convex.
\end{enumerate}
\end{lemma}

\begin{proof}
(1) Assume $E(\HH)$ is strictly convex. Let $h\in\Groker{\HH,\ll}$ and $0 \to h' \to h \to h'' \to 0$ be a short exact sequence. Then $\ll(h') + \ll(h'') = \ll(h) = 0$ and strict convexity of $E(\HH)$ implies $\ll(h') = \ll(h'') = 0$. Thus $\Groker{\HH,\ll}$ is closed under subobjects and quotients. As $\Groker{\HH,\ll}$ is always closed under extensions, it is a Serre subcategory.

(2) Assume $\ll\colon K(\HH) \to \LL$ is surjective. Then $E(\HH)$ is a full cone because the images $\ll(\HH)$ span $\LL$, and hence its integral dual cone is strictly convex and so is its convex hull $C(\HH)$.
\end{proof}

The next examples exhibit effective cones that are not strictly convex or polyhedral. Moreover, distinct hearts inside a derived category can have the same effective or heart cone.

\begin{example}[projective line]
\label{ex:cone:projective line}
Consider the category $\HH = \coh(\PP^1)$ of coherent sheaves on the projective line. Then $\LL \coloneqq K(\HH) \cong \Z^2$ and we fix $[\cO], [\cO_p]$ as an orthonormal basis, \ie $[A] = (\rk(A),\deg(A))$ for $A\in\HH$. Here $\cO_p$ denotes the length one skyscraper sheaf at a closed point $p\in\PP^1$ and $\cO = \cO_{\PP^1}$ is the structure sheaf. We use the basis to identify $\LL$ and $\LLdual$.

As there is an infinite chain of line bundles, the effective cone $E(\HH)$ is the set of lattice points in a half-plane missing a boundary ray: $\{(m,n) \in \Z^2 \mid m>0 \text{ or } m=0, n \geq 0 \}$. This cone is not polyhedral. The heart cone $C(\HH)$ is the ray $\R\times \zerocone$. See Example~\ref{ex:coherent sheaves heart cone} for a generalisation.
\end{example}

\begin{example}[curves tilted in points]
\label{ex:tilted curve}
Let $\HH = \coh(X)$ be the category of coherent sheaves on a smooth, projective curve $X$ over its numerical Grothendieck group $\LL \coloneqq N(X) \isom \Z^2$, $\ll(A) = (\rk(A),\deg(A))$, as in Remark~\ref{rem:Grothendieck kernel}.
Fix a subset $M \subseteq X$ of closed points.

Let $\TTT{M} = \clext{ \cO_p \mid p\in M }$ be the subcategory generated by skyscraper sheaves at points in $M$ and
$\FFF{M} = \clext{\Pic(X), \cO_q \mid q \notin M}$ be generated by all line bundles and the skyscraper sheaves at the remaining points.
Then $(\TTT{M},\FFF{M})$ is a torsion pair, see Appendix~\ref{appA:torsion pairs}, and gives rise to a new heart
  $\HH_M \coloneqq \FFF{M}[1] * \TTT{M}$
in $\Db(X)$, see Appendix~\ref{appA:tilting}.
Moreover, $(\TTT{M},\FFF{M})$ is a \emph{cotilting torsion pair} in $\HH$, that is each object of $\HH$ is a quotient of an object of $\FFF{M}$.
This implies that $\HH_M$ is the heart of a faithful bounded t-structure, \ie $\Db(\HH_M) \cong \Db(X)$; see \cite[Prop.~5.4.3]{Bondal-vdB}.
%
%
For later reference, we introduce the following notation and terminology:
\begin{align*}
  \reversedHeart{X} & \coloneqq \HH_X, && \text{the \defn{reversed geometric heart}}; \\
  \mixedHeart{X}    & \coloneqq \HH_M, && \text{a \defn{mixed geometric heart} where $\emptyset \neq M \subsetneq X$}.
\end{align*}
We don't specify the subset $M$ as it will never matter precisely which points are chosen.

Let $p, q\in X$ be linearly equivalent points, \ie $\cO_X(-p) \cong \cO_X(-q)$. If moreover $p\in M$, $q\notin M$ then $\cO_p, \cO_q[1] \in \mixedHeart{X}$.
The short exact sequences
$0 \to \cO_X(-p) \to \cO_X \to \cO_p \to 0$ and
$0 \to \cO_X(-q) \to \cO_X \to \cO_q \to 0$ show
$[\cO_p \oplus \cO_q[1]] = 0 \in K(\mixedHeart{X})$. Thus the null subcategory $\NN_{\mixedHeart{X},\ll}$ is non-zero and not closed under direct summands, so not a Serre subcategory.
\end{example}

\begin{example}[projective line tilted in points]
\label{ex:cone:tilted projective line}
We combine the previous two examples.
The Auslander--Reiten quiver of $\reversedHeart{\PP^1}$ looks like that of $\coh(\PP^1)$, except that the skyscraper sheaves sit before the line bundles. In $\mixedHeart{\PP^1}$, only some of the skyscraper sheaves are tilted before the line bundles. The notation is supposed to suggest that all/some skyscraper sheaves appear at the opposite end of the Auslander--Reiten quiver of the category.

The categories $\coh(\PP^1), \mixedHeart{\PP^1}$ and $\reversedHeart{\PP^1}$ have different effective cones
\begin{itemize}
\item $E(\coh{\PP^1}) = \{ (m,n) \mid m>0 \text{ or } m=0, n\geq0 \}$,
\item $E(\mixedHeart{\PP^1}) = \{ (m,n) \mid m\geq0 \}$,
\item $E(\reversedHeart{\PP^1}) = \{ (m,n) \mid m>0 \text{ or } m=0, n\leq0 \}$,
\end{itemize}
but the same heart cone $\R \times\{0\}$. 
%
%
The cone $E(\mixedHeart{\PP^1})$ is not strictly convex as it contains the subgroup $n=0$. It is the effective cone of infinitely many distinct hearts in $\Db(\PP^1)$.
The null subcategory, $\NN_{\mixedHeart{\PP^1}}$, of $\mixedHeart{\PP^1}$ mentioned in Remark~\ref{rem:Grothendieck kernel} is non-zero.
\end{example}

\begin{example}[coherent sheaves]
\label{ex:coherent sheaves heart cone}
Let $X$ be a connected, smooth and projective variety of dimension $d$ and $\HH = \coh(X)$ the category of coherent sheaves on $X$.
As in Remark~\ref{rem:Grothendieck kernel}, $\LL = N(X)$ is the numerical Grothendieck group.
%
Let $D \subset X$ be an effective ample divisor and $Z\subseteq X$ a closed subset of positive dimension; we can assume that $D$ intersects $Z$ transversally. The classes of the sheaves $\cO_Z(nD)$ for $n\in\Z$ give an entire coset $[\cO_Z]+\Z[D]$ in the effective cone $E(\HH)$. In contrast, if $p\in X$ is a point then only the monoid $\N[\cO_p]$ lies in $E(\HH)$. It follows that the heart cone $C(\HH) \cong \Rpos\cdot \rk$ is one-dimensional, generated by the rank function.
\end{example}

This example shows that the heart cones of smooth projective varieties are not interesting in themselves, as they are independent of the geometry. Example~\ref{ex:heart cones of algebraic categories} below shows that the heart cones of finite-dimensional algebras are always orthants, likewise regardless of the algebra. These cones become interesting when studied as part of the heart fans in Subsection~\ref{sub:heart fan}.

If $\LL = K(\HH)$, full cones correspond to an algebraic setting. An abelian category is \defn{algebraic} if it is noetherian and artinian, and has finitely many simple objects up to isomorphism.

\begin{proposition}
\label{prop:full cones}
Let $\HH$ be an abelian lattice category. Then the following are equivalent:
\begin{enumerate}[label = (\roman*)]
\item The full subcategory $\Groker{\HH} = \{ h\in \HH \mid [h]=0\}$ is Serre and $\HH/\Groker{\HH}$ is algebraic.
\item The effective cone $E(\HH)$ is smooth and full.
\item The heart cone $C(\HH)$ is full. 
\end{enumerate}
\end{proposition}

\begin{proof}
$(i) \timplies (ii)$:
Suppose $\Groker{\HH}$ is a Serre subcategory and $\HH/\Groker{\HH}$ is algebraic. The quotient functor induces a canonical isomorphism $K(\HH) = K(\HH/\Groker{\HH})$ identifying the effective cones of $\HH$ and $\HH/\Groker{\HH}$. Since $\HH/\Groker{\HH}$ is algebraic, the classes of its simple objects form a basis of $K(\HH/\Groker{\HH})=K(\HH)$. The effective cone $E(\HH/\Groker{\HH}) = E(\HH)$ is generated by this basis, so is smooth and full. 

$(ii) \timplies (iii)$:
By definition, a smooth, full cone is generated by a basis. Its dual cone is generated by the dual basis, hence is also smooth and full.

$(iii) \timplies (i)$:
If the heart cone $C(\HH)$ is full then its dual cone $\rdual(C(\HH))$ is strictly convex, hence so is the real subcone $E(\HH)_\R \subseteq \rdual\rdual(E(\HH))$ and then the integral cone $E(\HH)$ is too. By Lemma~\ref{lem:cones} $\Groker{\HH} = \Groker{\HH,\id}$ is a Serre subcategory.

Now fix $h\in \HH$. Any subobject of $h$ has class in $E(\HH) \cap ([h]-E(\HH))$. Since $E(\HH)$ is strictly convex, there are only finite many lattice points in this region, hence only finitely many possible classes for subobjects. Therefore the class of any ascending or descending chain of subobjects is eventually constant, and the quotients of successive terms in the chain are eventually in $\Groker{\HH}$. We conclude that $\HH/\Groker{\HH}$ is a length abelian category. As $K(\HH) = K(\HH/\Groker{\HH})$ is a lattice, the classes of simple objects form a basis and the category $\HH/\Groker{\HH}$ is even algebraic.
\end{proof}

\begin{example}[effective and heart cones of algebraic categories are smooth]
\label{ex:heart cones of algebraic categories}
Let $\HH$ be an algebraic abelian category, such as $\HH = \mod{A}$ for a finite-dimensional algebra $A$. Then $K(\HH)$ is a free abelian group of finite rank $r$. We put $\LL = K(\HH)$, as always in this situation unless explicitly stated otherwise. By the proposition, $E(\HH)$ and $C(\HH)$ are smooth. Concretely, the heart cone of an algebraic category only depends on the rank $r$ and is isomorphic to the orthant $\Rpos^r \subset \R^r \cong \VSdual$ generated by the classes of the simple objects.
\end{example}

\begin{example}[a round cone]
\label{ex:round cone}
Let $\HH$ be an abelian category generated by countably many isomorphism classes of simple objects such that $K(\HH) \cong \Z^\infty$ is free of countably infinite rank, \eg $\HH = \mod{\kk^\Z}$.
%
%
Set $\LL \coloneqq \Z^3$ and $S \coloneqq \{ (p,q,r) \in \LL \mid p^2+q^2\leq r^2 \ \text{and}\ r\geq 0\}$. Choose a surjective homomorphism $\ll \colon K(\HH) = \Z^\infty \to \LL$ mapping the classes of the simple objects onto the countable set $S$. Then $E(\HH)=E(S)$ is a cone with countably many exposed faces, one for each primitive Pythagorean triple $(p,q,r)$. Its convex hull in $\VS = \R^3$ is a non-closed cone whose closure is the the round cone $\{ x^2+y^2\leq z^2, z\geq 0\}$ which has uncountably many exposed faces. The heart cone $C(\HH) = \cdual{E(\HH)_\R}$ is a closed round cone.
\end{example}

\medskip
\noindent
\begin{minipage}[b]{0.84\textwidth}
\begin{example}[a non-exposed face]
\label{ex:non-exposed face}
Extending Example~\ref{ex:round cone}, we map the classes of the countably many simple objects into $\LL = \Z^3$ so that the intersection in $\VSdual$ of the heart cone $C(\HH)$ with the plane $z=1$ is as shown on the right, where $A$ and $B$ are in the dual lattice. Then the ray through the point $A$ intersects $C(\HH)$ in a non-exposed face, \ie not cut out by a plane section. 
\end{example}
\end{minipage}
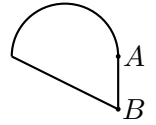
\hfill
\begin{tikzpicture}[scale=0.7]
  \draw[thick] (-1,0) -- (1,-1) -- (1,0);
  \draw[thick] (1,0) arc[start angle=0, delta angle=180, radius=1cm];
  \filldraw (1,0) circle (1pt); \filldraw (1,-1) circle (1pt);
  \node at (1.3,0) {$A$};       \node at (1.3,-1) {$B$};
\end{tikzpicture}

\subsection{Cofaces of effective cones}
\label{sec:cofaces of effective cones}

To a subcategory $\SS$ of $\HH$, we associate
\begin{align*}
      E(\SS) & \coloneqq E(\ll(\SS)),     && \text{the effective cone of $\SS$ and} \\
  E(\HH,\SS) & \coloneqq E(\HH) - E(\SS), && \text{the corresponding coface of $E(\HH)$.}
\end{align*}

\begin{definition}
\label{def:face subcategory}
A \defn{face subcategory} of $\HH$ is a subcategory of the form
\[ \SS_\tau \coloneqq \{h \in \HH \mid \ll(h) \in \tau\} \]
for a face $\tau \faceof E(\HH)$. Let $\faceSerre{\HH}$ be the poset of face subcategories of $\HH$ under inclusion.
\end{definition}

\begin{example}
The effective cone of $\coh(\PP^1)$ has three faces, of dimensions 2, 1 and 0; see Example~\ref{ex:cone:projective line}. The face subcategory $\SS_\tau$ of the 1-dimensional face $\tau$ is the category of torsion sheaves. Any non-zero extension-closed subcategory of $\SS_\tau$ has the same effective cone as $\SS_\tau$.
%
%
\end{example}

\begin{lemma} Let $\HH$ be an abelian category of $\LL$ and $\tau$ a face of the effective cone $E(\HH)$.
\label{lem:face subcategories}
  \begin{enumerate}
  \item The face subcategory $\SS_\tau$ is a Serre subcategory of $\HH$.
  \item There are lattice isomorphisms
    \[ \begin{array}{r @{\,} l @{\hspace{1em}} r @{\,} l}
      \faceSerre{\HH} &\isom \Faces{E(\HH)} ,       &   \SS  &\mapsto E(\SS) \\
      \Faces{E(\HH)}  &\isom \Cofaces{E(\HH)}\opp , & E(\SS) &\mapsto E(\HH,\SS) .
    \end{array} \]
  \end{enumerate}
\end{lemma}

\begin{proof}
  (1) If $0\to h' \to h \to h'' \to 0$ is a short exact sequence in $\HH$ then $\ll(h) = \ll(h') + \ll(h'')$ and all three vectors are in $E(\HH)$. The face property of $\tau$ becomes $h \in \SS_\tau \iff h' \in \SS_\tau, h'' \in \SS_\tau$.

  (2) The map $\Faces{E(\HH)} \to \Serre{\HH}$, $\tau \mapsto \SS_\tau$ is injective with image $\faceSerre{\HH}$ and inverse $\SS \mapsto E(\SS)$. The other lattice isomorphism follows from Lemma~\ref{lem:faces and cofaces}.
\end{proof}

\begin{remark}
In general, not every Serre subcategory is a face subcategory.
However, if $\HH$ is length and $\LL=K(\HH)$ then $\faceSerre{\HH} = \Serre{\HH}$ is the set of all Serre subcategories of $\HH$.

Shunya Saito establishes in \cite[\S2.3]{Saito} a bijection between Serre subcategories of an abelian (or exact) category and the faces of its Grothendieck monoid. To construct the effective cone we pass through the Grothendieck group and the lattice $\LL$, so its face poset is a coarser invariant classifying only certain Serre subcategories.
\end{remark}

\section{Cofan and fan of a bounded heart}
\label{sec:heart fan}

\noindent
Let $\DD$ be a triangulated category over a lattice $\LL$, \ie equipped with a choice of homomorphism $\ll \colon K(\DD) \to \LL$.
Henceforth, a \emph{heart $\HH$ in $\DD$} will always mean the heart of a bounded t-structure on $\DD$, although occasionally we write `bounded heart' for emphasis.
For any heart $\HH$ there is a canonical identification $K(\HH)=K(\DD)$, so that all hearts are categories over $\LL$.

The set $\Hearts{\DD}$ of all bounded hearts in $\DD$ is partially-ordered by 
\[
  \KK \geq \HH 
  \iff \KK \subseteq \bigcup_{n\in \N} \HH[n] *\cdots* \HH
  \iff \Homm{\DD}{\KK[{>}0]}{\HH} = 0.
\]
Happel--Reiten--Smal\o\ tilting yields a bijection $(\TT,\FF) \mapsto \FF[1]*\TT$ between torsion pairs $(\TT,\FF)$ in $\HH$ and hearts $\HH[1] \geq \KK \geq \HH$, see \cite[Lemma~1.1.2]{Polishchuk} or Appendix~\ref{appA:tilting}.

In this section we first define the heart cofan $\cofanZ(\HH)$, an integral cofan, \ie a collection of integral cones (monoids) in the lattice $\LL$; see Subsection~\ref{sub:cofans} for the basic notion. Afterwards, we define the heart fan $\fan(\HH)$, a real fan in the $\R$-vector space $\VSdual$, \ie a collection of convex cones. One can define an integral fan $\fan_\Z(\HH)$ dual to $\cofanZ(\HH)$, however we do not explore this here. The reason is that the integral heart fan loses information: the integral fan of an elliptic curve only contains the rational rays; the hearts whose heart cones are irrational are not visible; see Example~\ref{ex:fan:elliptic curve}. We take this as further evidence of the richer structure of the heart cofan.

Readers interested in the heart fan only may jump to Subsection~\ref{sub:heart fan} where we give a self-contained proof of the fan property, not relying on the notions of integral (co)fans.

\subsection{Heart cofan}

\begin{thmdef}
\label{thm:algebraic fan}
Let $\HH$ be a bounded heart in a triangulated category $\DD$ over $\LL$. Let 
\begin{align*}
  \cofanZ(\HH) & \coloneqq \bigcup_{\HH[1] \geq \KK \geq \HH} \Cofaces{E(\KK)} =
                           \{ E(\KK) - E(\SS) \mid \HH[1] \geq \KK \geq \HH \text{ and } \SS\in \faceSerre{\KK} \}
                 \text{ and} \\
  \cofanR(\HH) & \coloneqq \{ \conv(\sigma) \mid \sigma \in \cofanZ(\HH) \} =
                           \{ E(\KK)_\R - E(\SS)_\R \mid \HH[1] \geq \KK \geq \HH \text{ and } \SS\in \faceSerre{\KK} \} .
\end{align*}
Then
\begin{enumerate}
\item $\cofanZ(\HH)$ is an integral cofan in $\LL$.
\item $\cofanZ(\HH)$ does not depend on the ambient triangulated category $\DD$.
\item If $\HH$ is a length heart then $\cofanZ(\HH)$ is complete.
\item $\cofanR(\HH)$ is a real cofan in $\VS$, and the realification of $\cofanZ(\HH)$.
\end{enumerate}
We refer to $\cofanZ(\HH)$ as the \defn{(integral) heart cofan} and to $\cofanR(\HH)$ as the \defn{real heart cofan} of $\HH$, as a category over $\LL$.
We write $\cofanZ(\HH/\LL)$ or $\cofanR(\HH/\LL)$ to stress the lattice.
\end{thmdef}

\begin{proof}
The equality in the definition of $\cofanZ(\HH)$ is Lemma~\ref{lem:face subcategories}(2).
The equality for $\cofanR(\HH)$ follows from Subsection~\ref{sub:real cones} and convex hulls commuting with Minkowski sums and differences.

(1)
By Remark~\ref{rmk:generated cofans} it suffices to check that the sum of any two effective cones is a common coface of each. Suppose $\KK=\FF[1]*\TT$ and $\KK'=\FF'[1]*\TT'$ are hearts with $\HH=\TT*\FF=\TT'*\FF'$. We have $\TT' \subseteq \TT*(\TT'\cap \FF)$ and $\FF' \subseteq (\TT\cap \FF')*\FF$ because $\TT'$ is closed under quotients and $\FF'$ closed under subobjects. Since $E(\TT) \subseteq E(\KK)$ and $-E(\FF) = E(\FF[1]) \subseteq E(\KK)$ this yields
\begin{align*}
  E(\KK)+E(\KK') &= E(\KK) + E(\TT') - E(\FF') \\
                 &= E(\KK) + E(\TT'\cap \FF) - E(\TT\cap \FF') \\
                 &= E(\KK) - E( (\TT'\cap \FF)[1] ) - E(\TT\cap \FF').
\end{align*}
This is a localisation of $E(\KK)$ at a subcone since $(\TT'\cap \FF)[1], \TT\cap \FF' \subseteq \KK$, hence a localisation at a face by Lemma~\ref{lem:faces and cofaces}. Thus $E(\KK)+E(\KK')$ is a coface of $E(\KK)$. By symmetry it is also a coface of $E(\KK')$. Hence $\cofanZ(\HH)$ is a cofan.

(2) The heart cofan $\cofanZ(\HH)$ is determined by the effective cones $E(\KK)$. These have the form $E(\FF[1]*\TT) = E(\TT) - E(\FF)$ for torsion pairs $\HH = \TT*\FF$. This subset of $\LL$ only depends on the torsion pair: the tilted hearts $\FF[1]*\TT$ may depend on $\DD$ but their effective cones $E(\FF[1]*\TT)$ do not. Therefore $\cofanZ(\HH)$ is independent of $\DD$.

(3) This is precisely the content of Lemma~\ref{lem:from charge to torsion pair} below: given $v\in \LLdual$ we are looking for a heart $\KK$ with $v|_\KK \geq 0$. The lemma provides a torsion pair $(\TT,\FF)$ on $\HH$ with $v|_\TT \geq0$ and $v|_\FF \leq0$ so that $v|_\KK \geq0$ for $\KK \coloneqq \FF[1] * \TT$.

(4) Follows from Corollary~\ref{cor: real cofan}.
\end{proof}

\subsubsection*{Length hearts}

If $\HH$ is length and $v\in \VSdual$ then not only is the set of hearts $\KK$ with $v\in C(\KK)$ non-empty, we will give an explicit description of the set. Its elements correspond to torsion pairs $(\TT,\FF)$ in $\HH$ with $v|_\TT \geq 0$ and $v|_\FF\leq0$. The following lemma describes the minimal and maximal such torsion pairs. The corresponding torsion(free) classes are referred to as \emph{numerical} in \cite[Def.~2.11]{Asai21}, and accordingly we call these \defn{numerical torsion pairs}.

\begin{lemma}[{\cite[\S 3.1]{BKT14}}]
\label{lem:from charge to torsion pair}
Let $\HH$ be a length heart in $\DD$ and let $v\in \Hom{\LL}{\R}$. Then the following full subcategories of $\HH$ define torsion pairs $(\TTT{v},\FFFF{v})$ and $(\TTTT{v},\FFF{v})$ in $\HH$:
  \begin{align*}
     \TTT{v}  &\coloneqq \{ h \in \HH \mid v(h'')    > 0 \ \forall\ h \onto h'' \neq 0 \} , 
   & \FFFF{v} &\coloneqq \{ h \in \HH \mid v(h')  \leq 0 \ \forall\ h' \into h \}, \\
     \TTTT{v} &\coloneqq \{ h \in \HH \mid v(h'') \geq 0 \ \forall\ h \onto h'' \} , 
   & \FFF{v}  &\coloneqq \{ h \in \HH \mid v(h')     < 0 \ \forall\ 0 \neq h' \into h \}.
  \end{align*}
\end{lemma}

Since torsion subcategories are closed under quotients, and torsionfree subcategories are closed under subobjects,
\begin{align*}
v\in C(\FF[1]*\TT) & \iff \TTT{v} \subseteq \TT\subseteq \TTTT{v} \ \text{and}\ \FFF{v} \subseteq \FF\subseteq \FFFF{v} \\
                             & \iff \FFFF{v}[1]*\TTT{v} \geq \FF[1]*\TT \geq \FFF{v}[1]*\TTTT{v}.
\end{align*}
This condition becomes more explicit once we turn the torsion pairs $(\TTTT{v},\FFF{v})$ and $(\TTT{v},\FFFF{v})$ into a tripartite decomposition $\HH = \TTT{v} * \HHsst{v} * \FFF{v}$ where 
\[
   \HHsst{v} \coloneqq \TTTT{v} \cap \FFFF{v} = \{ h\in \HH \mid v(h)=0 \text{ and $v(h')\leq 0$ for all $h' \into h$} \}
\]
is the subcategory of $v$-semistable objects. This is a wide, hence abelian, subcategory of $\HH$.
Section~\ref{sec:wall-and-chambers} relates the heart (co)fan to King-semistability and wall-and-chamber decompositions.
The following characterisation of the set of hearts with $v$ in their heart cone is immediate. We thank Parth Shimpi for pointing out this result to us.

\begin{corollary}
\label{cor:heart cones containing v}
Suppose $\HH$ is length and $v\in \VSdual$. Then there is a bijection
\begin{align*}
\{ \text{torsion pairs in}\ \HHsst{v} \} & \lisom \{ \text{hearts $\KK$ between $\HH[1]$ and $\HH$ with $v\in C(\KK)$} \}\\
   (\TT\sst,\FF\sst)                     & \longmapsto \KK \coloneqq ( \FF\sst * \FFF{v})[1] * (\TTT{v} * \TT\sst).
\end{align*}
\end{corollary}

\subsubsection*{Face pairs}

Having given a criterion for when two length heart cones contain the same dual vector, we discuss when two cones in the heart cofan coincide.

\begin{definition}
\label{def:face pair}
Let $\DD/\LL$ be a triangulated category over a lattice and $\HH$ be a heart in $\DD$.

\begin{enumerate}
\item A \defn{face pair} $(\KK,\SS)$ consists of a heart $\KK$ in $\DD$ and a face subcategory $\SS$ of $\KK$.\\
      The set of face pairs is denoted $\FacePairs{\DD}$. We define the subset
      \[ \FacePairs{\HH} \coloneqq \{ (\KK,\SS) \text{ face pair with } \HH[1] \geq \KK \geq \HH \} . \]
\item We say face pairs $(\KK,\SS)$ and $(\KK',\SS')$ are \defn{coface equivalent} if $E(\KK,\SS)=E(\KK',\SS')$.
\end{enumerate}
\end{definition}

Each cone in the heart cofan $\cofanZ(\HH)$ has the form $E(\KK,\SS)$ for some $(\KK,\SS) \in \FacePairs{\HH}$.
The next result gives a homological interpretation of coface equivalence. 

\begin{proposition}
\label{prop:coface equivalence}
Let $(\KK,\SS), (\KK',\SS')\in \FacePairs{\HH}$. The following conditions are equivalent:

\begin{enumerate}[label=(\roman*)]
\item The associated faces coincide: $E(\KK,\SS)=E(\KK',\SS')$.
\item $\thick{}{\SS} = \thick{}{\SS'} \eqqcolon \UU$ and $\KK/\SS = \KK'/\SS'$ as hearts in the common quotient $\DD/\UU$.
\end{enumerate}
\end{proposition}

Let $\Thick{\DD}$ be the poset of thick subcategories of $\DD$ under inclusion. The proposition implies that we can attach a thick subcategory to each cone in the heart cofan:

\begin{definition}
\label{def:labelled cofan}
The \defn{thick label map} of an abelian category $\HH$ over $\LL$ is the monotone map
\[
\Theta \colon \cofanZ(\HH) \to \Thick{\DD}\opp,\quad \Theta(E(\KK,\SS)) = \thick{}{\SS} .
\] 
\end{definition}

\begin{proof}
Suppose $\thick{}{\SS} = \thick{}{\SS'}$ and $\KK/\SS = \KK'/\SS'$. Then for each $k\in \KK$ we can find $k'\in \KK'$ and $s'\in \thick{}{\SS'}$ with $\ll(k) = \ll(k')+\ll(s')$. It follows immediately from this and the fact that $\thick{}{\SS} = \thick{}{\SS'}$ that $E(\KK,\SS) = E(\KK) - E(\SS) \subseteq E(\KK') + E(\SS) + E(\SS') - E(\SS') \subseteq E(\KK') - E(\SS') = E(\KK',\SS')$. By symmetry they are equal.

Now suppose $E(\KK,\SS) = E(\KK',\SS')$. To show that $\thick{}{\SS}=\thick{}{\SS'}$ and $\KK/\SS=\KK'/\SS'$ we consider, for $k\in \KK$, the diagram
\begin{equation}
\label{diag:torsion pair triangles} \tag{$\ast$}
\begin{tikzcd}[row sep = scriptsize]
t_0'[1] \ar[d]        &          & t_1' \ar[d] \\
   f[1] \ar[d] \ar[r] & k \ar[r] & t    \ar[d] \\
f_0'[1]               &          & f_1'
\end{tikzcd}
\end{equation}
whose rows and columns are triangles arising from torsion pairs $\TT*\FF = \HH = \TT'*\FF'$ with $\KK = \FF[1]*\TT$ and $\KK' = \FF'[1]*\TT'$. The row is short exact in $\KK$ and the columns short exact in $\HH[1]$ and $\HH$ respectively. Note that $t_0'\in \TT'\cap \FF $ since $f\in \FF$ which is closed under subobjects, and $f_1'\in \TT\cap \FF'$ as $t\in \TT$ which is closed under quotients. 

The proof of Theorem~\ref{thm:algebraic fan} expressed the sum of effective cones as a coface by the formula $E(\KK)+E(\KK') = E(\KK) -E'$ where $E' \coloneqq E((\TT'\cap \FF)[1]) + E(\TT\cap \FF')$.
The assumption $E(\KK,\SS) = E(\KK',\SS')$ means this is a common coface of $E(\KK)$ and $E(\KK')$, hence a coface of the maximal common coface $E(\KK)+E(\KK') = E(\KK) - E'$. We obtain the coface relation $E(\KK) - E(\SS) = E(\KK,\SS) \cofaceof E(\KK) - E'$. Therefore the minimal face of $E(\KK)$ containing $E'$ is a face of $E(\SS)$. Because $\SS$ is a face subcategory this implies $(\TT'\cap \FF)[1], \TT\cap \FF' \subseteq \SS$. By symmetry we get $(\TT\cap \FF')[1], \TT'\cap \FF \subseteq \SS'$, so that $\TT\cap \FF' \subseteq \SS\cap \SS'[-1]$ and $\TT'\cap \FF \subseteq \SS'\cap \SS[-1]$.
In particular, returning to the diagram \eqref{diag:torsion pair triangles}, we have $t_0'[1] \in \SS \cap \SS'[1]$ and $f_1'\in \SS \cap \SS'[-1]$. 

For the proof of $\thick{}{\SS} = \thick{}{\SS'}$, let $k\in \SS$.
Then also $t, f[1]\in \SS$ because $\SS$ is a Serre subcategory.
Combined with $t_0'[1], f_1' \in \SS$ this shows $\ll(f_0'[1]), \ll(t_1') \in E(\SS)-E(\SS)$.
This Minkowski difference is the minimal face of $E(\KK,\SS)$ hence, by $E(\KK,\SS)=E(\KK',\SS')$, also the minimal face $E(\SS')-E(\SS')$ of $E(\KK',\SS')$. Therefore
\[
  \ll( f_0'[1] ), \ll( t_1') \in E(\KK') \cap (E(\SS') - E(\SS')) = E(\SS').
\]
Now using that $\SS'$ is a face subcategory, we find $f_0'[1], t_1' \in \SS'$. Combined with $t_0', f_1'[1]\in \SS'$, this in turn implies $f[1], t\in \thick{}{\SS'}$ and hence $k\in \thick{}{\SS'}$. Therefore $\thick{}{\SS} \subseteq \thick{}{\SS'}$ and by symmetry the thick subcategories coincide.

Let $\UU \coloneqq \thick{}{\SS}=\thick{}{\SS'}$ be the common thick subcategory and consider $k\in \KK$ in \eqref{diag:torsion pair triangles}. We have shown that $t_0'[1], f_1' \in \UU$. Therefore $k \in \FF'[1]*\TT'$ in the quotient $\DD/\UU$. We get the containment $\KK/\SS \subseteq \KK'/\SS'$ of bounded hearts in $\DD/\UU$, and so they are equal.
\end{proof}

\subsection{The heart fan}
\label{sub:heart fan}

We define the heart fan from HRS-tilting and show that it is the fan associated to the real heart cofan. For readers only interested in the heart fan, we give a self-contained proof that $\fan(\HH)$ is indeed a fan.
Recall from Definition~\ref{def:effective and heart cones} that the heart cone of a heart $\KK$ is
  $C(\KK) = \{ v\in\VSdual \mid v|_\KK \geq 0 \}$;
this is the dual cone of the real effective cone $E(\KK)_\R$.

\begin{thmdef}
\label{thm:heart fan}
Let $\HH$ be a bounded heart in a triangulated category $\DD$ over $\LL$. Let
\[ \fan(\HH) \coloneqq \bigcup_{\HH[1] \geq \KK \geq \HH} \Faces{C(\KK)} \]
be the set of all heart cones $C(\KK)$ for bounded hearts $\KK$ such that $\HH[1] \geq \KK \geq \HH$, together with their faces. Then
\begin{enumerate}
\item $\fan(\HH)$ is the fan in $\VSdual$ associated to the real cofan $\cofanR(\HH)$.
\item $\fan(\HH)$ does not depend on the ambient triangulated category $\DD$.
\item If $\HH$ is a length heart then $\fan(\HH)$ is complete, \ie $\supp{\fan(\HH)} = \VSdual$.
\end{enumerate}
We call $\fan(\HH)$ the \defn{heart fan} of $\HH$, as a category over $\LL$; we write $\fan(\HH/\LL)$ to stress the lattice.
\end{thmdef}

\begin{proof}
This follows from Lemma~\ref{lem:cofans induce fans} about fans associated to cofans and Theorem~\ref{thm:algebraic fan} about the heart cofan $\cofanZ(\HH)$. We give a standalone proof that $\fan(\HH)$ is a fan. By construction, $\fan(\HH)$ is the collection of cones generated by the real heart cones $C(\KK)$ for $\HH[1] \geq \KK \geq \HH$. As $\fan(\HH)$ is closed under taking faces, it suffices to show that the intersection of two generating cones $C(\KK)$, $C(\KK')$ in $\fan(\HH)$ is a face of each.

There are unique torsion pairs $(\TT,\FF)$ and $(\TT',\FF')$ in $\HH$ with $\KK = \FF[1] * \TT$ and $\KK' = \FF'[1] * \TT'$. Each $t'\in \TT'$ sits in a short exact sequence $0\to t_1 \to t' \to f_1 \to 0$ where $t_1 \in \TT$ and $f_1 \in \FF\cap \TT'$ because $\TT'$ is closed under quotients. Similarly each $f'\in \FF'$ sits in a short exact sequence $0\to t_2 \to f' \to f_2 \to 0$ where $t_2 \in \TT\cap \FF'$ and $f_2 \in \FF$, as $\FF'$ is closed under subobjects. Thus
\begin{align*}
v\in C(\KK) \cap C(\KK') 
& \iff v|_\TT \geq 0, v|_\FF \leq 0 \text{ and } v|_{\TT'} \geq 0, v|_{\FF'} \leq 0 \\
& \iff v\in C(\KK) \text{ and } v|_{\FF \cap \TT'} = 0 = v|_{\TT \cap \FF'} \\
& \iff v\in C(\KK) \text{ and } v \in \smash{E\orth}
\end{align*}
where $E \coloneqq E(\TT\cap \FF',\TT'\cap \FF)_\R$ is the cone in $\VS$ generated by the objects in $\TT\cap \FF'$ and $\TT'\cap \FF$. Hence $C(\KK) \cap C(\KK') = C(\KK) \cap E\orth$ is a face of $C(\KK)$. By symmetry it is also a face of $C(\KK')$.

Item (3) follows at once from Lemma~\ref{lem:from charge to torsion pair}. For (2), see the remark below.
\end{proof}

\begin{remark}
By Theorem~\ref{thm:algebraic fan}, the heart cofan $\cofanZ(\HH)$ can be constructed from $\HH$, without recourse to any ambient triangulated category, as the cofan consisting of all the cofaces of the cones $\{ E(\TT)-E(\FF) \mid (\TT,\FF) \text{ torsion pair in}\ \HH \}$. Likewise, the above theorem shows that the heart fan $\fan(\HH)$ is the fan generated by $\{ \cdual{ E(\TT)_\R - E(\FF)_\R} \mid (\TT,\FF) \text{ torsion pair in}\ \HH \}$.

This notwithstanding, avoiding the ambient triangulated category seems unnatural as it obscures the interpretation of cones in the heart fan as heart cones. This becomes more relevant in the sequel where we glue the heart fans together to construct the multifan of $\DD$.
\end{remark}

\subsection{The dual face heart fan}
\label{sec: dual face heart fan}

By Lemma~\ref{lem:face subcategories}, we have an explicit, categorical description of all cofaces of effective cones $E(\KK)$ by face pairs, hence of the entire heart cofan. We do not have a similar description for all faces of heart cones and this is one reason why we consider the cofan as the starting object.
However, we have a categorical description of the dual faces of heart cones, leading us to discuss the dual face heart fan. This circumstance is independent of whether we use integral or real cones; as we defined heart fans to consist of real cones, we will follow suit for their dual face subsets.

Dual faces were introduced in Definition~\ref{def:exposed and dual faces}. Given an abelian category $\HH$ over a lattice $\LL$ and a Serre subcategory $\SS \subseteq \HH$, the \defn{induced sublattice} $\LL_\SS$ is the image of the composite $K(\SS) \to K(\HH) \to \LL$. The quotient $\HH/\SS$ is a category over the lattice $\LL/\LL_\SS$ via the induced homomorphism $K(\HH/\SS) = K(\HH)/K(\SS) \to \LL/\LL_\SS$. 
First we observe that a dual face of a heart cone is the heart cone of a quotient category:

\begin{lemma}
\label{lem:dual faces are heart cones}
The dual face of $C(\HH)$ corresponding to a Serre subcategory $\SS$ of $\HH$ is the heart cone of the quotient $\HH/\SS$, considered as a cone in $\LLdual$ via the canonical map $(\LL/\LL_\SS)\vdual \into \LLdual$:
\[
  \rdual(E(\HH,\SS)_\R) = \rdual( E(\HH)_\R - E(\SS)_\R ) = C(\HH) \cap E(\SS)\orth = C(\HH/\SS) .
\]
\end{lemma}

\begin{proof}
The first equality is merely the definition of the coface $E(\HH,\SS)$.
The middle equality is immediate; this dual face is the image of $C(\HH/\SS) = \cdual{E(\HH/\SS)_\R}$ under the dual $(\LL/\LL_\SS)\vdual \into \LLdual$ of the canonical projection map $\LL \onto \LL/\LL_\SS$.
\end{proof}

By Example~\ref{ex:cofaces adjunction}, the dual faces of a heart cone are in bijection with the exposed faces of the effective cone and therefore also in bijection with certain face subcategories.

\begin{definition}
\label{def:kernel pair}
Let $\HH$ be a heart in a triangulated category $\DD$ be a over the lattice $\LL$.

\begin{enumerate}
\item The \defn{kernel subcategory} of $v\in C(\HH)$ is $\kerc{\HH}{v} \coloneqq \{ h\in\HH \mid v(h) = 0 \}$.
\item A \defn{kernel pair} $(\KK,\SS)$ is a heart $\KK$ in $\DD$ together with a kernel subcategory $\SS$ of $\KK$.\\
      We denote by $\KerPairs{\HH}$ the set of kernel pairs $(\KK,\SS)$ with $\HH[1] \geq \KK \geq \HH$.
\end{enumerate}
\end{definition}

Clearly $\KerPairs{\HH} \subseteq \FacePairs{\HH}$ because $\kerc{\HH}{v}$ is the face subcategory of the exposed face $E(\HH)\cap v\orth$.
We homologically characterise when two dual faces in the heart fan coincide.

\begin{corollary}
\label{cor:cone-equivalence}
The following conditions are equivalent for $(\KK,\SS), (\KK',\SS') \in \KerPairs{\HH}$:
\begin{enumerate}[label = (\roman*)]
\item $(\KK,\SS)$ and $(\KK',\SS')$ are coface equivalent, \ie $E(\KK,\SS)=E(\KK',\SS')$;
\item $C(\KK/\SS) = C(\KK'/\SS')$ as cones in $\VSdual$;
\item $\thick{}{\SS} = \thick{}{\SS'} \eqqcolon \UU$ and $\KK/\SS = \KK'/\SS'$ as hearts in $\DD/\UU$.
\end{enumerate}
\end{corollary}

\begin{proof}
Since kernel subcategories are the face subcategories of exposed faces, Corollary~\ref{cor:Galois correspondence} implies
$
E(\KK,\SS) = E(\KK',\SS') \iff \cdual{E(\KK,\SS)} = \cdual{E(\KK',\SS')} .
$
(This is not true for \emph{face} pairs!)
The equivalence
$
 \cdual{E(\KK,\SS)} = \cdual{E(\KK',\SS')} \iff C(\KK/\SS) = C(\KK'/\SS')
$
follows from Lemma~\ref{lem:dual faces are heart cones} after applying Lemma~\ref{lem:real-integral-faces-cofaces}.
The result then follows from Proposition~\ref{prop:coface equivalence}.
\end{proof}

\begin{remark}
\label{rem: equal hearts}
Distinct hearts may have the same effective cones; if they do, then their heart cones coincide too, but not conversely. The mixed geometric hearts $\mixedHeart{\PP^1}$ of Example~\ref{ex:cone:tilted projective line} all have the same effective cone which differs from the effective cones $E(\coh(\PP^1))$ and $E(\reversedHeart{\PP^1})$. However, they all have the same heart cone.

By contrast, if $0$ is a kernel subcategory of $\KK$ --- for example if $\KK$ is a length heart and $\LL = K(\KK)$ --- then the previous result shows $C(\KK)$ is not contained in $C(\KK')$ for any $\KK'\neq \KK$. In particular, if $C(\KK) = C(\KK')$ for length hearts $\HH[1] \geq \KK,\KK' \geq \HH$ and $\LL = K(\HH)$ then $\KK=\KK'$. This can also be deduced from Corollary \ref{cor:heart cones containing v}: if $v$ lies in the interior of $C(\KK)$ then the category $\HHsst{v}$ of $v$-semistable objects vanishes, so there is a unique heart whose heart cone contains $v$. 
\end{remark}

\begin{remark}[dual face fan]
\label{rem:dual face fan}
According to Remark~\ref{rmk:dual faces and fans}, there is a dual face fan\footnote{In the original \texttt{arXiv} submission, we defined the heart fan of an abelian category $\HH$ over a lattice $\LL$ to be this dual face fan. We believe that the new approach is simpler and better: defining the integral heart cofan first, and inducing from that real heart cofan, the heart fan and only then the heart dual face fan.}
attached to the real heart cofan. It is defined as the set $\rdual(\cofanR(\HH))$ of dual cones of the real cofan; equivalently, it is the subset of the heart fan $\fan(\HH)$ consisting of all maximal cones $C(\KK)$ for hearts $\HH[1] \geq \KK \geq \HH$ together with their dual faces. By Lemma~\ref{lem:dual faces are heart cones} and Corollary~\ref{cor:cone-equivalence}, we have an explicit and categorical description of this dual face fan as
\[ \rdual(\cofanR(\HH)) = \{ C(\KK/\SS) \mid \KK \in \Hearts{\DD}, \HH[1] \geq \KK \geq \HH, \SS \subseteq \KK \text{ kernel subcategory} \} . \]

The following commutative diagram summarises the constructions of the real cofan, fan and dual face fan of the heart $\HH$:
\[ 
\begin{tikzcd}
 \FacePairs{\HH} \ar[two heads,r,"E_\R"]          & \cofanR(\HH) \ar[d,"\rdual"] \\
 \FacePairs{\HH} \ar[equal,u] \ar[r,"C"]          & \fan(\HH) \\
  \KerPairs{\HH} \ar[two heads,r,"C"] \ar[hook,u] & \rdual(\cofanR(\HH)). \ar[hook,u] 
\end{tikzcd} 
\]
The top map is surjective but never injective because $E(\HH,\HH)_\R = E(\HH[1],\HH[1])_\R$.
The fibres of $E_\R$ correspond to coface equivalent face pairs, a categorical condition by Proposition~\ref{prop:coface equivalence}. 

The map $C$ on face pairs is not surjective in general; its image is the dual face fan which generates the heart fan. Example~\ref{ex:fan:countable_semisimple} gives a length (but not algebraic) abelian category $\HH$ whose heart fan is not a dual face fan; we do not know if $C$ is surjective for algebraic $\HH$.

The restriction of $C$ to the subset of kernel pairs also maps onto the dual face fan $\rdual(\cofanR(\HH))$; two kernel pairs have the same image under $C$ if and only if they are coface equivalent, which is a categorical condition by Corollary~\ref{cor:cone-equivalence}.
%
\end{remark}

\begin{remark}
\label{rmk:labelled fan}
By Corollary~\ref{cor:cone-equivalence}, there is a well-defined monotone \defn{thick label map} 
\[
\Theta \colon \rdual(\cofanR(\HH)) \to \Thick{\DD}\opp, \quad \Theta(C(\KK/\SS)) = \thick{}{\SS}
\]
assigning a thick subcategory to each cone in the dual face heart fan; see Definition~\ref{def:labelled cofan} for the cofan version of this map.
Recall from Remark~\ref{rmk:dual faces and fans} that the dual face heart fan $\rdual(\cofanR(\HH))$ has the same support as the heart fan $\fan(\HH)$. Thus each $v\in \supp{\fan(\HH)}$ lies in a unique minimal dual cone. It follows that there is a well-defined map
\[
\Theta \colon \supp{\fan(\HH)} \to \Thick{\DD}, \quad v \mapsto \Theta(C(\KK/\SS)) \quad\text{where $v\in C(\KK/\SS)$;}
\]
we abuse notation by using the same symbol and name as for the map on the cofan and fan, and consider the target simply as a set. More explicitly, $\Theta(v) = \thick{}{\kerc{\KK}{v}}$ for $v\in C(\KK)$.
\end{remark}

\subsection{Lattice change and autoequivalences}
\label{sub:lattice change}

\noindent
We briefly discuss the effect on heart cofans and fans when the lattice is changed. See Example~\ref{ex:elliptic curve:lattice change} for an instance.

\begin{proposition}
\label{prop:lattice change}
Let $\HH$ be an abelian category over $\LL$ and $f\colon \LL \to \Gamma$ a lattice homomorphism.
\begin{enumerate}
\item The heart cofan $\cofanZ(\HH/\Gamma)$ is the pushforward $f(\cofanZ(\HH/\LL))$.\\
      If $f$ is injective it is the isomorphic image of $\cofanZ(\HH/\LL)$ in the sublattice $f(\LL)$ of $\Gamma$.
\item The heart fan $\fan(\HH/\Gamma)$ is the pullback $(f_\R\vdual)^{-1}( \fan(\HH/\LL) )$.\\
      If $f$ is surjective then $\fan(\HH/\Gamma)$ consists of intersections of cones in $\fan(\HH/\LL)$ with the vector subspace $f_\R\vdual(\Gamma_\R\vdual)$ in $\VSdual$.
\end{enumerate}
\end{proposition}

\begin{proof}
(1) follows from the definition of heart cofan and Lemma~\ref{lem:cofans under linear maps}(1). Then (2) follows from (1), Lemma~\ref{lem:cofans under linear maps}(3), Remark~\ref{rem: real cofans under linear maps} and $(f\vdual)^{-1}(\sigma) = \sigma \cap f\vdual(\Gamma\vdual)$ for surjective $f$. 
\end{proof}

Now we discuss the action of autoequivalences; see Examples~\ref{ex:fan:derived-discrete} and \ref{ex:heart fan:coherent sheaves} for instances. We denote by $\Autt{\LL}{\HH}$ the set of pairs $(\Phi,\phi)$ where $\Phi\colon \HH \to \HH$ is an equivalence and $\phi\colon \LL \to \LL$ a group isomorphism such that $\ll \circ [\Phi] = \phi \circ \ll$. Then $\Autt{\LL}{\HH}$ is a group under composition. If $\ll$ is surjective then $\phi$ is uniquely determined by $\Phi$. If $\HH$ is a lattice category and $\LL = K(\HH)$ then $\Autt{\LL}{\HH} = \Aut{\HH}$ with $\phi = [\Phi]$. The same holds for $\HH = \coh(X)$ and $\LL = N(X)$ where $X$ is a smooth projective variety.

\begin{proposition}
\label{prop:aut action}
Let $\HH$ be an abelian category over the lattice $\LL$. Then $\Autt{\LL}{\HH}$ acts on
\begin{enumerate}
   \item the heart cofan $\cofanZ(\HH/\LL)$ by $(\Phi,\phi)\cdot \sigma = \phi(\sigma)$;
   \item the heart fan $\fan(\HH/\LL)$ by $(\Phi,\phi)\cdot \sigma = (\phi_\R\vdual)^{-1}(\sigma)$.
\end{enumerate}
\end{proposition}

\begin{proof}
This follows from Lemma~\ref{lem:cofans under linear maps} and the identifications
  $\phi(\cofanZ(\HH)) = \cofanZ(\Phi(\HH)) = \cofan(\HH)$ and
  $(\phi_\R\vdual)^{-1}( \fan(\HH) ) = \fan(\Phi(\HH)) = \fan(\HH)$
for any $(\Phi,\phi)$ in $\Autt{\LL}{\HH}$.
\end{proof}

\section{Examples of heart fans}
\label{sec:examples}

\noindent
We compute a few heart fans for triangulated categories with lattice $\LL$ of rank up to $2$. Except where indicated otherwise $\LL$ is the Grothendieck group. The pictures are collected on page~\pageref{fig:fans}.
Following \cite{Altmann}, vectors in $\LL = \Z^2$ are written $(m,n)$ and vectors in $\VSdual$ as $[a,b]$. 

Our examples will come from algebras or varieties. We write $\fan(A) \coloneqq \fan(\mod{A} / K(A))$ for the heart fan of modules over a finite-dimensional algebra $A$, and $\fan(X) \coloneqq \fan(\coh(X) / N(X))$ for the heart fan of coherent sheaves on a smooth projective variety $X$ where the lattice $\LL = N(X)$ is the numerical Grothendieck group.

\begin{example}[$\bm{\LL}$ of rank 1]
\label{ex:fan:rank 1}
All fans in $\R$ occur as heart fans.
The complete fans $\{ \R \}$ and $\{ \R_{\leq0}, \zerocone, \Rpos \}$ occur for $\HH = \mod{\kk}$ with $\ll = 0$ and $\ll = \id$, respectively.
The incomplete fan $\{ \Rpos, \zerocone \}$ arises for example as the intersection of the heart fan of $\coh(\PP^1)$ from Example~\ref{ex:fan:projective line} with a line; see Proposition~\ref{prop:lattice change}(2).
Finally, if $E(\HH)_\R = \R$ then the heart fan consists just of $\zerocone$, \eg if $\HH = \mixedHeart{\PP^1}$ of Example~\ref{ex:cone:tilted projective line} with $\ll\colon K(\HH) \to \LL = \Z, \ll(\cO) = 0, \ll(\cO_p) \neq 0$.
\end{example}

\begin{example}[$\bm{A_2}$ quiver]
\label{ex:fan:A2}
Let $\HH = \mod{\kk A_2}$ be the category of representations of the quiver $2 \too 1$. Denote the simple modules by $S_1$ and $S_2$, and let $0 \to S_1 \to E \to S_2 \to 0$ be the non-trivial extension. Fix $[S_1], [S_2]$ as orthonormal basis for $\LL = K(\HH) \cong \Z^2$, so that $[E]$ corresponds to the vector $(1,1)$. There are five hearts between $\HH$ and $\HH[1]$, all of which are algebraic, and accordingly five 2-dimensional cones in $\fan(\kk A_2)$. The hearts are listed below, together with their simple objects, the torsion pair which tilts from $\HH$ towards the heart, and the cone generators. They are of $A_2$ type, \ie equivalent to $\HH$, except for $\KK_\ssi$ which is semisimple.

\begin{center}
\scriptsize
\begin{tabular}{cl c r@{~}l c cc c r@{~}l c r@{~}l} 
  \mc{2}{c}{heart $\KK$} && \mc{2}{c}{simple objects} && $\TT$             & $\FF$           && \mc{2}{c}{$E(\KK)$}  && \mc{2}{c}{$C(\KK)$} \\ \midrule
  & $\HH$          && $S_1$,   & $S_2$           && $\HH$            & $0$             && $(1,0)$, & $(0,1)$   && $[1,0]$,  & $[0,1]$ \\
  & $\KK_1$        && $E$,     & $S_1[1]$        && $\clext{E, S_2}$ & $\clext{S_1}$   && $(1,1)$, & $(-1,0)$  && $[-1,1]$, & $[0,1]$ \\
  & $\KK_2$        && $S_2$,   & $E[1]$          && $\clext{S_2}$    & $\clext{S_1,E}$ && $(0,1)$, & $(-1,-1)$ && $[-1,0]$, & $[-1,1]$ \\
  & $\HH[1]$       && $S_1[1]$,& $S_2[1]$        && $0$              & $\HH$           && $(-1,0)$,& $(0,-1)$  && $[-1,0]$, & $[0,-1]$ \\
  & $\KK_\ssi$     && $S_1$,   & $S_2[1]$        && $\clext{S_1}$    & $\clext{S_2}$   && $(1,0)$, & $(0,-1)$  && $[1,0]$,  & $[0,-1]$
\end{tabular}
\end{center}
\end{example}

\begin{example}[$\bm{2}$-dimensional semisimple algebra]
\label{ex:fan:semisimple}
Let $\HH = \mod{\kk^2}$ be a semisimple category with two simple objects. Similarly to the previous example, but even easier, one can see that the heart fan $\fan(\kk^2)$ has four maximal cones, all of which are smooth.
\end{example}

\begin{example}[tube categories and derived-discrete algebras]
\label{ex:fan:derived-discrete}
Let $\HH$ be the standard stable tube of rank $2$. It is a hereditary abelian length category with two simple objects $S_1$ and $S_2$. It has no non-zero projective or injective objects and each indecomposable object is uniserial. One model is
$\HH = \nilrep(%
    \begin{tikzcd}[cramped, sep=small] 
      2 \ar[r, shift left=0.5ex] & 1 \ar[l, shift left=0.5ex]
    \end{tikzcd})$,
the category of finite-dimensional nilpotent representations of the two vertex quiver with one cycle. (A representation is nilpotent if it has a filtration where all factors are $S_1$ and $S_2$.)
%
%
See \cite{SS} for details on tube categories.

At the bottom of the tube, there are two short exact sequences $0 \to S_1 \to E \to S_2 \to 0$ and $0 \to S_2 \to F \to S_1 \to 0$.
There are six intermediate hearts, described below together with the corresponding torsion pairs:

\begin{center}
\tiny
\begin{tabular}{cl c r@{~}l c cc c r@{~}l c r@{~}l} 
  \mc{2}{c}{heart $\KK$} && \mc{2}{c}{simple objects} && $\TT$             & $\FF$           && \mc{2}{c}{$E(\KK)$}  && \mc{2}{c}{$C(\KK)$} \\ \midrule
  & $\HH$          && $S_1$,   & $S_2$           && $\HH$            & $0$             && $(1,0)$, & $(0,1)$   && $[1,0]$,  & $[0,1]$ \\
  & $\KK_1$        && $E$,     & $S_1[1]$        && $\clext{H \in \HH \mid \top H=S_2}$ & $\clext{S_1}$   && $(1,1)$, & $(-1,0)$  && $[-1,1]$, & $[0,1]$ \\
  & $\KK_2$        && $S_2$,   & $E[1]$          && $\clext{S_2}$    & $\clext{H \in \HH \mid \soc H = S_1}$ && $(0,1)$, & $(-1,-1)$ && $[-1,0]$, & $[-1,1]$ \\
  & $\KK_3$        && $S_1$,   & $F[1]$          && $\clext{S_1}$    & $\clext{H \in \HH \mid \soc H = S_2}$ && $(1,0)$, & $(-1,-1)$ && $[-1,0]$, & $[1,-1]$ \\
  & $\KK_4$        && $F$,     & $S_2[1]$        && $\clext{H \in \HH \mid \top H = S_1}$  & $\clext{S_2}$   && $(1,1)$, & $(-1,0)$  && $[1,-1]$, & $[1,0]$ \\
  & $\HH[1]$       && $S_1[1]$,& $S_2[1]$        && $0$              & $\HH$           && $(-1,0)$,& $(0,-1)$  && $[-1,0]$, & $[0,-1]$
\end{tabular}
\end{center}

\medskip
Denote by $\Lambda_{2,2,0}$ and $\Lambda_{1,2,0}$ the derived-discrete algebras with Gabriel quiver
 $\begin{tikzcd}[cramped, sep=small] 2 \ar[r, shift left=0.5ex] & 1 \ar[l, shift left=0.5ex] \end{tikzcd}$
and zero relations at both vertices or just one vertex, respectively.
The Auslander--Reiten quiver of $\HH$ truncates to those of $\Lambda_{2,2,0}$ and $\Lambda_{1,2,0}$.
Then $\fan(\HH) = \fan(\Lambda_{2,2,0}) = \fan(\Lambda_{1,2,0})$.
 
However, there are differences:
\label{todo:virtual2}
First, $\fan(\HH)$ is not a $g$-fan due to the lack of silting objects but it may be regarded as a `virtual $g$-fan'; see Remark~\ref{rem:virtual g-fan}.
Second, $\fan(\Lambda_{2,2,0}) = \fan(\Lambda_{1,2,0})$ are $g$-fans but behave differently equivariantly:
the quiver and the algebra $\Lambda_{2,2,0}$ have an automorphism $\alpha$ swapping the simple objects and inducing an exact autoequivalence $\alpha_*$.
By Proposition~\ref{prop:aut action}, $\alpha_*$ acts as a non-trivial involution on $\fan(\HH) = \fan(\Lambda_{2,2,0})$, fixing $\pm C(\HH)$ and swapping $C(\KK_1) \mapsto C(\KK_4)$ and $C(\KK_2) \mapsto C(\KK_3)$. In contrast, the same fan automorphism of $\fan(\Lambda_{1,2,0})$ is not induced by an exact autoequivalence of $\mod{\Lambda_{1,2,0}}$.
\end{example}

\begin{remark}
The previous three examples show finite, complete heart fans with $5$, $4$ and $6$ full cones, respectively.
Much more is possible: By \cite[Thm.~4.16]{AHIKM}, each finite, complete and sign-coherent fan in $\R^2$ can be realised as the $g$-fan of a $\tau$-tilting finite algebra, hence as a heart fan; see Section~\ref{sec:g-fans} for the relationship between $g$-fans and heart fans. 

A fan in $\R^2$ is \emph{sign-coherent} if
  (a) there is a basis $e_1, e_2$ of $\Z^2$ such that the cones generated by $e_1, e_2$ and by $-e_1, -e_2$ are in the fan;
  (b) each maximal cone is full, smooth and lies in a quadrant with respect to the basis $e_1, e_2$;
  (c) any ray is contained in exactly two full cones.
\end{remark}

\begin{example}[Kronecker quiver]
\label{ex:fan:kronecker}
Let $\HH = \modsf( \kk( \begin{tikzcd}[cramped, sep=small] 2 \ar[r, shift left=0.5ex] \ar[r, shift right=0.5ex] & 1 \end{tikzcd} ) )$.
Denote the simple modules by $S_1$, $S_2$ and let
  $M_n = [\begin{tikzcd}[cramped, sep=small] \kk^{n-1} \ar[r, shift left=0.5ex] \ar[r, shift right=0.5ex] & \kk^n\end{tikzcd}]$
for $n\geq 1$ be the indecomposable modules in the postprojective component ($M_1 = S_1 = P_1$ and $M_2 = P_2$).
Similarly, let 
$N_n = [\begin{tikzcd}[cramped, sep=small] \kk^n \ar[r, shift left=0.5ex] \ar[r, shift right=0.5ex] & \kk^{n-1}\end{tikzcd}]$
for $n\geq1$ be the indecomposable modules in the preinjective component ($N_1 = S_2 = I_2$ and $N_2 = I_1$).

Fix $[S_1], [S_2]$ as orthonormal basis for $\LL = K(\HH) \cong \Z^2$. There are infinitely many hearts between $\HH$ and $\HH[1]$, of which two series $\{\KK_n\}$ and $\{\KK'_n\}$ are algebraic. Their limit and inverse limit are the reversed geometric heart $\KK_\infty = \reversedHeart{\PP^1}$ and the geometric heart $\KK_{-\infty} = \coh(\PP^1)$, respectively. Between the $\KK_\infty$ and $\KK_{-\infty}$, there are mixed hearts $\KK_U$ parametrised by subsets $U \subseteq \PP^1$; see Example~\ref{ex:cone:tilted projective line}.
In particular, $\KK_\infty = \KK_{\PP^1}$ and $\KK_{-\infty} = \KK_\emptyset$. 
Finally, there is one semisimple heart $\KK_\ssi$.

The table below summarises the simple objects in each heart, describes the torsion pair $(\TT,\FF)$ used to obtain the heart by tilting from $\HH$, the generators of the effective cone of the heart and the generators of the dual cone of the heart.

\begin{center}
\scriptsize
\begin{tabular}{cl c r@{~}l c cc c r@{~}l c r@{~}l} 
  \mc{2}{c}{heart $\KK$} && simple & objects && $\TT$            & $\FF$             && \mc{2}{c}{$E(\KK)$}  && \mc{2}{c}{$C(\KK)$} \\ \midrule
  & $\HH$           && $M_1$,    & $N_1$    && $\HH$             & $0$               && $(1,0)$, & $(0,1)$   && $[1,0]$, & $[0,1]$  \\
  & $\KK_1$         && $M_2$,    & $M_1[1]$ && $\clext{M_2,N_1}$ & $\clext{M_1}$     && $(2,1)$, & $(-1,0)$  && $[0,1]$, & $[-1,2]$ \\
  & $\KK_2$         && $M_3$,    & $M_2[1]$ && $\clext{M_3,N_1}$ & $\clext{M_1,M_2}$ && $(3,2)$, & $(-2,-1)$ && $[-1,2]$,& $[-2,3]$ \\
  & ~\tightverticaldots
  & $\KK_\infty$    && $R_p$:         & $p\in\PP^1$                     && $\clext{R_{p}, N_{n\geq1}}$ & $\clext{M_{n\geq1}}$ && & && $[-1,1]$ & ray \\
  & $\KK_U$         && $R_p, R_q[1]$: & $p \in U, q\in\PP^1\setminus U$ && $\clext{R_{p}, N_{n\geq1}}$ & $\clext{M_{n\geq1}, R_q}$ && & && $[-1,1]$ & ray \\
  & $\KK_{-\infty}$ && $R_q[1]$:      & $q\in\PP^1$                     && $\clext{N_{n\geq1}}$ & $\clext{M_{n\geq1},R_q}$ && & && $[-1,1]$ & ray \\
  & ~\tightverticaldots  
  & $\KK'_2$        && $N_2$,    & $N_3[1]$ && $\clext{N_2,N_1}$ & $\clext{M_1,N_3}$ && $(1,2)$, & $(-2,-3)$ && $[-2,1]$,& $[-3,2]$ \\  
  & $\KK'_1$        && $N_1$,    & $N_2[1]$ && $\clext{N_1}$     & $\clext{M_1,N_2}$ && $(0,1)$, & $(-1,-2)$ && $[-1,0]$, & $[-2,1]$ \\
  & $\HH[1]$        && $M_1[1]$, & $N_1[1]$ && $0$               & $\HH$             && $(0,-1)$,& $(-1,0)$  && $[0,-1]$,& $[-1,0]$ \\
  & $\KK_\ssi$      && $M_1$,    & $N_1[1]$ && $\clext{M_1}$     & $\clext{N_1}$     && $(1,0)$, & $(0,-1)$  && $[1,0]$, & $[0,-1]$
\end{tabular}
\end{center}
\end{example}

\begin{example}[projective line]
\label{ex:fan:projective line}
Let $\HH = \coh(\PP^1)$ and write $\LL = K(\HH) = \Z[\cO] \oplus \Z[\cO_p]$ as in Example~\ref{ex:cone:projective line}.
Let $(\TT,\FF)$ be a torsion pair in $\HH$.
If $\TT$ contains all line bundles then $\TT = \HH$.
If $\TT$ contains no line bundles, then it is of the form $\TTT{M} = \clext{\cO_q \mid q\in M}$ for a subset $M \subseteq \PP^1$. Tilting preserves the heart cone: $C(\HH) = C(\mixedHeart{\PP^1}) = C(\reversedHeart{\PP^1})$; see Example~\ref{ex:tilted curve}.

Thus assume $\TT$ contains a line bundle $\cO(n)$ of minimal degree $n$. Then all skyscraper sheaves are in $\TT$ because torsion classes are closed under quotients. Hence $\TT$ includes the line bundles $\cO(n+1), \cO(n+2), \ldots$ as torsion classes are closed under extensions. Write $\TTT{n} = \clext{\cO(n), \cO_x}$ for this torsion class. The corresponding torsionfree class is $\FF_n = \clext{\cO(m) \mid m<n}$. Denote by $\KK_n$ the heart obtained from the positive tilt of $\HH$ at $(\TTT{n},\FF_n)$. The simple objects of $\KK_n$ are $\cO(n)$ and $\cO(n-1)[1]$. Therefore the effective cone $E(\KK_n)$ is generated by $(1,n)$ and $(-1,1-n)$ and its dual cone $C(\KK_n)$ has generators $[-n,1]$ and $[1-n,1]$.

All this implies that the maximal cones of the fan $\fan(\PP^1)$ are the opposite rays $C(\HH)$ and $C(\HH[1])$ together with a series of 2-dimensional cones $C(\KK_n)$. These fill out a half-plane in $\VSdual$.
The fan $\fan(\PP^1)$ is incomplete because $v(\cO_x) \geq 0$ for any heart $\KK$ and any $v\in C(\KK)$.
Moreover, the tilted hearts from Example~\ref{ex:tilted curve} form an infinite family with the same heart fan: $\fan(\coh(\PP^1)) = \fan(\mixedHeart{\PP^1}) = \fan(\reversedHeart{\PP^1})$.
%
\end{example}

\begin{example}[elliptic curve]
\label{ex:fan:elliptic curve}
Let $\HH = \coh(X)$ for an elliptic curve $X$. The Grothendieck group $K(X) = K(\HH) \cong \Z \oplus \Pic(X) = \Z \oplus \Z \oplus X$ is large. Let $\LL = N(X) \cong \Z^2$ be the numerical Grothendieck group, \ie the quotient of $K(X)$ by the kernel of the Euler pairing, and let
  $\ll\colon K(X) \to N(X)$, $\ll(A) = (\rk(A),\deg(A))$.
We are going to show that the heart fan consists of all rays in the upper half-plane, together with the origin:
  $\fan(X) = \{ e^{i\pi\theta} \Rpos \mid \theta \in [0,1] \} \cup \{ \zerocone \}$.

Denote the slope of a coherent sheaf by $\mu(A) = \deg(A)/\rank(A) \in \Q \cup \{\infty\}$, and consider the full subcategories $\HH_\mu = \coh(X)_\mu$ of semistable sheaves of slope $\mu$. Two general facts about semistable sheaves are: $\HH_\mu$ is a finite length abelian category; these categories are ordered: $\Hom{A_1}{A_2}=0$ if $A_1$ and $A_2$ are semistable with $\mu(A_1) > \mu(A_2)$.
Peculiar to the elliptic curve case is that then $\Hom{A_2}{A_1} \neq 0$ holds; moreover, indecomposable sheaves are semistable, and $\HH_\mu \cong \HH_\infty$ for all $\mu\in\Q$.
As $\HH_\infty$ is the subcategory of torsion sheaves, splitting as a direct sum indexed by points of $X$, each $\HH_\mu$ is a direct sum of rank one tubes. See \cite[II.14]{Polishchuk} for details.
%
%
%
%

For \emph{real} numbers $\delta \in \R$, the additive subcategories $\TTT{\delta} = \clext{\HH_\mu \mid \mu > \delta }$ and $\FF_\delta = \clext{\HH_\mu \mid \mu \leq \delta }$ define mutually different torsion pairs in $\HH$. Let $\KK_\delta = \FF_\delta[1] * \TTT{\delta}$ be the corresponding tilts of $\HH$.
The real cone $E(\KK_\delta)_\R$ is a half-plane containing $\zerocone\times\Rpos$ whose boundary line has slope $\delta$ and so the heart cone $C(\KK_\delta) = \cdual{E(\KK_\delta)_\R}$ is a ray of slope $-1/\delta$ in $\VSdual \cong \R^2$. This extends to $\delta = \pm\infty$ with $\KK_{-\infty} = \HH = \KK_\infty[-1]$ and $C(\KK_{-\infty}) = \Rpos \times \zerocone = -C(\KK_\infty)$. 
We are going to explain that these are all heart cones in $\fan(X)$. They do fill the upper half-plane with rays.

So let $(\TT,\FF)$ be an arbitrary torsion pair in $\HH$ with tilted heart $\KK = \FF[1] * \TT$. Let $A\in\TT$ be an indecomposable object. The category $\TT$ is closed under quotients and extensions, hence contains the top of $A$ and thus the entire tube of $A$. Moreover, for any $\mu > \mu(A)$, there are non-zero morphisms to objects in $\HH_\mu$ and in particular to the simple objects of $\HH_\mu$ which are surjections. Hence $\TT$ contains all $\HH_\mu$ with $\mu > \mu(A)$. Put $\delta \coloneqq \inf\{\mu \in \R\cup\{\pm\infty\} \mid \TT \cap \HH_\mu \neq 0\}$. The real effective cone $E(\KK)_\R$ is a half-plane with the same slope as $E(\KK_\delta)_\R$. The boundary of $E(\KK)_\R$ depends on whether some, all or no tubes from $\HH_\delta$ belong to $\TT$, as in Example~\ref{ex:tilted curve}. Regardless of this, the heart cone is always the same ray $C(\KK) = C(\KK_\delta)$ of slope $-1/\delta$.
\end{example}

\begin{example}
\label{ex:elliptic curve:lattice change}
We continue the previous example: let $N(X) = \Z^2 \to \Gamma \coloneqq \Z$, $(r,d) \mapsto d$. Then $E(\HH/\Gamma)_\R = \Rpos$, and $\fan(\HH/\Gamma)$ is a complete fan on $\R$ although $\HH$ is not a length category. This is an illustration of general lattice change as in Proposition~\ref{prop:lattice change}(2).
\end{example}


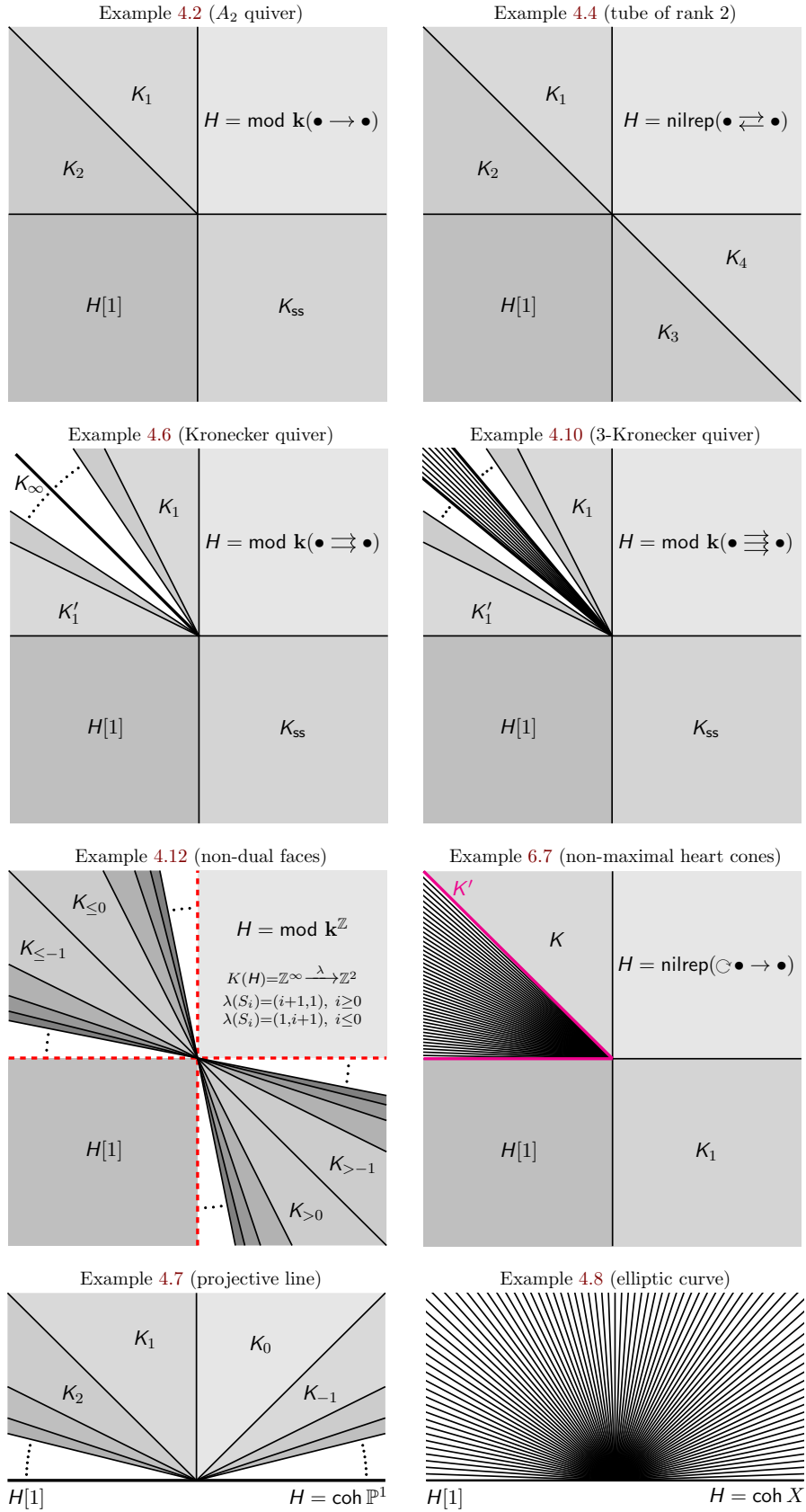
\begin{figure}[p]%
\scalebox{0.79}{%
  \newcommand{\fps}[2]{{\small Example~\ref{#1} (#2)}}
  \newcommand{\fplf}{\\[1.5ex]} 
\begin{tabular}{@{}c @{\hspace{0.02\textwidth}} c@{}}
  \fps{ex:fan:A2}{$A_2$ quiver}                               & \fps{ex:fan:derived-discrete}{tube of rank 2} \\
       \input{diag_heart_fan_A2.tex}                          &      \input{diag_heart_fan_derived-discrete.tex} \fplf
  \fps{ex:fan:kronecker}{Kronecker quiver}                    & \fps{ex:fan:3-kronecker}{3-Kronecker quiver} \\
       \input{diag_heart_fan_Kronecker.tex}                   &      \input{diag_heart_fan_3-Kronecker.tex} \fplf
  \fps{ex:fan:countable_semisimple}{non-dual faces}           & \fps{ex:fan:nilrep}{non-maximal heart cones} \\
        \input{diag_heart_fan_countable_semisimple.tex}       &      \input{diag_heart_fan_nilrep.tex} \fplf
  \fps{ex:fan:projective line}{projective line}               & \fps{ex:fan:elliptic curve}{elliptic curve} \\ 
       \input{diag_heart_fan_projective_line.tex}             &      \input{diag_heart_fan_elliptic_curve.tex}%
\end{tabular}%
}%
\caption{Heart fans of abelian categories $\HH$ over $\LL = \Z^2$. Some maximal cones $C(\KK)$ are labelled by the heart $\KK$. Faces which are not dual faces shown as red dashed lines. Non-maximal heart cones are shown as magenta lines.}%
\label{fig:fans}%
\end{figure}%

\begin{example}[wild quivers]
\label{ex:fan:3-kronecker}
The heart fan of the elliptic curve contains a region entirely filled by rays. This behaviour occurs regularly, \eg for all path algebras of wild quivers; see \cite[Prop.~3.32]{DHKK}. It is related to having a dense region in the phase diagram; see Section~\ref{sec:stability conditions}.

Specifically, and nice to draw due to their rank 2 Grothendieck groups, this holds for the $n$-Kronecker quivers $Q_n$ with two vertices and $n\geq 3$ parallel arrows; see \cite{Chen}. Let $\HH_n = \mod{\kk Q_n}$. By Kac’s Theorem, the class in $K(\HH_n) = \Z^2$ of an indecomposable representation of $Q_n$ is a positive root $(a,b) \in \N^2$ of the associated Euler form
  $q(a,b) = a^2 -nab + b^2$.
The real roots $q(a,b)=1$ correspond to indecomposables in the postprojective and preinjective components of the Auslander--Reiten quiver (there is one indecomposable up to isomorphism for each positive real root).
Let $P_i$ be the indecomposable objects in the postprojective component with $P_1$ and $P_2$ being the indecomposable projectives with respective classes $(0,1)$ and $(1,n)$. The Auslander--Reiten sequences
  $0 \to P_i \to P_{i+1}^n \to P_{i+2} \to 0$
imply that $[P_i] = (a_i, a_{i+1})$ where the $a_i$ satisfy the recurrence $a_{i+2} = n a_{i+1} - a_i$ with $a_0=0$, $a_1=1$. There is a dual picture for the preinjective component where $[I_i] = (a_{i+1}, a_i)$.

The imaginary roots $q(a,b)<0$ correspond to indecomposable objects in the regular component (there are infinitely many indecomposables for each positive imaginary root). Rays through positive imaginary roots are dense in the cone
\[ n - \sqrt{n^2-4} \leq 2x/y \leq n + \sqrt{n^2-4} . \]
The sequences $a_{i+1}/a_i$ and $a_i/a_{i+1}$ converge to the roots $(n\pm\sqrt{n^2-4})/2$ of $x^2-nx+1$.

If $n=2$ then $a_i=i$ and the classes of every regular module lie on the ray of slope 1. For $n>2$, the classes of the regular modules lie in a cone with non-empty interior.
\end{example}

\begin{example}[coherent sheaves]
\label{ex:heart fan:coherent sheaves}
Let $X$ be a connected, smooth and projective variety, $\HH = \coh(X)$ the category of coherent sheaves on $X$ and $\LL = N(X)$ the numerical Grothendieck group of $X$ as in Example~\ref{ex:coherent sheaves heart cone}. Since $\HH$ is noetherian, all torsion pairs in $\HH$ are of the form $(\TT,\TT\orth)$ where $\TT\subseteq\HH$ is an additive subcategory closed under quotients and extensions. If $\TT\neq0$ then it must contain the skyscraper sheaf $\cO_p$ of a point $p\in X$, hence the effective cone of the tilted heart $\KK = \TT\orth[1] * \TT$ contains the point class $[\cO_p] \in E(\KK)$ (as $X$ is connected, all points have the same class in $N(X)$). Choose an orthonormal basis for $\VS = N(X) \otimes \R$ containing the point class and denote its dual basis vector by $v_0 \in \VSdual$. Then $C(\KK)$ is contained in the half-space of $\VSdual$ given by $v_0 \geq 0$. With $C(\HH[1]) = - C(\HH)$ being a ray by Example~\ref{ex:coherent sheaves heart cone}, this shows that the support of the heart fan $\fan(X)$ is contained in a half-space.

A line bundle on $X$ induces an exact autoequivalence of $\coh(X)$ and hence, according to Proposition~\ref{prop:aut action}(2), an automorphism of the fan $\fan(X)$. In this way, the quotient group $\Pic(X)/\Pic^0(X)$ acts faithfully on the fan. For example, the action of $\Pic(\PP^1)$ on $\fan(\PP^1)$ fixes the two rays $\pm C(\HH)$ and acts transitively on the remaining maximal cones: twisting by $\cO(1)$ maps $C(\KK_i) \mapsto C(\KK_{i+1})$ for all $i\in\Z$, using the notation of Example~\ref{ex:fan:projective line}.
\end{example}

\begin{example}[an infinite-dimensional semisimple algebra and non-dual faces]
\label{ex:fan:countable_semisimple}
Let $\HH = \mod{\kk^\Z}$ be the semisimple category with simple objects $S_i$ for $i\in \Z$. Let $\LL=\Z^2$ and define the homomorphism $\ll \colon K(\HH) \to \LL$ by $\ll(S_i) \coloneqq (i+1,1)$ and $\ll(S_{-i}) \coloneqq (1,1+i)$ for $i\in\N$.

Subsets $I\subseteq\Z$ are in bijection with torsion classes in $\HH$; the latter are $\clext{S_i \mid i\in I}$. The corresponding tilt $\KK_I$ of $\HH$ is also semisimple with simple objects $S_i$ for $i\in I$ and $S_i[1]$ for $i \notin I$. It follows that $C(\KK_{\Z \setminus I}) = -C(\KK_I)$. For example $C(\KK_\Z)=C(\HH)$ is the first quadrant and $C(\KK_\emptyset)=C(\HH[1])$ is the third quadrant. The only other non-zero heart cones are those of the hearts $\KK_{\leq n} \coloneqq \KK_I$ where $I = \{ i\in \Z \mid i\leq n\}$ and $\KK_{>n} \coloneqq \KK_I$ where $I = \{ i\in \Z \mid i> n\}$. Here
$C(\KK_{\leq n})$ is generated by $[n,1]$, $[n-1, 1]$ if $n<0$ and by $[-1, n+1]$, $[-1, n+2]$ if $n\geq0$, and $C(\KK_{>n}) = -C(\KK_{\leq n})$. In all other cases $C(\KK_I)=0$; this does not contradict Proposition~\ref{prop:full cones} because $\LL\neq K(\KK_I)$. 

The effective cone $E(\HH) = \{ a,b>0 \} \cup \zerocone$ is not polyhedral, and the heart fan $\fan(\HH)$ contains non-dual faces. Hence the dual face fan $\rdual(\cofanR(\HH))$ of Remark~\ref{rem:dual face fan} is a proper subset of $\fan(\HH)$. The latter contains four extra cones, namely the rays through $[\pm1,0]$ and $[0,\pm 1]$, which are faces of $C(\HH)$ or $C(\HH[1])$ but not dual faces. The four rays are shown in red in Figure~\ref{fig:fans}.
\end{example}

\section{Finite heart fans}
\label{sec:finite heart fans}

\noindent
In this section we answer the natural question of when a heart fan is complete and finite, in the context of an \emph{algebraic} category.
Other heart fans can be finite because of a particular choice of lattice; see Example~\ref{ex:elliptic curve:lattice change} for a complete fan of non-length category. 
The subsequent sections do not depend on the results in this one, so it can be safely skipped. 

We begin with a technical result which is required for the following corollary. It applies for example if $\HH$ is an algebraic category. The full-dimensionality of $C(\KK)$ is crucial: if $\HH$ is the module category of the Kronecker quiver, then one possible tilt is the mixed geometric heart $\KK = \mixedHeart{\PP^1}$ of Example~\ref{ex:cone:tilted projective line} which has $\NN_\KK \neq 0$; however $C(\KK)$ is not full in this case.

\begin{lemma}
\label{lem:ghostfree tilting}
Let $\HH$ and $\KK$ be hearts in $\DD$ with $\HH[1] \geq \KK \geq \HH$. Assume that $K(\DD)$ is a lattice. If $C(\KK)$ is full and $\NN_\HH = 0$ then $\NN_\KK = 0$.
\end{lemma}

\begin{proof}
Let $(\TT,\FF)$ be the torsion pair in $\HH = \TT * \FF$ such that $\KK = \FF[1] * \TT$ is its tilt. Assume $\NN_\KK \neq 0$, \ie there is $0 \neq k \in\KK$ with $0 = [k] \in K(\KK)$. The decomposition sequence of $k\in\KK$ with respect to the torsion pair $(\FF[1],\TT)$ is $0 \to f[1] \to k \to t \to 0$, where $f\in\FF$ and $t\in\TT$.
The decomposition is nontrivial because $[k]=0$, and $[f] = 0$ or $[t] = 0$ in $K(\HH) = K(\KK)$ only when $f = 0$ or $t = 0$ by $\NN_\HH = 0$.
Then the effective cone $E(\KK)$ contains the opposite generators $[f[1]] = [k]-[t] = -[t]$ and $[t]$, so is not strictly convex, contradicting Proposition~\ref{prop:full cones} as $C(\KK)$ was assumed to be full.
\end{proof}

\begin{corollary}
\label{cor:full-alg bijection}
Suppose $K(\DD)$ is a lattice and $\HH$ is a heart with $\NN_\HH=0$. Then $\KK \mapsto C(\KK)$ is a bijection between algebraic hearts $\HH[1] \geq \KK \geq \HH$ and full cones in $\fan(\HH)$.
\end{corollary}

\begin{proof}
As $\DD$ is a lattice category, so is every bounded heart in $\DD$. Thus, by Proposition~\ref{prop:full cones}, $C(K)$ is full if and only if $\NN_\KK$ is a Serre subcategory of $\KK$ and $\KK/\NN_\KK$ is algebraic. Since $\NN_\HH=0$, Lemma~\ref{lem:ghostfree tilting} implies $\NN_\KK=0$ too. Hence $C(\KK)$ is full if and only if $\KK$ is algebraic. Finally, by Remark~\ref{rem: equal hearts}, $C(\KK)=C(\KK')$ for algebraic hearts if and only if $\KK=\KK'$. 
\end{proof}

The following result is a natural generalisation of finiteness and completeness results for $g$-fans of finite-dimensional algebras, see Remark~\ref{rem:finite g-fan}.

\begin{theorem}
\label{thm:finite heart fan}
Let $\HH$ be an algebraic heart and $\LL = K(\HH)$. Then the following are equivalent:
\begin{enumerate}[(i)]
\item The heart fan $\fan(\HH)$ is finite.
\item Each heart $\KK$ in $\DD$ with $\HH[1] \geq \KK \geq \HH$ is algebraic.
\item There are finitely many hearts $\KK$ in $\DD$ with $\HH[1] \geq \KK \geq \HH$.
\end{enumerate}
\end{theorem}

\begin{proof}
$(i) \timplies (ii)$: First we show that each facet $\tau$ of a full cone $\sigma$ in $\fan(\HH)$ is contained in another full cone. Since $\fan(\HH)$ is complete we can choose a sequence in $\supp{\fan(\HH)} \setminus \sigma$ converging to a point in the relative interior of $\tau$. Since $\fan(\HH)$ is finite we can pass to a subsequence lying in a single cone $\sigma' \in \fan(\HH)$ which is not a face of $\sigma$. It follows that $\tau \cap\sigma' \neq \emptyset$. Since the intersection is a face of each, and contains a point of the relative interior of $\tau$, we deduce that $\tau$ is a proper face of $\sigma'$. Therefore $\sigma'$ is a full cone.

Next we show that the set of algebraic hearts between $\HH[1]$ and $\HH$ is closed under irreducible simple tilts (which do not result in hearts outside this interval). Suppose $\KK$ is an algebraic heart between $\HH[1]$ and $\HH$ and $s$ is a simple object of $\KK$. Let $\SS$ be the Serre subcategory generated by $s$. Then $C(\KK/\SS) = C(\KK) \cap [s]\orth$ is a facet of $C(\KK)$. By the first part there is another full cone containing this facet, and by Corollary~\ref{cor:full-alg bijection} this cone has the form $C(\KK')$ for a unique algebraic heart $\KK'$ between $\HH[1]$ and $\HH$. Choose $v$ in the relative interior of the facet. By Corollaries~\ref{cor:heart cones containing v} and \ref{cor:cone-equivalence}, the intersection $\HHsst{v} = \HH \cap \thick{}{\SS}$ is either $\SS$ or $\SS[-1]$, depending on whether $s\in \HH$ or $s[-1] \in \HH$. In each case the only torsion pairs in $\HHsst{v}$ are the trivial ones with vanishing torsion or torsionfree class. So Corollary \ref{cor:heart cones containing v} implies there are precisely two heart cones containing $v$, corresponding to a heart and a simple tilt of that heart. Therefore $\KK'$ is either the positive or negative tilt of $\KK$ at $s$, and in particular this tilt is algebraic.

Finally, we show that every heart $\KK=\FF[1]*\TT$ between $\HH[1]$ and $\HH$ is algebraic. Since $\HH$ is algebraic either $\FF=0$, in which case $\KK=\HH$ and there is nothing to prove, or $\FF$ contains a simple object, $s$, of $\HH$. Let $\HH_1$ be the positive tilt of $\HH$ at $s$. By the above $\HH_1$ is algebraic, and by construction $\HH_1[1] \geq \KK \geq \HH_1$ because $\FF \supseteq \clext{s}$. Repeating, we construct a sequence of algebraic hearts $\HH[1] \geq \KK \geq \cdots > \HH_1 > \HH_0 \eqqcolon \HH$. Since each corresponds to a distinct full cone in the heart fan, and there are only finitely many such cones, the sequence must terminate with $\HH_n=\KK$ for some $n\in \N$. In particular, $\KK$ is algebraic.

$(ii) \timplies (iii)$: Consider the poset of hearts between $\HH$ and $\HH[1]$. By assumption each of these hearts is algebraic. It follows that each element has valency the rank of $K(\HH)$ because the covering relations in the poset correspond to the irreducible simple tilts. If the length of any chain in the poset is also bounded then the poset is finite as claimed. 

Each chain in the poset induces a sequence $(\TT_i,\FF_i)$ of torsion pairs in $\HH$ with $\FF_0 \subseteq \FF_1 \subseteq \cdots$. Since $\HH$ is length, $(\TT,\FF) \coloneqq (\bigcap_i \TT_i\, ,\, \bigcup_i \FF_i)$ is also a torsion pair because $\TT$ is closed under extensions and quotients, $\FF$ is closed under extensions and subobjects, and $\TT = {}\orth\FF$. The corresponding heart $\KK = \FF[1]*\TT$ is, by assumption, algebraic. Therefore each of its finitely many simple objects is either in $\TT$ or $\FF[1]$. It follows that each is in either $\TT_i$ or $\FF_i[1]$ for some $i\in \N$. Hence $\KK \subseteq \FF_i[1]*\TT_i$ are nested bounded hearts, and are therefore equal. It follows that the chain has finite length.

$(iii) \timplies (i)$: The heart fan consists of finitely many heart cones and their faces. In particular, there are finitely many full heart cones, whose union is the support of the heart fan. It follows that each cone in the heart fan is a face of one of these full cones. By Proposition \ref{prop:full cones} each full cone is simplicial, in particular has finitely many faces. Therefore the heart fan is finite.
\end{proof}

\section{Stability conditions and the universal phase diagram}
\label{sec:stability conditions}

\noindent
The main result of this section is an interpretation of the heart fan $\fan(\HH)$ as a `universal phase diagram' for stability functions on $\HH$ with the HN property. We first recall standard notions of stability theory in the abelian setting and then specify what we mean by phase diagram.

\subsection{King and Bridgeland stability}

Alastair King's notion depends on a vector in the real dual of the Grothendieck group; Tom Bridgeland's depends on one in the complex dual.

\begin{definition}
\label{def:King stability}
Let $\HH$ be an abelian category and $v \in \Hom{K(\HH)}{\R}$.
\begin{enumerate}
\item $h\in \HH$ is \defn{$v$-semistable} if $v(h)=0$ and $v(h')\leq 0$ for all subobjects $h' \subseteq h$.
\item $h\in \HH$ is \defn{$v$-stable} if it is $v$-semistable and if $h'\into h$, $v(h')=0$ implies $h'=0$ or $h'=h$.
\item $\HHsst{v}$ denotes the full subcategory of $v$-semistable objects in $\HH$.
\end{enumerate}
\end{definition}

This is \cite[Def.~1.1]{King} although for compatibility with \cite{bridgeland scattering} we have reversed the sign convention.
The subcategory $\HHsst{v} \subseteq \HH$ is wide, \ie closed under extensions, kernels and cokernels. The simple objects of $\HHsst{v}$ are precisely the $v$-stable objects of $\HH$.

Fix a surjective group homomorphism $\ll\colon K(\HH) \to \LL$ onto a lattice. Let $\IH_- \coloneqq \IH \cup \R_{<0} = \{me^{\pi i \phi} \in \C \mid \phi\in(0,1], m>0\}$ be the upper half-plane extended by the negative real axis. A \defn{stability function} on $\HH$ is a group homomorphism $Z \colon \VS \to \C$ with $Z(\ll(h)) \in \IH_-$ for all $0 \neq h\in \HH$; we write $Z(h) \coloneqq Z(\ll(h))$. The \defn{phase} $\phi = \phi(h) \in (0,1]$ of $0\neq h\in\HH$ is defined by $Z(h) = m e^{\pi i \phi} \in \C$. An object $h\in\HH$ is \defn{$Z$-semistable} if $\phi(h') \leq \phi(h)$ for all non-zero subobjects $h' \into h$.
The full subcategory $P(\phi)$ of $Z$-semistable objects of phase $\phi$ is a wide abelian subcategory of $\HH$; it is called the \defn{slice} of phase $\phi$.
Simple objects of $P(\phi)$ are called \defn{$Z$-stable}.
If $\psi>\phi$ then $\Hom{P(\psi)}{P(\phi)} = 0$.
%
%
For an interval $I\subseteq (0,1]$ let $P(I)$ be the extension closure of the subcategories $P(\phi)$ with $\phi\in I$. We abuse notation by writing $P(a,b)$ for $P( (a,b))$ and so on.

A stability function $Z$ has the \defn{Harder--Narasimhan (HN) property} if for each $h\in\HH$ there is a finite chain of inclusions $0 = h_0 \into h_1 \into \cdots \into h_n = h$ such that (a) each quotient $0 \neq h_i/h_{i-1} \in P(\phi_i)$ for some $\phi_i \in (0,1]$ and (b) $\phi_1 > \phi_2 > \cdots > \phi_n$.
%



\subsection{Phase diagrams}

Let $\HH$ be a heart in a triangulated category $\DD$ and $Z \colon \VS \to \C$ a stability function on $\HH$.
We define an open half-plane in the $\R$-dual vector space $\C^* = \Homm{\R}{\C}{\R}$ --- not to be confused with the multiplicative group $\C^\times$ --- by
\[
\IH^* \coloneqq \{ w\in \C^* \mid w(e^{\pi i \phi})=0, w(ie^{\pi i \phi})>0 \ \text{some}\ \phi\in(0,1] \}.
\]
This dual half-plane comes with a map $\theta \colon \IH^* \to (0,1]$ defined by $e^{\pi i\theta(w)} \in \ker(w)$.

For any $v\in\VSdual$, there is a subcategory $\HHsst{v} \subseteq \HH$ of $v$-semistable objects as in Definition~\ref{def:King stability}. (Technically we should write $\HHsst{v \ll}$ but as elsewhere we suppress the map $\ll \colon K(\HH) \to \LL$.) The two notions of semistable subcategories agree under the dual map $Z^* \colon \C\vdual \to \VSdual$:

\begin{lemma}
The subcategory $\HHsst{Z^*(w)} = P(\theta(w))$ for any $w\in \IH^*$.
\end{lemma}

\begin{proof}
Put $v \coloneqq Z^*(w) \in \VSdual$ and $\theta \coloneqq \theta(w) \in (0,1]$.
By Definition~\ref{def:King stability}, a non-zero object $h\in \HH$ is $v$-semistable if $v(h)=0$ and $v(h')\leq 0$ for all non-zero subobjects $h'\subseteq h$. Since $w\in \IH^*$ and $Z$ is a stability function on $\HH$ this is equivalent to $\phi(h)=\theta$ and $\phi(h')\leq \theta$ for all non-zero subobjects $h'\subseteq h$, \ie it is equivalent to $h\in P(\theta)$ being $Z$-semistable.
\end{proof}

\begin{definition}
The \defn{phase diagram} of the stability function $Z \colon \VS \to \C$ is the map
\[
 \PD{Z} \colon \IH^* \to \Wide{\HH}, \quad w \mapsto P(\theta(w)) .
\]
\end{definition}

\begin{remark}
The map $\PD{Z}$ records $Z$-semistable objects of each phase in $(0,1]$ and is a stronger invariant than the `naive phase diagram' $\{ \phi \in (0,1] \mid P(\phi) \neq 0 \}$, tracking non-zero slices.
\end{remark}

Recall the thick label map
  $\Theta \colon \supp{\fan(\HH)} \to \Thick{\DD}$, $v \mapsto \thick{}{\kerc{\KK}{v}}$
from Remark~\ref{rmk:labelled fan}, defined by choosing a heart $\HH[1] \geq \KK \geq \HH$ with $v\in C(\KK)$, where $\kerc{\KK}{v} = \{k \in \KK \mid v(k) = 0\}$.

\begin{proposition}
\label{prop:universal phase diagram}
Let $Z$ be a stability function on $\HH$ with the HN property.
Then the phase diagram of $Z$ factorises over the support of the heart fan of $\HH$:
\[
  \begin{tikzcd}
    \IH^* \ar[rr,"\PD{Z}"] \ar[dr,"Z^*"'] && \Wide{\HH} \\
     & {} \supp{\fan(\HH)} \ar[ur,"\Theta\cap\HH"']
  \end{tikzcd}
\]
\ie $\PD{Z}(w) = \Theta(Z^*(w))\cap \HH$ for all $w\in \IH^*$.
\end{proposition}

\begin{proof}
Let $P$ be the slicing induced by $Z$ on $\Db(\HH)$. As $Z$ satisfies the HN property, for any $\psi\in(0,1)$ there are two torsion pairs
  $\HH = P(\psi,1] * P(0,\psi] = P[\psi,1] * P(0,\psi)$,
and therefore the positive tilts $P(\psi,\psi+1]$ and $P[\psi,\psi+1)$ of $\HH = P(0,1]$ exist.

If $w\in \IH^*$ then $w(e^{i\pi\theta(w)})=0$ and $w(ie^{i\pi\theta(w)})>0$, hence $Z^*(w) = w\circ Z$ is in the heart cone of the positive tilt
 $\KK\coloneqq P[\theta(w), \theta(w)+1)$ of $\HH$, and so is in the support $\supp{\fan(\HH)}$. Moreover $\kerc{\KK}{Z^*(w)} = P(\theta(w))$. Thus $\Theta(Z^*(w))= \thick{}{P(\theta(w))}$, hence $\Theta(Z^*(w)) \cap \HH = P(\theta(w))$.
\end{proof}

\begin{remark}
The previous result can be viewed as a description of `half' of the pullback of the heart fan along $Z^* \colon \C^* \to \VSdual$. The other `half' is similarly related to $Z$-semistable objects for the opposite ordering on the real numbers.

There is an essentially equivalent description of `half' of the push-forward of the heart cofan along $Z \colon \VS \to \C$. Its cones in the half-space $\phi\leq\arg(z)\leq\phi+1$ are $Z(P(\phi,\phi+1])$, its coface $Z(P[\phi,\phi+1])$ and when $\phi>0$ also $Z(P[\phi,\phi+1))$. These are the same if $P(\phi)=0$, but are otherwise distinct.
\end{remark}

\begin{example}[non-maximal heart cones]
\label{ex:fan:nilrep}
A heart cone can be a proper face of another one: below we exhibit heart cones $\KK''$, $\KK'$, $\KK$ such that $\zerocone = C(\KK'') \pfaceof C(\KK') \pfaceof C(\KK)$ have dimensions 0, 1 and 2, respectively.
We are grateful to Parth Shimpi for showing us a rank 3 example in the context of contraction algebras \cite{Shimpi} with such behaviour.

Let $Q$ be the quiver with two vertices connected by an arrow and with a loop:
$\begin{tikzcd}[column sep = small, cramped]
\arrow[loop left] 1 \ar[r] & 2
\end{tikzcd}$.
The category $\HH = \nilrep(\kk Q)$ of nilpotent representations of the quiver $Q$ is an algebraic hereditary abelian category with two simple objects $S_1$ and $S_2$. We use $[S_1], [S_2]$ as an orthonormal basis for $\LL = K(\HH) \cong \Z^2$.
For more details, see \cite[p. 428]{Koenig-Yang} and \cite[Ex.\ 2.9]{BCSPW24}.

The indecomposable nilpotent representations of $Q$ are of the form $M(n_0; n_1, \ldots, n_k)$ with $0 \leq k-1 \leq n_0 \leq k$, $n_1 \geq 0$ and $n_j \geq 1$ for each $j > 1$, in which $n_0$ is the number of occurrences of $S_2$ in a composition series, $n_1$ is the number of occurrences of $S_1$ above the first $S_2$ in the composition series, $n_2$ is the number of occurrences of $S_1$ above the second $S_2$ and so on:
\[
  M(0;1) = S_1, \quad M(1;0) = S_2,
  \quad
  M(1;2) = \raisebox{-1.7ex}{%
             \begin{tikzpicture}[every node/.style={scale=0.7, font=\sffamily}]
               \node at (0.0,-0.44) {2}; \node at (0,-0.22) {1}; \node at (0,0) {1};
             \end{tikzpicture}},
  \quad
  M(2;2,1) = \raisebox{-2.4ex}{%
             \begin{tikzpicture}[every node/.style={scale=0.7, font=\sffamily}]
              \node at (0,0) {1}; \node at (0,-0.22) {1}; \node at (-0.2,-0.44) {2}; \node at (0.2,-0.44) {1}; \node at (0.2,-0.66) {2};
            \end{tikzpicture}}%
  .
\]

There are four obvious algebraic hearts between $\HH$ and $\HH[1]$. Besides the trivial ones $\HH$, $\HH[1]$, the simple tilt at $S_1$ produces a heart $\KK_1$ whose heart cone is the fourth quadrant, and the simple tilt at $S_2$ produces a heart $\KK$ with heart cone generated by $[0,1]$ and $[-1,1]$. 
The simple objects of $\KK$ are $S_2[1]$ and $E \coloneqq M(1;1)$. By a phase gap argument as in \cite[Ex.\ 2.9]{BCSPW24}, one can show that the simple tilt of $\KK$ at $E$ is not algebraic.

We prove that the remaining region in $\VSdual$ is filled by 1-dimensional heart cones.
Consider the stability function $Z\colon \VS \to \C$ defined by $Z(S_1) = -1$ and $Z(S_2) = i$. Because $\HH$ is algebraic, $Z$ satisfies the support and HN properties. If $M$ is an indecomposable nilpotent representation then $Z(M)$ lies in the region between the ray $e^{3i\pi/4} \Rpos$ and the negative real axis, and one can check that any subrepresentation has smaller phase so that $M$ is $Z$-semistable.

Fix $\phi = p/q \in [3/4,1] \cap \Q$ where $p,q\in\N$, and choose an indecomposable $M\in\HH$ with $\arg Z(M) = \pi\phi$.
Define a torsion pair $(\TTT{\phi},\FFF{\phi})$ by $\TTT{\phi} \coloneqq P(0,\phi)$ and $\FFF{\phi} \coloneqq P(\phi,1] * \clext{M}$.
The effective cone of the associated tilt $\KK_\phi$ is the half-plane $E(\KK_\phi) = \{(m,n) \in \VS \mid q n \geq p m\}$.
Therefore, the heart cone is $C(\KK_\phi) = \{ [x,y] \mid p y = -q x, x \leq 0\} = e^{i\pi\phi} \Rpos$, the ray of slope $\phi$.
By Theorem~\ref{thm:heart fan}, the heart fan $\fan(\HH)$ is complete, and hence the above region is filled with rays.
At its boundary, we get a heart $\KK' \coloneqq \smash{\KK_{3/4}}$ whose heart cone $C(\KK')$ is a face of $C(\KK)$.
Moreover, $\KK'$ is the non-algebraic heart obtained by tilting $\KK$ at the simple object $E = M(1;1)$.

Finally, we present a heart $\KK''$ with $C(\KK'') = \zerocone$. Let $\KK''$ be the tilt of $\HH$ at the torsion pair defined by the torsionfree class
  $\FF'' \coloneqq \clext{\sub{M(1;2)}}$. 
Then the representations $T_n \coloneqq M(n; 3, 1, 1, \ldots, 1)$ with $k = n$ lie in $\TT'' \coloneqq {}\orth\FF''$. The cone $E(\FF'')$ is bound by the ray of slope $2$, while for $n$ sufficiently large $T_n$ lies in the region bound by slopes $1$ and $2$. Hence, $E(\TT'')$ and $E(\FF'')$ intersect in a full cone, $E(\KK'') = E(\TT'') - E(\FF'') = \LL$ and $C(\KK'') = \zerocone$.

The torsion pairs $(\TTT{\phi},\FFF{\phi})$ are numerical; by contrast, $(\TT'',\FF'')$ is not.
\end{example}

\begin{example}[smooth projective curve]
\label{ex:smooth projective curve}
Let $\HH = \coh(X)$ be the category of coherent sheaves over a smooth, projective curve $X$ of genus $g$ over $\LL = N(X)$ as in Example~\ref{ex:tilted curve}.
The map
  $Z \colon \VS \to \C$, $Z(F) = -\deg(F) + i\,\rk(F)$
is the primordial geometric stability function; it satisfies the support and HN properties.

Line bundles are stable; if $g=0$ no other stable objects exist. This is $\fan(\PP^1)$ of Example~\ref{ex:fan:projective line}.

Let $g>0$ and $d,r\in\Z$ be coprime with $r>0$. Then there exists a stable vector bundle $F$ on $X$ of rank $r$ and degree $d$, \ie $F \in P(\phi)$ where $\phi \coloneqq \arg Z(F) \in (0,1)$.
%




Let $\KK_\phi$ be the tilt of $\HH$ for the torsion pair $(P(0,\phi],P(\phi,1])$. Then $\clE{\KK_\phi}$ is a half-plane and $C(\KK_\phi)$ a ray. Hence the heart fan $\fan(\HH)$ contains a dense family of 1-dimensional cones in the half-plane defined by the condition $\rk\geq0$ in $\VSdual$. We cannot use the completeness argument of the previous example but we can argue as in Example~\ref{ex:fan:elliptic curve}. We find that the heart fan $\fan(X)$ consists of rays covering a half-plane for all $g\geq1$.
\end{example}

\section{The wall-and-chamber structure of an algebraic heart}
\label{sec:wall-and-chambers}

\noindent
We relate the heart fan of an \emph{algebraic} abelian category $\HH$ to the wall-and-chamber structure of \cite{BST} defined via King semistability (Definition~\ref{def:King stability}). More precisely, we use the dual face heart fan $\rdual(\cofanR(\HH))$ of Subsection~\ref{sec: dual face heart fan}. It is complete by Remark~\ref{rmk:dual faces and fans} and Theorem~\ref{thm:heart fan}, \ie $\supp{\rdual(\cofanR(\HH))} = \supp{\fan(\HH)} = \VSdual$ where $\LL = K(\HH)$ and $\DD = \Db(\HH)$.

We begin by showing that among all kernel pairs describing the same cone in the dual face heart fan there is a unique one closest to the reference heart $\HH$. Recall from Definition~\ref{def:kernel pair} that a kernel pair $(\KK,\SS) \in \KerPairs{\HH}$ consists of a heart $\KK$ with $\HH[1] \geq \KK \geq \HH$ and a face subcategory $\SS = \kerc{\HH}{v} = \{ h\in\HH \mid v(h) = 0 \}$ for some $v\in C(\HH)$.
Below, $\relintC{\KK^v \rslash \SS^v}$ denotes the relative interior of $C(\KK^v\rslash\SS^v)$, \ie the interior of the cone inside its linear hull.

\begin{proposition}
\label{prop: distinguished kernel pair}
Let $\HH$ be an algebraic abelian category and $\LL = K(\HH)$ its Grothendieck group. For each $v\in \VSdual$ there is a unique kernel pair $(\KK^v,\SS^v)$ in $\DD = \Db(\HH)$ such that
\begin{enumerate}
\item $C(\KK^v \rslash \SS^v)$ is the minimal dual face containing $v$,
\item $\KK^v$ is minimal amongst hearts $\HH[1] \geq \KK \geq \HH$ with $v\in C(\KK)$ and
\item $\SS^v = \HHsst{v} \subseteq \HH \cap \KK^v$ is the subcategory of $v$-semistable objects in $\HH$.
\end{enumerate}
Moreover, if $w \in \relintC{\KK^v \rslash \SS^v}$ then $(\KK^w,\SS^w)=(\KK^v,\SS^v)$.
\end{proposition} 

We call kernel pairs of the form $(\KK^v,\SS^v)$ \defn{distinguished}.
Recall the surjection from kernel pairs to the dual face heart fan
  $\KerPairs{\HH} \onto \rdual(\cofan_\R(\HH))$, $(\KK,\SS) \mapsto C(\KK/\SS)$;
see Remark~\ref{rem:dual face fan}. By the proposition, this map has a section sending a dual cone $C(\KK/\SS) \in \rdual(\cofan_\R(\HH)$ to the distinguished kernel pair $(\KK^v,\SS^v)$ with $C(\KK/\SS) = C(\KK^v/\SS^v)$ or, equivalently, $v\in\relintC{\KK^v \rslash \SS^v}$.

\begin{proof} Let $\HH = \TTTT{v} * \FFF{v}$ be the torsion pair of Lemma~\ref{lem:from charge to torsion pair}, \ie
$\TTTT{v} \coloneqq \{ h \in \HH \mid v(h'') \geq 0 ~\forall h \onto h'' \}$ and
$\FFF{v}  \coloneqq \{ h \in \HH \mid v(h') < 0 ~\forall\, 0 \neq h' \into h \}.$
Let $\KK^v \coloneqq \FFF{v}[1] * \TTTT{v}$ be the tilted heart; by construction $v\in C(\KK^v)$.
The subcategory
  $\SS^v \coloneqq \{ t\in \TTTT{v} \mid v(t)=0 \} = \{ k\in \KK^v \mid v(k)=0 \} \subseteq \KK^v$
is the face subcategory of the exposed face $E(\KK^v) \cap v\orth$, \ie a kernel subcategory of $\KK^v$. Then $C(\KK^v \rslash \SS^v)$ is a dual face of $C(\KK^v)$ containing $v$, and is minimal with this property.

If $v\in C(\KK)$ with $\HH[1] \geq \KK \geq \HH$ then $\KK = \FF[1]*\TT$ for a torsion pair $\HH = \TT*\FF$ where $v|_\FF\leq 0$. 
Considering a decomposition of an object of $\FF_v$ with respect to $(\TT,\FF)$ and using that $\FF_v$ is closed under subobjects, we get $\FFF{v} \subseteq \FF$, and hence $\KK^v \leq \KK$. Thus $\KK^v$ is minimal amongst hearts $\HH[1] \geq \KK \geq \HH$ with $v\in C(\KK)$. In particular, $v$ determines $\KK^v$ uniquely.

To show that $\SS^v \subseteq \TTTT{v} \subseteq \HH$ is the subcategory of $v$-semistable objects, let $h\in \SS^v$ and $h'\into h$ be a subobject. Then $0 \leq v(h/h') = v(h)-v(h') = -v(h')$ by $h\in \TTTT{v}$, and hence $h\in \HHsst{v}$.
Conversely if $h\in \HHsst{v}$ then it has a short exact decomposition sequence $0\to t\to h \to f\to 0$ in $\HH$ with $t\in \TTTT{v}$ and $f\in \FFF{v}$. Since $h$ is semistable, we have $v(t)\leq 0$ for its subobject $t$. On the other hand, $v(t) \geq 0$ for all objects of $\TTTT{v}$. Therefore $v(t)=0$ and $v(f)=v(h)-v(t)=0$ too, so that $f=0$ by definition of $\FFF{v}$. Hence $h = t \in \SS^v$.

Finally, if $C(\KK^v \rslash \SS^v)$ is the minimal dual face containing $w$ then $\KK^w = \KK^v$ by minimality. Thus $\SS^w = \SS^v$ too, as each is the kernel subcategory of the dual face $C(\KK^v \rslash \SS^v)$. 
\end{proof}

\begin{remark}
If the minimal dual face containing both $v$ and $w$ is $C(\KK/\SS)$ for some $\HH[1] \geq \KK \geq \HH$ then the subcategories of semistable objects $\HHsst{v} = \HHsst{w}$ are the same. The converse does not hold in general; for example $\HHsst{v}=0$ for any $v$ in the interior of a full cone in $\fan(\HH)$.
\end{remark}

\begin{example}[distinguished kernel pair in Kronecker quiver heart fan]
Consider the heart fan of the module category over the Kronecker quiver, $\HH = \mod{\kk(\Kronecker)}$.
This fan contains a limiting ray denoted $C(\KK_\infty)$ in Example~\ref{ex:fan:kronecker}. This ray is the dual cone of the infinitely many hearts of $\Db(\HH)$ listed in Example~\ref{ex:cone:tilted projective line}, among them the reversed geometric heart $\reversedHeart{\PP^1}$.

The distinguished kernel pair $(\KK^v, \SS^v)$ of Proposition~\ref{prop: distinguished kernel pair} for $0 \neq v\in C(\KK_\infty)$ is  $\reversedHeart{\PP^1}$ with its face subcategory of (shifted) torsion sheaves: as the geometric and mixed hearts are positive tilts of  $\reversedHeart{\PP^1}$, this heart is indeed minimal amongst $\KK$ with $\HH[1] \geq \KK \geq \HH$ and $v\in C(\KK)$.
\end{example}

\begin{lemma}
\label{lem:stable object}
An object of $\HHsst{v}$ is stable if and only it is simple as an object of $\KK^v$.
\end{lemma}

\begin{proof}
Suppose $h\in \HHsst{v}$. If $0\to h' \to h \to h/h' \to 0$ is short exact in $\KK^v$ then, since $\HHsst{v}=\SS^v$ is a Serre subcategory of $\KK^v$, both $h'$ and $h/h'$ are in $\HHsst{v}$ so the sequence is also short exact in $\HH$. It follows that if $h$ is stable it is simple in $\KK^v$.

Conversely, if $h''$ is a subobject of $h$ in $\HH$ with $v(h'')=0$ then there is a short exact sequence $0\to h'' \to h\to h/h''\to 0$ in $\HH$ with all terms in $\HHsst{v}$, which is therefore also short exact in $\KK^v$. Thus if $h$ is simple in $\KK^v$ then it is stable.
\end{proof}

\subsection{Wall-and-chamber structures}

Introduced in \cite{BST} for finite-dimensional algebras, we explain how to reconstruct walls and chambers from the heart fan $\fan(\HH)$ of an algebraic abelian category $\HH$. The \defn{stability space} \cite[Def.~3.2]{BST} of an object $h\in\HH$ is 
\[
  \cD(h) \coloneqq \{ v \in \VSdual \mid h\in \HHsst{v} \}.
\]
This is a \defn{wall} if $\codim_\VSdual \cD(h) = 1$. We write $\cD \coloneqq \bigcup_{h\neq0} \cD(h)$ for the union of stability spaces of non-zero objects. A \defn{chamber} is a connected component of the complement $\VSdual \setminus \smash{\overline{\cD}}$.
The \defn{wall-and-chamber structure} is the set of stability spaces $\cD(h)$ for indecomposable objects $h$ together with the set of chambers.
While walls are defined integrally in $\LL^*$ (more precisely, a wall $\cD(h)$ is determined by $\cD(h)\cap\LLdual$), chambers are only defined in $\VSdual$.

\begin{theorem}
\label{thm: wall and chamber}
The stability space $\cD(h)$ of $h\in \HH$ is a rational, polyhedral cone. Moreover,
\[
\Fsup(h) \coloneqq \{ C(\KK/\SS) \mid (\KK,\SS) \text{ distinguished kernel pair with } h\in \SS\} 
\]
is a subset of $\fan(\HH)$ forming a dual face fan with support
$\supp{\Fsup(h)} = \cD(h)$. 

The chambers are the open cones $\relintC{\KK}$ for the algebraic hearts $\HH[1] \geq \KK \geq \HH$, and each chamber is the interior of a full, smooth cone.
\end{theorem}

\begin{remark}
This provides an alternative proof of \cite[Prop.~3.15]{BST} and \cite[Thm.~13.7]{Asai21} for $\HH = \mod{A}$ for a finite-dimensional algebra $A$ which combined show that the chambers are in natural bijection with $\tau$-tilting pairs, \ie with algebraic hearts in the interval between $\HH$ and $\HH[1]$. Always $\supp{\Fsup(h)} \subseteq \supp{\fan(\HH)} \cap [h]\orth$ but $\Fsup(h)$ need not contain all dual faces in this hyperplane section of the heart fan; it only contains those in heart cones $C(\KK)$ where $h \in \KK$.
\end{remark}

\begin{proof}
\proofstep{Stability spaces are rational, polyhedral cones}
The set $S \coloneqq \{ [h'] \mid h' \subseteq h \} \subset K(\HH)$ of classes of subobjects of $h$ is finite because $\HH$ is a length category, and therefore it generates a rational, polyhedral cone $E(S)_\R \subset \VS$. Hence its dual cone $C(S) \subset \VSdual$ is rational, polyhedral as well (this is Gordan's Lemma in convex geometry), as is $\cD(h) = -C(S) \cap [h]\orth$. 
 
\proofstep{$\Fsup(h)$ is a dual face subfan}
Let $(\KK,\SS)$ be a distinguished kernel pair with $h\in \SS$. Then $h$ is $v$-semistable for any $v\in \relintC{\KK/\SS}$ by Proposition~\ref{prop: distinguished kernel pair}; as semistability is a closed condition, this holds for any $v \in C(\KK/\SS)$.
Hence every dual face of $C(\KK/\SS)$ in $\fan(\HH)$ is also in $\Fsup(h)$. So $\Fsup(h)$ is a dual face fan, and it is a subset of $\fan(\HH)$ anyway.

\proofstep{Stability space as support}
The support $\supp{\Fsup(h)} \subseteq \VSdual$ is the union of all cones in $\Fsup(h)$. We have $\cD(h) \supseteq \supp{\Fsup(h)}$ by the previous step.
Conversely, suppose $v \in \cD(h)$. Then $h$ is $v$-semistable and, as $\rdual(\cofanR(\HH))$ is complete, there is a minimal dual face $C(\KK/\SS)$ containing $v$ where $(\KK,\SS)$ is a distinguished kernel pair with $h \in \SS=\HHsst{v}$ by Proposition~\ref{prop: distinguished kernel pair}, and hence $\cD(h) \subseteq \supp{\Fsup(h)}$.

\proofstep{Chambers are interiors of full cones}
Suppose $\KK$ is an algebraic heart with $\HH[1] \geq \KK \geq \HH$. Then $C(\KK)$ is a full cone, so the relative interior $\relintC{\KK} \subseteq \VSdual \setminus \cD$ is an open subset.
Each proper face of $C(\KK)$ has the form $C(\KK/\SS)$ for $\SS\neq 0$ and thus is contained in $\cD$. Hence $\relintC{\KK}$ is a chamber. 

Conversely, suppose $v\in \VSdual$ is in a chamber $\cC$. As $\rdual(\cofanR(\HH))$ is complete, $v$ is contained in a minimal dual face $C(\KK/\SS)$ where $(\KK,\SS)$ is a distinguished kernel pair with $\SS=\HHsst{v} = 0$. It suffices to show that $\relintC{\KK}$ is open in $\VSdual$ for then $\KK$ is algebraic by Proposition~\ref{prop:full cones} and Lemma~\ref{lem:ghostfree tilting}, and $\cC = \relintC{\KK}$ by the first part.

Let $U \subset \VSdual$ be a convex open subset with $v \in U \subset \cC$.
Then $\SS^w = \HHsst{w} = 0$ for all $w\in U$.
We claim $\TTTT{w} = \TTTT{v}$ for all $w\in U$. Suppose contrariwise $\TTTT{w} \neq \TTTT{v}$ for some $w\in U$. As $\HH$ is algebraic, we may assume $t\in \TTTT{v} \setminus \TTTT{w}$, and that every proper quotient of $t$ is in $\TTTT{w}$, \ie
  $v(t) \geq 0 > w(t)$ 
and $v(t''), w(t'')\geq 0$ for every proper quotient $t''$ of $t$. Thus there is $u\in U$ on the line segment between $v$ and $w$ with $u(t)=0$ and $u(t'')\geq 0$ for all quotients $t''$ of $t$. Then $u(t') = u(t)-u(t/t') = -u(t/t') \leq 0$ for all subobjects $t'$ of $t$, hence $0 \neq t\in \HHsst{u}$, contradicting $u\in U \subset \cC$. We conclude that $\TTTT{w} = \TTTT{v}$ for all $w\in U$ after all. Therefore $\KK^w = \KK$ for all $w\in U$, \ie $U \subset \relintC{\KK}$. Thus $\relintC{\KK}$ is open and we are done.
\end{proof}

\begin{definition}
\label{def:stability fan}
The \defn{stability fan} of $\HH$ is the dual face fan 
\[
\fanst(\HH) \coloneqq \{ C(\KK/\SS) \in \fan(\HH) \mid (\KK,\SS) \text{ kernel pair with } \SS\neq 0 \} .
\]
\end{definition}

\begin{remark}
\label{rem:stability fan}
It suffices to use the set of distinguished kernel pairs in the definition. The stability fan is a dual face fan by Proposition~\ref{prop: distinguished kernel pair}, and is the union of the dual face fans $\Delta(h)$ for $0 \neq h \in \HH$. Its support $\supp{\fanst(\HH)} = \bigcup_{h\neq 0} \cD(h)$ is the union of stability spaces of non-zero modules, justifying the name. Each stability space $\cD(h)$ is a rational polyhedral cone, even though $\fanst(\HH)$ may contain cones which are not rational, see Remark~\ref{rem:non-rational stability fan}.
\end{remark}

\begin{example}[3-Kronecker quiver]
\label{ex:3-Kronecker fans}
Let
  $A \coloneqq \kk(
       \begin{tikzcd}[cramped, sep=small] 
         \bullet \ar[r, shift left=0.65ex] \ar[r, shift right=0.65ex] \ar[r] & \bullet
       \end{tikzcd})$
be the path algebra of the quiver with two vertices and three parallel arrows; this is the smallest instance of Example~\ref{ex:fan:3-kronecker}. Its stability fan $\fanst(A)$ is the subfan of the heart fan $\fan(A)$ where all full cones and all non-rational cones are removed. By Theorem~\ref{thm:g-fan}, the $g$-fan $\fan\gfan_A$ is the subfan of the heart fan $\fan(A)$ generated by full cones. We show these three fans side by side, with their supports shaded grey. The magenta rays in the heart fan symbolise the non-rational maximal cones which are absent in the stability fan.

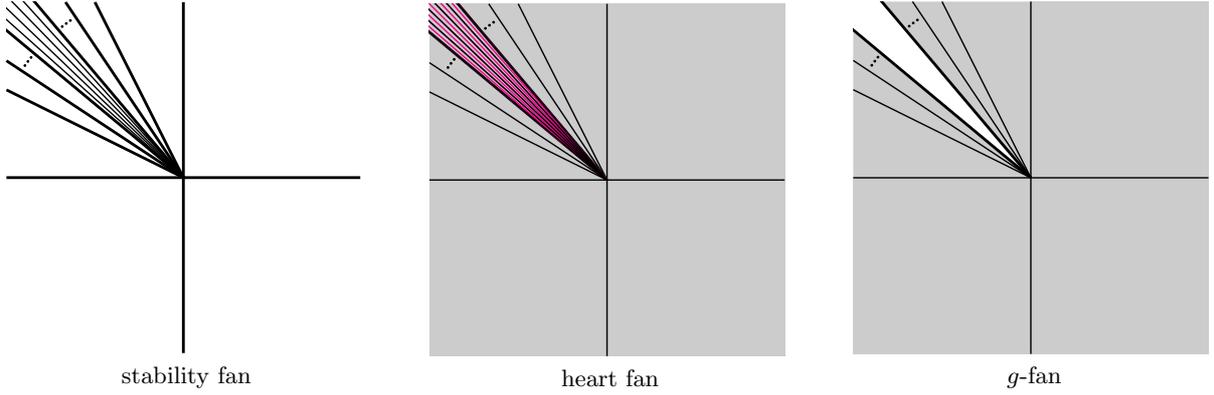
\newcommand{\sstex}[2]{%
  \parbox{0.30\textwidth}{%
    \resizebox{\linewidth}{!}{\input{#1}} \\
    {\footnotesize\centering #2\par}
  }
}

\noindent
\begin{center}
  \parbox{0.30\textwidth}{%
    \resizebox{\linewidth}{!}{\input{diag_3-Kronecker__stability_fan.tex}} \\
    {\footnotesize\centering stability fan\par}
  }

\hfill
  \parbox{0.30\textwidth}{%
    \resizebox{\linewidth}{!}{\input{diag_3-Kronecker__heart_fan.tex}} \\
    {\footnotesize\centering heart fan\par}
  }

\hfill
  \parbox{0.30\textwidth}{%
    \resizebox{\linewidth}{!}{\input{diag_3-Kronecker__g-fan.tex}} \\
    {\footnotesize\centering $g$-fan\par}
  }

\end{center}
\end{example}

\begin{remark}
\label{rem:non-rational stability fan}
The previous example may suggest that the stability fan is the heart fan with all full cones and all non-rational cones removed. This is not true in general: if $Q = Q' \amalg Q''$ is the disjoint union of the $A_2$ and the 3-Kronecker quivers then its heart fan is the product of the heart fans of the components,
  $\fan(\kk Q) = \fan(\kk Q') \times \fan(\kk Q'')$.
The product of a non-rational ray from $\fan(\kk Q'')$ and any ray from $\fan(\kk Q')$ is a non-rational cone with non-zero distinguished subcategory $\SS$, thus is in the stability fan.
\end{remark}

\section{$g$-fans}
\label{sec:g-fans}

\noindent
Let $A$ be a finite-dimensional algebra over an algebraically closed field $\kk$. We write $\fan(A)$ for the heart fan of the category of finite-dimensional $A$-modules $\mod{A}$ over $\LL \coloneqq K(A) \cong \Z^n$. The main result of this section is that the $g$-fan of $A$ is the subfan of $\fan(A)$ consisting of the full cones and their faces.
In order to define the $g$-fan of $A$ we first recall some silting theory.

\subsection{Projective-simple duality}
We may assume that $A$ is basic, \ie $A \cong P(1)\oplus\cdots\oplus P(n)$, where the $P(i)$ are pairwise non-isomorphic indecomposable projective $A$-modules.
We consider two Hom-finite triangulated categories associated to $A$ together with their Grothendieck groups:
\begin{align*}
       \DD^p &\coloneqq \Kb(\proj{A})             && \text{the perfect derived category of $A$,} \\[-1ex]
       \DD^b &\coloneqq \Db(\mod{A})              && \text{the bounded derived category of $A$,} \\[-1ex]
        K(A) &\coloneqq K(\DD^b) = K(\mod{A})     && \text{the Grothendieck group of $A$,}       \\[-1ex]
  \Ksplit{A} &\coloneqq K(\DD^p) = K(\proj{A})    && \text{the split Grothendieck group of $A$}.
\end{align*}
The classes $[P(1)],\ldots,[P(n)]$ form a basis of $\Ksplit{A}$ and the classes of simple modules $[S(1)],\ldots,[S(n)]$ form a basis of $K(A)$. These bases are dual under the non-degenerate pairing
\[ \Ksplit{A} \times K(A) \to \Z, \quad (P,M) \mapsto \pairing{P}{M} = \dim \Homm{A}{P}{M} . \]
Having set $\LL \coloneqq K(A)$, we identify $\LLdual = \Ksplit{A}$.

\subsection{Two-term (pre)silting objects}

Consider the full subcategory of two-term complexes of projective $A$-modules, $\Ktwo \coloneqq \cat{K}^{[-1,0]}(\proj{A}) = \add{A} * \add{A}[1] \subset \DD^p$. It is an extriangulated category with $K(\CC) = K(\DD^p)$.
A complex $P\cp \in \Ktwo$ is called
\begin{itemize}
\item (two-term) \defn{presilting} if $\Homm{\DD^p}{P\cp}{P\cp[{>}0]} = 0$;
\item (two-term) \defn{silting} if it is presilting and the thick subcategory generated by $P\cp$ is $\DD^p$.
\end{itemize}
While we don't use the adjective `two-term', we stress that we are working in $\CC$ throughout, so that our complexes have non-zero terms in homological degrees $1$ and $0$ only.

\subsection{$\bm{g}$-vectors, $\bm{g}$-cones and the $\bm{g}$-fan}

The \defn{$\bm{g}$-vector} of $P\cp \in \Ktwo$ is its class $[P\cp] \in K(\CC)$; writing $P\cp = (P_1 \to P_0)$ and decomposing $P_0 = \bigoplus_{i=1}^n P(i)^{a_i}$ and $P_1 = \bigoplus_{i=1}^n P(i)^{b_i}$, this is the classical $g^{P\cp} = (a_1 - b_1, \ldots, a_n - b_n) \in \Z^n \cong K(\CC)$, identified using the basis $P(1),\ldots,P(n)$.

For $P\cp \in \CC$ presilting with indecomposable summands $Q^1\cp, \ldots, Q^t\cp$, the $g$-vectors $[Q^1\cp], \ldots, [Q^t\cp]$ are linearly independent by \cite{Demonet-Iyama-Jasso} and generate the \defn{$g$-cone} 
\[
C(P\cp) = \{ a_1 [Q^1\cp] + \cdots + a_t [Q^t\cp] \mid a_i \geq 0 \} \subset \VSdual = K(\CC)\otimes\R.
\]
They form a basis if $t=n$, \ie if $P\cp$ is silting.

The \defn{$\bm{g}$-fan}
$\fan_A\gfan \coloneqq \{ C(P\cp) \mid P\cp \in \Ktwo \text{ is presilting}\}$
of a finite-dimensional algebra $A$ was introduced in \cite[\S6]{Demonet-Iyama-Jasso}.
For comparison with the heart fan, we note that any fan $\fan$ has a subfan $\fan\full$ generated by its full cones and for the heart fan $\fan(A)$ this is
\[
  \fan(A)\full = \{ \kappa \in \fan(A) \mid \exists \sigma\in\fan(A) \text{ full}, \kappa \faceof \sigma \}
  = \{ C(\KK/\SS) \in \fan(A) \mid \KK \text{ algebraic} \} ,
\]
where the identification of full cones with algebraic hearts follows from Proposition~\ref{prop:full cones} together with Lemma~\ref{lem:ghostfree tilting}.
The $g$-fan coincides with the full subfan of the heart fan $\fan(A)$:

\begin{theorem}
\label{thm:g-fan}
The $g$-fan of the finite-dimensional algebra $A$ is $\fan_A\gfan = \fan(A)\full$.
\end{theorem}

\begin{remark}
\label{rem:finite g-fan}
In Theorem~\ref{thm:heart fan}, we showed that the heart fan of any length category is complete. Thus the heart fan of $\mod{A}$ is a natural completion of the $g$-fan.
The following finiteness conditions are all equivalent for a finite-dimensional algebra $A$:
\begin{enumerate}[(i)]
\item \label{f:heart fan finite} the heart fan is finite;
\item \label{f:g-fan finite}     the $g$-fan is finite, \ie $A$ is $\tau$-tilting finite (the set of silting complexes in $\Ktwo$ is finite);
\item \label{f:g-fan complete}   the $g$-fan is complete;
\item \label{f:heart is g-fan}   heart fan and $g$-fan coincide;
\item \label{f:algebraic hearts} all two-term hearts relative to $\HH$ are algebraic;
\item \label{f:finite hearts}    there are finitely many two-term hearts relative to $\HH$.
\end{enumerate}
Theorem~\ref{thm:g-fan} shows
\ref*{f:heart fan finite} $\timplies$ \ref*{f:g-fan finite} and
\ref*{f:g-fan complete} $\timplies$ \ref*{f:heart is g-fan} $\tiff$ \ref*{f:algebraic hearts},
where the last equivalence uses that full heart cones correspond to algebraic hearts by Corollary~\ref{cor:full-alg bijection}.
\ref*{f:g-fan finite} $\tiff$ \ref*{f:g-fan complete} is \cite[Thm.~4.7]{Asai21},
\ref*{f:algebraic hearts} $\timplies$ \ref*{f:finite hearts} follows from \cite[Thms.~2.7 \& 3.2]{Adachi-Iyama-Reiten} and \cite[Thm.~6.1]{Koenig-Yang} using 
  \cite[Thm.~3.8]{Demonet-Iyama-Jasso},
and
\ref*{f:finite hearts} $\implies$ \ref*{f:heart fan finite} is just the heart fan definition.

Of these conditions, \ref*{f:heart fan finite}, \ref*{f:algebraic hearts}, \ref*{f:finite hearts} don't refer to the algebra and so make sense for arbitrary abelian categories.
For algebraic abelian categories, they are indeed equivalent by Theorem~\ref{thm:finite heart fan}.
\end{remark}

\begin{remark}
Lutz Hille \cite{Hille} defined a fan using classical (partial) tilting modules over a finite-dimensional algebra.
His fan comprises a collection of subcones of the effective cone of the module category, and therefore sits in the positive orthant. 
Despite a similar homological approach, Hille's construction is unrelated to ours and the $g$-fan.
\end{remark}

To prove the proposition, we need to describe the connection of silting theory to Serre subcategories via simple-minded collections and $c$-vectors.

\subsection{Koenig-Yang correspondences and Serre subcategories}

In \cite{Koenig-Yang}, Steffen Koenig and Dong Yang obtain correspondences between algebraic hearts $\KK$ in $\DD^b$ and silting subcategories of $\DD^p$ which restrict to two-term versions by \cite{Bruestle-Yang}. Recall that $\{X_1,\ldots,X_r\}$ is a \defn{simple-minded collection} in a $\kk$-linear triangulated category if (a) $\Home{<0}{X_i}{X_j} = 0$ for all $i, j$ and (b) $\Hom{X_i}{X_i} = \kk$, $\Hom{X_i}{X_j} = 0$ for $i\neq j$ and (c) the collection generates the triangulated category; see \cite[\S3.2]{Koenig-Yang}.

\begin{proposition}[{\cite[Thm.~6.1]{Koenig-Yang} \& \cite[Cor.~4.3]{Bruestle-Yang}}]
\label{prop:Koenig-Yang}
Let $A$ be a finite-dimensional algebra and $\HH = \mod{A}$. There are bijections between the following sets:
\begin{enumerate}
\item silting objects in $\Ktwo$, up to additive equivalence;
\item algebraic hearts $\KK$ with $\HH[1] \geq \KK \geq \HH$;
\item simple-minded collections in $\HH[1] * \HH$. 
\end{enumerate}
\end{proposition}

The map from (2) to (3) sends $\KK \mapsto \{\text{simple objects of } \KK\}$; the inverse map from (3) to (2) sends a simple-minded collection $\{X_1,\ldots, X_n\} \mapsto \clext{X_1,\ldots,X_n}$, its extension closure in $\DD^b$.
The objects in (1) and (3) obey a \emph{generalised projective-simple duality}: let $P\cp = \bigoplus_{i=1}^n Q^i\cp$ be a silting object and $\{X_1,\ldots,X_n\}$ the corresponding simple-minded collection, then
  $\Hom{Q^i\cp}{X_i} = \kk$ and $\Hom{Q^i\cp}{X_j} = 0$ for $i\neq j$.
Moreover, $\KK \simeq \mod{\End{}{P\cp}}$.
By Proposition~\ref{prop:Koenig-Yang} and \cite[Prop. 5.3]{Geigle-Lenzing}, there are bijections
\[
      \{ \text{direct summands of } P\cp \} 
 \bij \{ \text{subsets of } \{X_1, \ldots, X_n\} \}
 \bij \{ \text{Serre subcategories of } \KK \}.
\]

\subsection{$\bm{c}$-vectors}

Following \cite[p.~5035]{Asai20}, if $\{X_1, \ldots, X_n \}$ is a simple-minded collection, then the \defn{$\bm{c}$-vector} of $X_i$ is its class $[X_i] \in K(A)$. Expressed in the basis of $K(A)$ of simple modules, this gives the classical $c^{X_i} \in \Z^n$.
Since $\{X_1, \ldots, X_n\}$ is a complete list of non-isomorphic simple objects for an algebraic heart in $\DD^b$, the $c$-vectors $\{[X_1],\ldots,[X_n]\}$ form a basis of $\VS$.
By generalised simple-projective duality $\pairing{[Q^i\cp]}{[X_j]} = \delta_{ij}$, and the $c$-vector basis is dual to the $g$-vector basis of $\LLdual$ coming from the silting object $P\cp = \bigoplus_{i=1}^n Q^i\cp$ corresponding to $\{X_1, \ldots, X_n\}$.

\begin{remark}[Heart cofans as $\bm{c}$-vector cofans]
\label{rem:virtual g-fan}
Let $\HH$ and $\HH[1] \geq \KK \geq \HH$ be algebraic hearts. The simple objects in $\KK$ form a simple-minded collection in $\Db(\HH)$ and their classes generate the integral cone $E(\KK)$; hence the heart cofan $\cofanZ(\HH)$ or its real version $\cofanR(\HH)$ can be thought of as convex-geometric objects corresponding to $c$-vectors. This is possible even if the projective objects necessary to define $g$-vectors don't exist.
Therefore, the heart cofan $\cofanZ(\HH)$ may be considered as the `$c$-cofan' and the corresponding dual fan $\fan(\HH)$ as a `virtual $g$-fan'.
\end{remark}

We now assemble the ingredients from silting theory to prove Theorem~\ref{thm:g-fan}.

\begin{proof}[Proof of Theorem~\ref{thm:g-fan}]
Let $\KK$ be an algebraic heart with $\HH[1] \geq \KK \geq \HH$, write $P\cp = \bigoplus_{i=1}^n Q^i\cp$ for the corresponding two-term silting complex in $\Ktwo$, and $\{X_1,\ldots,X_n\}$ for the corresponding two-term simple-minded collection in $\HH[1] * \HH$ under the bijections in Proposition~\ref{prop:Koenig-Yang}. Then $C(\KK) = C(P\cp)$ follows from the definitions and the duality $\pairing{[Q^i\cp]}{[X_j]} = \delta_{ij}$.
%

Now, suppose $\SS \subseteq \KK$ is a Serre subcategory. By \cite[Prop.~5.3]{Geigle-Lenzing}, we can re-order the two-term simple-minded collection $\{X_1,\ldots,X_n\}$ so that $\SS = \clext{X_1,\ldots,X_t}$ for some $1 \leq t \leq n$. 
We order the summands of the corresponding two-term silting complex $P\cp$ accordingly and write $Q\cp = \bigoplus_{i=t+1}^n Q^i\cp$.
It follows immediately from the generalised projective-simple duality that $C(\KK/\SS) = C(Q\cp)$ using $C(\KK/\SS) = C(\KK) \cap E(\SS)\orth$ from Lemma~\ref{lem:dual faces are heart cones}.
%
\end{proof}

\appendix
 
\section{Homological algebra}
\label{appA:homological algebra} 

\noindent
A good source emphasising the similarity of abelian and triangulated categories is \cite[Ch.~I]{Beligiannis-Reiten}.
Categories are assumed to be essentially small, \ie having a set of objects up to isomorphism. Subcategories are assumed to be full and strict, \ie closed under isomorphisms.
The $n$-fold shift (or: translation or suspension) functor of triangulated categories is denoted by $[n]$.

\subsection{Abelian categories.}
\label{appA:abelian}
An abelian category $\HH$ is called \defn{(of finite) length} if each object of $\HH$ admits a finite filtration with simple factors. Equivalently, $\HH$ is noetherian (all ascending chains of subobjects stabilise) and artinian (all descending chains stabilise). If in addition $\HH$ has finitely many isomorphism classes of simple objects, it is called \defn{algebraic}.

A non-empty full subcategory $\UU \subseteq \HH$  is \defn{wide} if it is closed under kernels, cokernels and extensions, and is \defn{Serre} if it is closed under subobjects, quotients and extensions. The quotient category $\HH/\UU$ by a Serre subcategory $\UU$  is abelian.
Denote  the sets of Serre and wide subcategories by $\Serre{\HH} \subseteq \Wide{\HH}$ respectively.
In Definition~\ref{def:face subcategory}, we introduce $\faceSerre{\HH}$ as the subset of face subcategories.

\subsection{Triangulated categories.}
Let $\DD$ be a triangulated category.
A subcategory $\UU$ of $\DD$ is called \defn{thick} if it is a triangulated subcategory of $\DD$ and is closed under direct summands. We denote the set of thick subcategories by $\Thick{\DD}$. If $S$ is any collection of objects of $\DD$ then we denote by $\thick{}{S}$ the smallest thick subcategory of $\DD$ containing $S$.

\subsection{Extension closure.}
Let $\CC$ be an abelian or triangulated category, so we can speak of extensions in $\CC$.
The \defn{ordered extension closure} $\MM_1 * \MM_2$ of subcategories $\MM_1,\MM_2 \subseteq \CC$ is the subcategory of all extensions of objects of $\MM_2$ by objects of $\MM_1$.
E.g.\ for $\CC$ triangulated $\MM_1 * \MM_2$ consists of $c\in\CC$ sitting in exact triangles $m_1 \to c \to m_2 \to m_1[1]$ with $m_1\in \MM_1, m_2\in \MM_2$.
%

The \defn{extension closure} of a subcategory $\MM \subseteq \CC$ is $\clext{M} \coloneqq \MM * \MM$. 
For $\HH\subset\CC$ a heart in a triangulated category and $\MM\subseteq\HH$, the two notions of $\clext{\MM}$ agree.

\subsection{Torsion pairs.}
\label{appA:torsion pairs}
Two full subcategories $(\TT,\FF)$ of an abelian category $\HH$ form a \defn{torsion pair} if (a) $\Hom{\TT}{\FF}=0$ and (b) $\HH = \TT * \FF$. In this case, $\TT$ is called the \defn{torsion class} and $\FF$ the \defn{torsionfree class} of the pair.
Torsion classes are closed under quotients and extensions; torsionfree classes are closed under subobjects and extensions.
If $(\TT,\FF)$ is a torsion pair then $\FF = \TT\orth \coloneqq \{ h \in \HH \mid \Hom{\TT}{h} = 0\}$ and $\TT= {}\orth\FF$; \eg $\FF \subseteq \TT\orth$ by (a) and $\FF \supseteq \TT\orth$ by (b).

If $\HH$ is noetherian and $\TT \subseteq \HH$ is a subcategory closed under quotients and extensions then $(\TT,\TT\orth)$ is a torsion pair. Dually, if $\HH$ is artinian and $\FF \subseteq \HH$ closed under subobjects and extensions then $({}\orth\FF,\FF)$ is a torsion pair. An example is $\FF = \clext{s}$ for a simple object $s\in\HH$.


\subsection{T-structures.}
\label{appA:t-structures}
Two subcategories $(\TT,\FF)$, called \defn{aisle} and \defn{co-aisle} of a triangulated category $\DD$ form a \defn{t-structure} if (a) $\Hom{\TT}{\FF}=0$ and (b) $\CC = \TT * \FF$ and (c) $\TT[1] \subseteq \TT$. The t-structure is \defn{bounded} if moreover (d) $\bigcup_{n\in\Z} \TT[n] = \DD = \bigcup_{n\in \Z} \FF[n]$.
Any t-structure $(\TT,\FF)$ on $\DD$, bounded or not, induces an abelian subcategory of $\DD$ as $\HH \coloneqq \TT\cap\FF[1]$ called the \defn{heart}. If the t-structure is bounded then $\TT$ and $\FF$ can be reconstructed from the heart $\HH$ as $\TT = \clext{\HH[\mathit{\geq}0]}$ and $\FF = \clext{\HH[\mathit{<}0]}$. In this article, {\em heart} always means the {\em heart of a bounded t-structure}.

Setting $(\TT,\FF) \geq (\TT',\FF') \iff \FF \supseteq \FF'$ defines a partial order on the bounded t-structures of a triangulated category $\DD$.
Equivalently,
  $\HH \geq \HH' \iff \clext{\HH[{\leq}0]} \supseteq \clext{\HH'[{\leq}0]}$
gives a partial order on hearts in $\DD$.
The convention we follow here makes $\HH[1] \geq \HH$.

A heart $\KK$ is called \defn{two-term relative to $\HH$} if $\HH[1] \geq \KK \geq \HH$.
The next fact is well known.

\begin{lemma}
For bounded hearts $\HH$ and $\KK$ in $\DD$, the following conditions are equivalent:
\begin{enumerate}[label=(\roman*)]
\item $\HH[1] \geq \KK \geq \HH$;
\item $\KK \subseteq \HH[1] * \HH$;
\item $\Hom{\HH[{\geq} 2]}{\KK} = 0$ and $\Hom{\KK[{\geq} 1]}{\HH} = 0$.
\end{enumerate}
\end{lemma}

\begin{proof}
  $(i) \timplies (ii)$: The inequalities say $\KK \subseteq \clext{\KK[{\leq} 0]} \subseteq \clext{\HH[{\leq} 1]}$ and $\KK \subseteq \clext{\KK[{\geq} 0]} \subseteq \clext{\HH[{\geq} 0]}$, hence $\KK \subseteq \clext{\HH[1],\HH}$. To show $\HH * \HH[1] \subseteq \HH[1] * \HH$, take $d \in \HH * \HH[1]$, \ie $d$ is the cone of a morphism $\alpha \colon h_1 \to h_0$ in $\HH$. The octahedral axiom for the composition $h_1 \to \im(\alpha) \to h_0$ produces the exact triangle $\ker(\alpha)[1] \to d \to \coker(\alpha) \to \ker(\alpha)[2]$, exhibiting $d \in \HH[1] * \HH$.

$(ii) \timplies (iii)$: $\HH$ is a heart, so $\Hom{\HH[{\geq} 1]}{\HH} = 0$; with $\KK \subseteq \HH[1] * \HH$ this gives $(iii)$.

$(iii) \timplies (i)$: We have $\Hom{\HH[{\geq} 2]}{\KK[{\leq} 0]} = 0$, so that $\clext{\KK[{\leq} 0]} \subseteq (\HH[{\geq} 2])\orth = \clext{\HH[{\leq} 1]}$, where the last equality holds because $\HH$ is a bounded heart. Hence, $\HH[1] \geq \KK$. 
Similarly, $\Hom{\KK[{\geq} 1]}{\HH[{\leq} 0]} = 0$ gives $\clext{\HH[{\leq} 0]} \subseteq (\KK[{\geq} 1])\orth = \clext{\KK[{\leq} 0]}$ and $\KK \geq \HH$.
\end{proof}

\subsection{Tilting.}
\label{appA:tilting}
Let $\DD$ be a triangulated category, $\HH$ a heart in $\DD$ and $(\TT,\FF)$ a torsion pair in the abelian category $\HH$, so that $\HH = \TT * \FF$. From this setup we get two new hearts of $\DD$ as
\[ \FF[1] * \TT \text{, the \defn{positive tilt} of $\HH$}; \qquad
   \FF * \TT[-1] \text{, the \defn{negative tilt} of $\HH$}. \]
The terminology is justified by $\HH[1] \geq \FF[1] * \TT \geq \HH \geq \FF * \TT[-1] \geq \HH[-1]$.
Positive/negative tilts are also called left/right.
This fact is crucial for us, see \cite[Lemma~1.1.2]{Polishchuk}:
If $\KK$ and $\HH$ are two hearts in a triangulated category $\DD$ with $\HH[1] \geq \KK \geq \HH$, then there is a unique torsion pair $(\TT,\FF)$ on $\HH = \TT * \FF$ with $\KK = \FF[1] * \TT$; it is given by $\TT = \HH \cap \KK$ and $\FF = \HH \cap \KK[-1]$.

%

The \defn{positive irreducible simple tilt} of an artinian category $\HH$ at a simple object $s\in\HH$ is given by the torsion pair $({}\orth s, \clext{s})$.
%

\input{bibliography_fan.tex}

\bigskip
{\small\texttt{
\noindent
\begin{tabular}{lll}
  \textit{\textrm{Contact:}}
  & nathan.broomhead@plymouth.ac.uk, & d.pauksztello@lancaster.ac.uk, \\
  & david.ploog@uis.no,              & jonathan.woolf@liverpool.ac.uk
\end{tabular}}}%

\end{document}

%% file: header-fan.tex
\DeclareMathAlphabet{\mathsfit}{T1}{\sfdefault}{\mddefault}{\sldefault}

\usepackage[shortlabels]{enumitem}
\usepackage[small]{caption} 
\usepackage{xcolor}
\usepackage{adjustbox}

\usepackage{amsthm,amssymb,amsmath}
\usepackage{mathtools} 
\usepackage{stmaryrd}  
\usepackage{bm}        

\usepackage{tikz,tikz-cd,tikz-3dplot}

\usepackage{booktabs}
\newcommand{\mc}{\multicolumn}

\usepackage{graphicx}

\usepackage{hyperref}
  \hypersetup{
      colorlinks,
      linkcolor={red!50!black},
      citecolor={blue!50!black},
      urlcolor={blue!80!black}
  }
\newcommand{\harxiv}[1]{\href{http://arxiv.org/abs/#1}{\texttt{arXiv:#1}}}

\usepackage{fullpage}

\newcommand{\defn}[1]{\emph{\textbf{#1}}}

\newtheoremstyle{boldhead}
  {}
  {}
  {}
  {}
  {\bfseries}
  {.}
  { }
  {\thmname{#1}\thmnumber{ #2}\thmnote{ (#3)}}

\theoremstyle{plain}
\newtheorem{theorem}{Theorem}[section]

\newtheorem{lemma}[theorem]{Lemma}
\newtheorem{corollary}[theorem]{Corollary}
\newtheorem{proposition}[theorem]{Proposition}
\newtheorem{thmdef}[theorem]{Theorem-Definition}
\newtheorem{lemdef}[theorem]{Lemma-Definition}
\newtheorem{introtheorem}{Theorem}

\theoremstyle{boldhead}

\newtheorem{remark}[theorem]{Remark}
\newtheorem{example}[theorem]{Example}

\newtheorem{definition}[theorem]{Definition}

\newcommand{\proofstep}[1]{\smallskip\noindent\emph{#1:}}

\newcommand{\bib}[6]{{\bibitem[#1]{#2} #3: {\emph{#4},} #5#6.}} 
\newcommand{\bibno}[1]{} 

\newcommand{\ie}{i.e.\ }
\newcommand{\eg}{e.g.\ }

\newcommand{\ssi}{\mathsf{ss}}                 
\newcommand{\sst}{{}^\mathsf{ss}}              
\newcommand{\HHsst}[1]{\cat{H}^\mathsf{\mkern1mu ss}_{#1}} 
\newcommand{\fanst}{\fan^\mathsf{\mkern-3mu s}}            

\newcommand{\gfan}{^g}                         
\newcommand{\full}{^{\mathrm{full}}}           

\newcommand{\rslash}{\!/}                      
\newcommand{\PD}[1]{\Phi_{#1}}                 

\DeclareMathOperator{\conv}{{\mathrm{conv}}}   
\DeclareMathOperator{\codim}{{\mathrm{codim}}} 
\DeclareMathOperator{\id}{{\mathrm{id}}}       

\DeclarePairedDelimiterX{\norm}[1]{\lVert}{\rVert}{#1} 
\DeclarePairedDelimiterX{\abss}[1]{\lvert}{\rvert}{#1} 
\DeclarePairedDelimiterX{\supp}[1]{\lvert}{\rvert}{#1} 

\newcommand{\kk}{{\mathbf{k}}}                 
\newcommand{\Pic}{\mathrm{Pic}}
\newcommand{\cp}{_\bullet}                     
\newcommand{\Ksplit}[1]{K^{\text{split}}(#1)}  

\newcommand{\timplies}{\Longrightarrow}        

\newcommand{\tiff}{\Longleftrightarrow}

\newcommand{\vdual}{^*}                        
\newcommand{\rdual}   {\mathbb{D}}             
\newcommand{\ldual}[1]{\mathbb{L}_{#1}}        
\newcommand{\cdual}[1]{\rdual(#1)}             
\newcommand{\orth}{^\perp}                     
\newcommand{\sat}[1]{\mathrm{sat}(#1)}         

\renewcommand{\ll}{\lambda}                    
\newcommand{\LL}{\Lambda}                      
\newcommand{\LLdual}{\LL\vdual}                
\newcommand{\VS}{\Lambda_\R}                   
\newcommand{\VSdual}{{\Lambda_\R\vdual}}       
\newcommand{\faceof}    {\leq}                 

\newcommand{\pfaceof}   {<}                    
\newcommand{\cofaceof}  {\preccurlyeq}         

\newcommand{\zerocone}{\{0\}}                  
 
\newcommand{\cones}[1]{\mathrm{Cones}(#1)}     
\newcommand{\Cones}[1]{\mathrm{Cones}(#1)}     

\newcommand{\Faces}  [1]{\catt{Faces}(#1)}     
\newcommand{\DFaces} [1]{\catt{DualFaces}(#1)} 
\newcommand{\ExFaces}[1]{\catt{ExFaces}(#1)}   
\newcommand{\Cofaces}[1]{\catt{Cofaces}(#1)}   
\newcommand{\CEFaces}[1]{\catt{CoExFaces}(#1)} 

\newcommand{\fCones}[1]{\catt{Cones}_f(#1)}    
\newcommand{\cCones}[1]{\catt{Cones}_c(#1)}    

\newcommand{\opp}{{}^{\mathrm{op}}}            

\newcommand{\mf}[1] {\mathrm{mf}(#1)}          
\newcommand{\mface}[1]{#1_{\mathrm{min}}}      

\newcommand{\Serre}[1]{\catt{Serre}(#1)}             
\newcommand{\faceSerre}[1]{\catt{Serre}_{\LL\!}(#1)} 
\newcommand{\kerc}[2]{\catt{ker}_{#1}(#2)}           

\newcommand{\Hearts}   [1]{\catt{Hearts}(#1)}            
\newcommand{\FacePairs}[1]{\catt{FacePairs}_{\LL\!}(#1)} 
\newcommand{\KerPairs} [1]{\catt{KerPairs}_{\LL\!}(#1)}  

\newcommand{\Reversed}{{\ltimes}}
\newcommand{\Mixed}{{\Join}}
\newcommand{\reversedHeart}[1]{\catt{coh}^{\Reversed \!}(#1)}
\newcommand{\mixedHeart}   [1]{\catt{coh}^{\Mixed \!}(#1)}

\newcommand{\rk}[1]{\mathrm{rk}\, #1} 
\newcommand{\rank}{\mathrm{rk}}

\newcommand{\relintC}[1]{C^\circ(#1)}
\newcommand{\clE}[1]{\mkern4mu\overline{\mkern-4muE}(#1)} 

\DeclareMathAlphabet{\mathpzc}{OT1}{pzc}{m}{it} 

\newcommand{\Kronecker}{\bullet\rightrightarrows\bullet}
\newcommand{\cofan}{\nabla}         
\newcommand{\cofanZ}{\nabla_{\!\Z}} 
\newcommand{\cofanR}{\nabla_{\!\R}} 
\newcommand{\fan}{\Delta}           
\newcommand{\Fsup}{\Delta}          

\newcommand{\Hom}[2]{{\mathrm{Hom}}(#1,#2)}
\newcommand{\Home}[3]{{\mathrm{Hom}}^{#1}(#2,#3)}
\newcommand{\Homm}[3]{{\mathrm{Hom}}_{#1}(#2,#3)}

\newcommand{\End}[2]{{\mathrm{End}}_{#1}(#2)}
\newcommand{\Aut}[1]{\mathrm{Aut}(#1)}
\newcommand{\Autt}[2]{\mathrm{Aut}_{#1}(#2)}

\renewcommand{\top}{\mathrm{top}\,}
\DeclareMathOperator{\soc}{\mathrm{soc}}

\newcommand{\onto}{\twoheadrightarrow}
\newcommand{\into}{\hookrightarrow}
\newcommand{\too}{\longrightarrow}
\newcommand{\isom} {\xrightarrow{\raisebox{-0.5ex}[0ex][0ex]{$\sim$}}}
\newcommand{\lisom}{\overset{\sim}{\too}}

\newcommand{\bij}{\stackrel{1-1}{\longleftrightarrow}}

\newcommand{\pairing}[2]{\chi(#1,#2)} 

\newcommand{\im}{\mathrm{im}}         
\newcommand{\coker}{\mathrm{coker}} 

\newcommand{\catt}[1]{\mathsf{#1}}  
\newcommand{\cat}[1]{\mathsfit{#1}} 

\newcommand{\CC}{\cat{C}}
\newcommand{\DD}{\cat{D}}
\newcommand{\FF}{\cat{F}}
\newcommand{\HH}{\cat{H}}
\newcommand{\KK}{\cat{K}}
\newcommand{\MM}{\cat{M}}
\newcommand{\NN}{\cat{N}}

\renewcommand{\SS}{\cat{S}}
\newcommand{\TT}{\cat{T}}
\newcommand{\UU}{\cat{U}}

\newcommand{\TTT} [1]{\cat{T}_{\! #1}}
\newcommand{\FFF} [1]{\cat{F}_{#1}}
\newcommand{\TTTT}[1]{\smash{\mkern4mu \overline{\mkern-4mu\cat{T}}_{\mkern-6mu #1}}} 
\newcommand{\FFFF}[1]{\smash{\mkern4mu \overline{\mkern-4mu\cat{F}}_{\mkern-4mu #1}}}

\DeclareMathOperator{\coh}{\mathsf{coh}}           
\newcommand{\Db}{\catt{D}^b}                       
\newcommand{\Kb}{\catt{K}^b}                       
\newcommand{\thick}[2]{\catt{thick}_{#1}(#2)}      
\newcommand{\Thick}[1]{\catt{Thick}{(#1)}}         
\newcommand{\Wide} [1]{\catt{Wide}{(#1)}}          
\newcommand{\clext}[1]{{\langle #1 \rangle}}       
\newcommand{\modsf}   {\catt{mod}}                 
\renewcommand{\mod}[1]{\modsf(#1)}                 
\newcommand{\nilrep}  {\catt{nilrep}}              

\newcommand{\add}  [1]{\catt{add}(#1)}
\newcommand{\proj} [1]{\catt{proj}(#1)}
\newcommand{\sub}  [1]{\catt{sub}(#1)}

\newcommand{\Ktwo}{\CC}                            
\newcommand{\tightverticaldots}{ \rotatebox{90}{...} \\[-1ex] } 
\newcommand{\Groker}[1]{\catt{N}_{#1}}             

\newcommand{\stab}[1]{\mathrm{Stab}(#1)}           

\renewcommand{\phi}{\varphi}
\renewcommand{\epsilon}{\varepsilon}
\renewcommand{\emptyset}{\varnothing}
\renewcommand{\setminus}{\backslash}

\newcommand{\N}{\mathbb{N}}
\newcommand{\Z}{\mathbb{Z}}
\newcommand{\Q}{\mathbb{Q}}
\newcommand{\R}{\mathbb{R}}
\newcommand{\IR}{{\,\mathbb{R}}}
\newcommand{\Rpos}{\mathbb{R}_{\geq0}}
\newcommand{\C}{\mathbb{C}}
\newcommand{\IH}{\mathbb{H}} 
\newcommand{\PP}{\mathbb{P}} 

\newcommand{\cC}{\mathcal{C}}
\newcommand{\cD}{\mathcal{D}}

\newcommand{\cO}{\mathcal{O}}

%% file: intro_Kronecker.tex
\begin{tikzpicture}[scale=3.5]
  \fill[violet!10] (0,0) -- (-1,1) -- (-1,0) -- (-1,-1) -- cycle;  
  \fill[black!10] (0,0) rectangle (1,1);
  \fill[black!15] (0,0) -- (0,1) -- (-0.5,1) -- cycle;
  \fill[black!20] (0,0) -- (-0.5,1) -- (-0.666,1) -- cycle;
  \fill[violet!35] (0,0) -- (-1,0.5) -- (-1,0.666) -- cycle;
  \fill[violet!30] (0,0) -- (-1,0) -- (-1,0.5) -- cycle;  
  \fill[violet!25] (0,0) rectangle (-1,-1);
  \fill[black!17] (0,0) rectangle (1,-1);  
  \draw[thick] (-1,0) -- (1,0);
  \draw[thick] (0,-1) -- (0,1);
  \draw[thick] (0,0) -- (-0.5,1);
  \draw[thick] (0,0) -- (-0.666,1);
  \draw[ultra thick,\raycolor] (0,0) -- (-0.98,0.98);
  \node[\raycolor] at (-1.15,0.93) {$C(\HH_2)$};
  \draw[thick] (0,0) -- (-1,0.666);
  \draw[thick] (0,0) -- (-1,0.5);
  \filldraw (125.5:1.1) circle (0.005);
  \filldraw (127.5:1.1) circle (0.005);
  \filldraw (129.5:1.1) circle (0.005);  
  \filldraw (131.5:1.1) circle (0.005); 
  \filldraw (133.5:1.1) circle (0.005);
  \filldraw (136.7:1.1) circle (0.005);
  \filldraw (138.7:1.1) circle (0.005);
  \filldraw (140.7:1.1) circle (0.005);  
  \filldraw (142.7:1.1) circle (0.005); 
  \filldraw (144.7:1.1) circle (0.005);
  \node at  (45:0.7)  {$C(\HH_1)$};
  \node[violet] at (225:0.7)  {$-C(\HH_1) = C(\HH_1[1])$};
\end{tikzpicture}

%% file: intro_projective_line.tex
\begin{tikzpicture}[scale=3.5]
  \fill[violet!10] (0,0) -- (1,0) -- (1,1) -- (0,1) -- cycle;
  \fill[violet!25] (0,0) -- (0,1) -- (1,1) -- cycle;
  \fill[violet!30] (0,0) -- (1,1) -- (1,0.5) -- cycle;
  \fill[violet!35] (0,0) -- (1,0.5) -- (1,0.333) -- cycle;
  \fill[violet!45] (0,0) -- (1,0.333) -- (1,0.25) -- cycle;
  \fill[black!15] (0,0) -- (0,1) -- (-1,) -- cycle;
  \fill[black!20] (0,0) -- (-1,1) -- (-1,0.5) -- cycle;
  \fill[black!25] (0,0) -- (-1,0.5) -- (-1,0.333) -- cycle;
  \fill[black!30] (0,0) -- (-1,0.333) -- (-1,0.25) -- cycle;
  \draw[ultra thick] (-1,0) -- (0,0);
  \draw[ultra thick,\raycolor] (0,0) -- (1,0);
  \draw[thick] (0,0) -- (0,1);
  \draw[thick] (0,0) -- (1,1);
  \draw[thick] (0,0) -- (1,0.5);
  \draw[thick] (0,0) -- (1,0.333);
  \draw[thick] (0,0) -- (1,0.25);
  \draw[thick] (0,0) -- (-1,1);  
  \draw[thick] (0,0) -- (-1,0.5);
  \draw[thick] (0,0) -- (-1,0.333);
  \draw[thick] (0,0) -- (-1,0.25);
  \node at (-0.55,-0.10) {$-C(\HH_2) = C(\HH_2[1])$};
  \node[\raycolor] at (0.85,-0.08) {$C(\HH_2)$};
  \node[violet] at (65:0.8)  {$C(\HH_1[1])$};  
  \filldraw (178.0:0.9) circle (0.005);
  \filldraw (176.0:0.9) circle (0.005);
  \filldraw (174.0:0.9) circle (0.005);
  \filldraw (172.0:0.9) circle (0.005);
  \filldraw (170.0:0.9) circle (0.005);
  \filldraw (168.0:0.9) circle (0.005);
  \filldraw (  2.0:0.9) circle (0.005);
  \filldraw (  4.0:0.9) circle (0.005);
  \filldraw (  6.0:0.9) circle (0.005);
  \filldraw (  8.0:0.9) circle (0.005);
  \filldraw ( 10.0:0.9) circle (0.005);
  \filldraw ( 12.0:0.9) circle (0.005);
\end{tikzpicture}

%% file: diag_heart_fan_A2.tex
\begin{tikzpicture}[scale=3.5]
  \fill[black!10] (0,0) rectangle (1,1);
  \fill[black!15] (0,0) -- (0,1) -- (-1,1) -- cycle;
  \fill[black!20] (0,0) -- (-1,0) -- (-1,1) -- cycle;
  \fill[black!25] (0,0) rectangle (-1,-1);
  \fill[black!17] (0,0) rectangle (1,-1);
  \draw[thick] (-1,0) -- (1,0);
  \draw[thick] (0,-1) -- (0,1);
  \draw[thick] (0,0) -- (-1,1);

  \node at  (45:0.7) {$\HH = \modsf\ \kk(%
       \begin{tikzcd}[cramped, sep=small] 
         \bullet \ar[r] & \bullet
       \end{tikzcd}%
     )$};

  \node at (225:0.7) {$\HH[1]$};
  \node at (-45:0.7) {$\KK_\ssi$};
  \node at (115:0.7) {$\KK_1$};
  \node at (160:0.7) {$\KK_2$};
\end{tikzpicture}

%% file: diag_heart_fan_derived-discrete.tex
\begin{tikzpicture}[scale=3.5]
  \fill[black!10] (0,0) rectangle (1,1);
  \fill[black!15] (0,0) -- (0,1) -- (-1,1) -- cycle;
  \fill[black!20] (0,0) -- (-1,0) -- (-1,1) -- cycle;
  \fill[black!25] (0,0) rectangle (-1,-1);
  \fill[black!15] (0,0) -- (1,0) -- (1,-1) -- cycle; 
  \fill[black!20] (0,0) -- (1,-1) -- (0,-1) -- cycle; 
  \draw[thick] (-1,0) -- (1,0);
  \draw[thick] (0,-1) -- (0,1);
  \draw[thick] (1,-1) -- (-1,1); 
  \node at  (45:0.7) {%
    $\HH = \nilrep(%
      \begin{tikzcd}[cramped, sep=small] 
        \bullet \ar[r, shift left=0.5ex] & \bullet \ar[l, shift left=0.5ex]
      \end{tikzcd})$%
    };
  \node at (225:0.7) {$\HH[1]$};
  \node at (115:0.7) {$\KK_1$};
  \node at (160:0.7) {$\KK_2$};
  \node at (-65:0.7) {$\KK_3$};
  \node at (-20:0.7) {$\KK_4$};
\end{tikzpicture}

%% file: diag_heart_fan_Kronecker.tex
\begin{tikzpicture}[scale=3.5]
  \fill[black!10] (0,0) rectangle (1,1);
  \fill[black!15] (0,0) -- (0,1) -- (-0.5,1) -- cycle;
  \fill[black!20] (0,0) -- (-0.5,1) -- (-0.666,1) -- cycle;
  \fill[black!20] (0,0) -- (-1,0.5) -- (-1,0.666) -- cycle;
  \fill[black!15] (0,0) -- (-1,0) -- (-1,0.5) -- cycle;  
  \fill[black!25] (0,0) rectangle (-1,-1);
  \fill[black!17] (0,0) rectangle (1,-1);  
  \draw[thick] (-1,0) -- (1,0);
  \draw[thick] (0,-1) -- (0,1);
  \draw[thick] (0,0) -- (-0.5,1);
  \draw[thick] (0,0) -- (-0.666,1); 
  \draw[ultra thick] (0,0) -- (-0.97,0.97);
  \draw[thick] (0,0) -- (-1,0.666);
  \draw[thick] (0,0) -- (-1,0.5);
  \filldraw (125.5:1.1) circle (0.005);
  \filldraw (127.5:1.1) circle (0.005);
  \filldraw (129.5:1.1) circle (0.005);  
  \filldraw (131.5:1.1) circle (0.005); 
  \filldraw (133.5:1.1) circle (0.005);
  \filldraw (136.7:1.1) circle (0.005);
  \filldraw (138.7:1.1) circle (0.005);
  \filldraw (140.7:1.1) circle (0.005);  
  \filldraw (142.7:1.1) circle (0.005); 
  \filldraw (144.7:1.1) circle (0.005);

  \node at  (45:0.7)  {$\HH = \modsf\ \kk(%
       \begin{tikzcd}[cramped, sep=small] 
         \bullet \ar[r, shift left=0.5ex] \ar[r, shift right=0.5ex] & \bullet
       \end{tikzcd}%
     )$};

  \node at (-45:0.7)  {$\KK_\ssi$};  
  \node at (225:0.7)  {$\HH[1]$};
  \node at (103:0.7)  {$\KK_1$};
  \node at (170:0.7)  {$\KK'_1$};  
  \node at (138:1.21) {$\KK_\infty$};
\end{tikzpicture}

%% file: diag_heart_fan_3-Kronecker.tex
\begin{tikzpicture}[scale=3.5]
  \fill[black!10] (0,0) rectangle (1,1);
  \fill[black!15] (0,0) -- (0,1) -- (-0.5,1) -- cycle;
  \fill[black!20] (0,0) -- (-0.5,1) -- (-0.666,1) -- cycle;
  \fill[black!20] (0,0) -- (-1,0.5) -- (-1,0.666) -- cycle;
  \fill[black!15] (0,0) -- (-1,0) -- (-1,0.5) -- cycle;  
  \fill[black!25] (0,0) rectangle (-1,-1);
  \fill[black!17] (0,0) rectangle (1,-1);  
  \draw[thick] (-1,0) -- (1,0);
  \draw[thick] (0,-1) -- (0,1);
  \draw[thick] (0,0) -- (-0.5,1);
  \draw[thick] (0,0) -- (-0.666,1);
  \begin{scope}
    \clip (-1,0) rectangle (1,1);
    \draw[ultra thick] (0,0) -- (130:1.5);
    \draw[ultra thick] (0,0) -- (140:105);
    \foreach \X in {0,1,...,10} {
      \draw[thick] (0,0) -- (130+\X:3);
    }
  \end{scope}
  \draw[thick] (0,0) -- (-1,0.666);
  \draw[thick] (0,0) -- (-1,0.5);
  \filldraw (125.5:1.1) circle (0.005);
  \filldraw (127.0:1.1) circle (0.005);
  \filldraw (128.5:1.1) circle (0.005);  
  \filldraw (141.5:1.1) circle (0.005);  
  \filldraw (143.0:1.1) circle (0.005); 
  \filldraw (144.5:1.1) circle (0.005);
  
  \node at  (45:0.7)  {$\HH = \modsf\ \kk(%
       \begin{tikzcd}[cramped, sep=small] 
         \bullet \ar[r, shift left=0.85ex] \ar[r, shift right=0.85ex] \ar[r] & \bullet
       \end{tikzcd}%
     )$};

  \node at (-45:0.7)  {$\KK_\ssi$};  
  \node at (225:0.7)  {$\HH[1]$};
  \node at (103:0.7)  {$\KK_1$};
  \node at (170:0.7)  {$\KK'_1$};  
\end{tikzpicture}

%% file: diag_heart_fan_countable_semisimple.tex
\begin{tikzpicture}[scale=3.5]
  \fill[black!10] (0,0) rectangle (1,1);
  \fill[black!25] (0,0) rectangle (-1,-1);
  
  \fill[black!20] (0,0) -- (-1,1) -- (-0.5,1) -- cycle;
  \fill[black!30] (0,0) -- (-0.5,1) -- (-0.333,1) -- cycle;
  \fill[black!40] (0,0) -- (-0.333,1) -- (-0.25,1) -- cycle;
  \fill[black!50] (0,0) -- (-0.25,1) -- (-0.2,1) -- cycle;
   
  \fill[black!20] (0,0) -- (-1,1) -- (-1,0.5) -- cycle;
  \fill[black!30] (0,0) -- (-1,0.5) -- (-1, 0.333) -- cycle;
  \fill[black!40] (0,0) -- (-1, 0.333) -- (-1,0.25) -- cycle;
  \fill[black!50] (0,0) -- (-1,0.25) -- (-1,0.2) -- cycle;
   
  \fill[black!20] (0,0) -- (1,-1) -- (0.5,-1) -- cycle;
  \fill[black!30] (0,0) -- (0.5,-1) -- (0.333,-1) -- cycle;
  \fill[black!40] (0,0) -- (0.333,-1) -- (0.25,-1) -- cycle;
  \fill[black!50] (0,0) -- (0.25,-1) -- (0.2,-1) -- cycle;
   
  \fill[black!20] (0,0) -- (1,-1) -- (1,-0.5) -- cycle;
  \fill[black!30] (0,0) -- (1,-0.5) -- (1, -0.333) -- cycle;
  \fill[black!40] (0,0) -- (1, -0.333) -- (1,-0.25) -- cycle;
  \fill[black!50] (0,0) -- (1, -0.25) -- (1,-0.2) -- cycle;

  \draw[ultra thick, red, dashed] (-1,0) -- (1,0);
  \draw[ultra thick, red, dashed] (0,-1) -- (0,1);
  
  \draw[thick] (-1,1) -- (1,-1);
  \draw[thick] (-1,0.5) -- (1,-0.5);  
  \draw[thick] (-1,0.333) -- (1,-0.333);
  \draw[thick] (-1,0.25) -- (1,-0.25);
  \draw[thick] (-1,0.2) -- (1,-0.2);
   
  \draw[thick] (-0.5,1) -- (0.5,-1);  
  \draw[thick] (-0.333,1) -- (0.333,-1);
  \draw[thick] (-0.25,1) -- (0.25,-1);
  \draw[thick] (-0.2,1) -- (0.2,-1);

  \filldraw (93:0.8) circle (0.005);
  \filldraw (96:0.8) circle (0.005);
  \filldraw (99:0.8) circle (0.005);  
  
  \filldraw (-87:0.8) circle (0.005);
  \filldraw (-84:0.8) circle (0.005);
  \filldraw (-81:0.8) circle (0.005);  

  \filldraw (171:0.8) circle (0.005);
  \filldraw (174:0.8) circle (0.005);
  \filldraw (177:0.8) circle (0.005);   
  
  \filldraw (-3:0.8) circle (0.005);
  \filldraw (-6:0.8) circle (0.005);
  \filldraw (-9:0.8) circle (0.005);

  \node at  (0.5,0.7) {$\HH = \modsf\ \kk^\Z$};
  \node[centered] at  (0.5,0.44) {$\scriptstyle K(\HH) = \Z^\infty \xrightarrow{~\ll~} \Z^2$};
  \node[centered] at  (0.5,0.3)  {$\scriptstyle \ll(S_i) = (i+1,1), ~ i\geq0$};
  \node[centered] at  (0.5,0.2)  {$\scriptstyle \ll(S_i) = (1,i+1), ~ i\leq0$};  
  \node at (225:0.7) {$\HH[1]$};
   
  \node at (125:1)   {$\KK_{\leq 0}$};
  \node at (145:1)   {$\KK_{\leq -1}$};
  \node at (-35:1)   {$\KK_{> -1}$};
  \node at (-55:1)   {$\KK_{> 0}$};

    \end{tikzpicture}

%% file: diag_heart_fan_nilrep.tex
\begin{tikzpicture}[scale=3.5]
  \fill[black!10] (0,0) rectangle (1,1);
  \fill[black!15] (0,0) -- (0,1) -- (-1,1) -- cycle;
  \fill[black!25] (0,0) rectangle (-1,-1);
  \fill[black!17] (0,0) rectangle (1,-1);
  \draw[thick] (-1,0) -- (1,0);
  \draw[thick] (0,-1) -- (0,1);
  \draw[thick] (0,0) -- (-1,1);

  \begin{scope}
    \clip (-1,0) rectangle (1,1);
    \foreach \X in {0,1,...,44} {
      \draw[thick] (0,0) -- (135+\X:3);
    }
  \end{scope}
  
  \draw[ultra thick, magenta] (-1,0) -- (0,0);
  \draw[ultra thick, magenta] (-1,1) -- (0,0);

  \node at  (45:0.7) {$\HH = \nilrep(\rotatebox[origin=c]{-90}{$\circlearrowright$} \, \bullet \rightarrow \bullet)$};
  \node at (225:0.7) {$\HH[1]$};
  \node at (-45:0.7) {$\KK_1$};
  \node at (115:0.7) {$\KK$};
  \node[magenta] at (131:1.21) {$\KK'$};
  \node[magenta] at (185:0.9)  {};
\end{tikzpicture}

%% file: diag_heart_fan_projective_line.tex
\begin{tikzpicture}[scale=3.5]
  \fill[black!10] (0,0) -- (0,1) -- (1,1) -- cycle;
  \fill[black!15] (0,0) -- (1,1) -- (1,0.5) -- cycle;
  \fill[black!20] (0,0) -- (1,0.5) -- (1,0.333) -- cycle;
  \fill[black!25] (0,0) -- (1,0.333) -- (1,0.25) -- cycle;
  \fill[black!15] (0,0) -- (0,1) -- (-1,) -- cycle;
  \fill[black!20] (0,0) -- (-1,1) -- (-1,0.5) -- cycle;
  \fill[black!25] (0,0) -- (-1,0.5) -- (-1,0.333) -- cycle;
  \fill[black!30] (0,0) -- (-1,0.333) -- (-1,0.25) -- cycle;
  \draw[ultra thick] (-1,0) -- (1,0);
  \draw[thick] (0,0) -- (0,1);
  \draw[thick] (0,0) -- (1,1);
  \draw[thick] (0,0) -- (1,0.5);
  \draw[thick] (0,0) -- (1,0.333);
  \draw[thick] (0,0) -- (1,0.25);
  \draw[thick] (0,0) -- (-1,1);  
  \draw[thick] (0,0) -- (-1,0.5);
  \draw[thick] (0,0) -- (-1,0.333);
  \draw[thick] (0,0) -- (-1,0.25);
  \node at (-0.9,-0.10) {$\HH[1]$};
  \node at (0.75,-0.08) {$\HH = \coh{\PP^1}$};
  \node at ( 65:0.8) {$\KK_0$};
  \node at (110:0.8) {$\KK_1$};
  \node at ( 35:0.8) {$\KK_{-1}$};
  \node at (145:0.8) {$\KK_2$};
  \filldraw (178.0:0.9) circle (0.005);
  \filldraw (176.0:0.9) circle (0.005);
  \filldraw (174.0:0.9) circle (0.005);
  \filldraw (172.0:0.9) circle (0.005);
  \filldraw (170.0:0.9) circle (0.005);
  \filldraw (168.0:0.9) circle (0.005);
  \filldraw (  2.0:0.9) circle (0.005);
  \filldraw (  4.0:0.9) circle (0.005);
  \filldraw (  6.0:0.9) circle (0.005);
  \filldraw (  8.0:0.9) circle (0.005);
  \filldraw ( 10.0:0.9) circle (0.005);
  \filldraw ( 12.0:0.9) circle (0.005);
\end{tikzpicture}

%% file: diag_heart_fan_elliptic_curve.tex
\begin{tikzpicture}[scale=3.5]
  \begin{scope}
    \clip (-1,0) rectangle (1,1);
    \foreach \X in {0,2,...,180} {
      \draw[thick] (0,0) -- (\X:3) -- cycle;
  }
  \end{scope}
  \draw[very thick] (-1,0) -- (1,0);
  \node at (-0.9,-0.10) {$\HH[1]$};
  \node at (0.75,-0.08) {$\HH = \coh{X}$};
\end{tikzpicture}

%% file: diag_3-Kronecker__stability_fan.tex
\begin{tikzpicture}[scale=3.5]
  \draw[ultra thick] (-1,0) -- (1,0);
  \draw[ultra thick] (0,-1) -- (0,1);
  \draw[ultra thick] (0,0) -- (-0.5,1);
  \draw[ultra thick] (0,0) -- (-0.666,1);
  \begin{scope}
    \clip (-1,0) rectangle (1,1);
    \draw[ultra thick] (0,0) -- (130:1.5);
    \draw[ultra thick] (0,0) -- (140:105);
    \foreach \X in {0,1,...,5} {
      \draw[thick] (0,0) -- (130+2*\X:3);
    }
  \end{scope}
  \draw[ultra thick] (0,0) -- (-1,0.666);
  \draw[ultra thick] (0,0) -- (-1,0.5);
  \filldraw (125.5:1.1) circle (0.005);
  \filldraw (127.0:1.1) circle (0.005);
  \filldraw (128.5:1.1) circle (0.005);  
  \filldraw (141.5:1.1) circle (0.005);  
  \filldraw (143.0:1.1) circle (0.005); 
  \filldraw (144.5:1.1) circle (0.005);
\end{tikzpicture}

%% file: diag_3-Kronecker__heart_fan.tex
\begin{tikzpicture}[scale=3.5]
  \fill[black!20] (-1,-1) rectangle (1,1);
  \draw[thick] (-1,0) -- (1,0);
  \draw[thick] (0,-1) -- (0,1);
  \draw[thick] (0,0) -- (-0.5,1);
  \draw[thick] (0,0) -- (-0.666,1);
  \begin{scope}
    \clip (-1,0) rectangle (1,1);
    \draw[ultra thick] (0,0) -- (130:1.5);
    \draw[ultra thick] (0,0) -- (140:105);
    \foreach \X in {0,1,...,10} {
      \draw[thick, magenta] (0,0) -- (130+\X:3);
    }
    \foreach \X in {0,1,...,5} {
      \draw[thick] (0,0) -- (130+2*\X:3);
    }
  \end{scope}
  \draw[thick] (0,0) -- (-1,0.666);
  \draw[thick] (0,0) -- (-1,0.5);
  \filldraw (125.5:1.1) circle (0.005);
  \filldraw (127.0:1.1) circle (0.005);
  \filldraw (128.5:1.1) circle (0.005);  
  \filldraw (141.5:1.1) circle (0.005);  
  \filldraw (143.0:1.1) circle (0.005); 
  \filldraw (144.5:1.1) circle (0.005);
\end{tikzpicture}

%% file: diag_3-Kronecker__g-fan.tex
\begin{tikzpicture}[scale=3.5]
  \fill[black!20] (0,0) rectangle (1,1);
  \fill[black!20] (0,0) -- (0,1) -- (-0.83,1) -- cycle;
  \fill[black!20] (0,0) -- (-1,0) -- (-1,0.83) -- cycle;  
  \fill[black!20] (0,0) rectangle (-1,-1);
  \fill[black!20] (0,0) rectangle (1,-1);  
  \draw[thick] (-1,0) -- (1,0);
  \draw[thick] (0,-1) -- (0,1);
  \draw[thick] (0,0) -- (-0.5,1);
  \draw[thick] (0,0) -- (-0.666,1);
  \begin{scope}
    \clip (-1,0) rectangle (1,1);
    \draw[ultra thick] (0,0) -- (130:1.5);
    \draw[ultra thick] (0,0) -- (140:105);
  \end{scope}
  \draw[thick] (0,0) -- (-1,0.666);
  \draw[thick] (0,0) -- (-1,0.5);
  \filldraw (125.5:1.1) circle (0.005);
  \filldraw (127.0:1.1) circle (0.005);
  \filldraw (128.5:1.1) circle (0.005);  
  \filldraw (141.5:1.1) circle (0.005);  
  \filldraw (143.0:1.1) circle (0.005); 
  \filldraw (144.5:1.1) circle (0.005);
\end{tikzpicture}
